\documentclass[10pt,twoside]{article}
\usepackage[T1]{fontenc}
\usepackage{textcomp}
\usepackage{mathrsfs}
\usepackage{amsmath}
\usepackage{amssymb}
\usepackage{fancyhdr}
\usepackage{latexsym}
\usepackage{bbding}
\usepackage{wasysym}
\usepackage[pagewise]{lineno}
\usepackage{multicol,graphics}
\usepackage{xcolor}
\usepackage{esint}

\RequirePackage{times}

\def\XXint#1#2#3{{\setbox0=\hbox{$#1{#2#3}{\int}$ }
		\vcenter{\hbox{$#2#3$ }}\kern-.6\wd0}}

\def \d { \mathbf{d} }
\def \w { \mathbf{w} }
\def \m { \mathbf{m} }

\def \u { \mathbf{u} }
\def \z { \mathbf{z} }

\newcommand{\Dv}{{\rm div\,}}

\newcommand{\beq}{\begin{equation}}
\newcommand{\eeq}{\end{equation}}
\newcommand{\ben}{\begin{eqnarray}}
\newcommand{\een}{\end{eqnarray}}
\newcommand{\beno}{\begin{eqnarray*}}
\newcommand{\eeno}{\end{eqnarray*}}

\allowdisplaybreaks[4]

\newtheorem{Theorem}{Theorem}[section]

\newtheorem{Proposition}[Theorem]{Proposition}
\newtheorem{Lemma}[Theorem]{Lemma}
\newtheorem{Corollary}[Theorem]{Corollary}
\newtheorem{Remark}[Theorem]{Remark}

\newenvironment{proof}[1][Proof]{\noindent\textbf{#1.} }{\hfill $\Box$}

\numberwithin{equation}{section}

\begin{document}
	\title{Global existence of solutions to the 2D compressible nematic liquid crystal flows  with  large initial data and vacuum}

\author{Q{\sc iao} L{\sc iu}\thanks{School of Mathematics and Statistics, HNP-LAMA, Central South University, Changsha, Hunan 410083, P. R. China.
\\
        e-mail: \texttt{liuqiaomath@csu.edu.cn; liuqao2005@163.com.}}
    \,		
     ~ Q{\sc ingkui} M{\sc a}\thanks{School of Mathematics and Statistics, HNP-LAMA, Central South University, Changsha, Hunan 410083, P. R. China.
    \\
        e-mail: \texttt{qingkuima@163.com.}}
\,
 {\sc and } F{\sc uyi} X{\sc u}\thanks{School of Mathematics and Statistics, Shandong University of Technology, Zibo, Shandong 255049, P. R. China.
    \\
        e-mail:\texttt{zbxufuyi@163.com.}}}

 \date{}

\maketitle

\begin{abstract}
For periodic initial data that allow vacuum, we establish the global-in-time existence of solutions to the 2D compressible nematic liquid crystal flows, under the condition that the bulk viscosity coefficient is sufficiently large and a geometric angular constraint is imposed on the initial direction field. Our results remain valid for arbitrarily large initial energy. Furthermore, under appropriate additional assumptions imposed on the pressure, we also establish the uniqueness of the constructed solutions.
\\
\\


\end{abstract}
\vspace{4mm}

{\noindent \textbf{Keywords:} Compressible nematic liquid crystal flows;  global-in-time existence; large initial data; uniqueness; vacuum}
\medskip
\medskip

{\noindent \textit{2020 AMS Subject Classification}:  76A15, 76N10,  35Q35}


\section{Introduction and main results}

The motion of a general viscous compressible isentropic fluid on the 2D periodic domain  $\mathbb{T}^2$, coupled with the dynamic evolution of the nematic liquid crystal, is governed by the following isentropic compressible nematic liquid crystal flows:
\begin{align} \label{1.1}
\left\{
\begin{aligned}
&\rho_t+\Dv(\rho \u)=0,\\
&(\rho \u)_t+\Dv(\rho \u \otimes \u)+\nabla P=\mu \Delta \u+(\mu+\lambda) \nabla \Dv \u-\Dv\big( \nabla\d\odot\nabla\d-\frac{1}{2}|\nabla\d|^2\mathbb{I}_2\big),\\
&\d_t+\u\cdot\nabla\d=\Delta\d+|\nabla \d|^2 \d,
\end{aligned}
\right.
\end{align}
where $\rho(t,x)\ge0$, $\u(t,x)=(\u^1,\u^2)(t,x)$, $\d(t,x)=(\d^1,\d^2)(t,x)\in \mathbb{S}^1$ (the unit sphere in $\mathbb{R}^2$) and $P = P ( \rho)$ represent the density, the velocity, the macroscopic average of orientation, and the pressure, respectively. The shear viscosity coefficient $\mu$ and the bulk one $\lambda$ satisfy the following physical restrictions
\begin{equation}\label{1.2}
\mu>0, \quad\nu=2\mu+\lambda> 0.
\end{equation}
$\mathbb{I}_2$ is the $2 \times 2$ identity matrix. The term $\nabla \mathbf{d} \odot \nabla \mathbf{d}$ denotes the $2 \times 2$ matrix whose $(j, i)$-th entry is given by $\partial_{j} \mathbf{d} \cdot \partial_{i} \mathbf{d}$, for $1 \leq i, j \leq 2$, i.e., $\nabla \mathbf{d} \odot \nabla \mathbf{d} = (\nabla \mathbf{d})^\top \nabla \mathbf{d}$, where $(\nabla \mathbf{d})^\top$ denotes the transpose of the $2 \times 2$ matrix $\nabla \mathbf{d}$.
Here and in what follows, we focus only on the isentropic flows, i.e.,  the pressure $P = P(\rho)$ is a nonnegative $C^1$-function of the density and
satisfies $P'(\rho)\geq 0$ for all $\rho\ge0$, and is given by
\begin{equation*}
  P(\rho)=a\rho^{\gamma},
\end{equation*}
where $a > 0$ and $\gamma \geq 1$ are some physical parameters. We also assume that \eqref{1.1} is subject to the following
 initial condition
\begin{equation}\label{1.3}
 (\rho, \rho\u,  \d)(0,x)=(\rho_0, \rho_0\u_0, \d_0)(x), \quad x \in \mathbb{T}^2,
\end{equation}
and $|\d_0|=1$. For simplicity, the Lebesgue measure of $\mathbb{T}^2$ is normalized to unity, i.e.,  $|\mathbb{T}^2|=1$.    
It is obvious that the total mass and momentum of smooth solutions to \eqref{1.1}
are conserved through evolution, that is, for all $t \geq 0$,
\begin{equation}\label{1.5}
\begin{aligned}
&\int_{\mathbb{T}^2}\rho (t,x)\;\mathrm{d}x=\int_{\mathbb{T}^2}\rho_0(x)\;\mathrm{d}x,
\quad\text{and}\quad\int_{\mathbb{T}^2}(\rho \u)(t,x)\;\mathrm{d}x=
\int_{\mathbb{T}^2}\rho_0(x)\u_0(x)\;\mathrm{d}x.
\end{aligned}
\end{equation}
We assume that the initial total mass is equal to one.
By the Galilean invariance of system \eqref{1.1},
we may also assume, without loss of generality, that
the initial total momentum vanishes. Thus,
\begin{align}\label{1.6}
&\int_{\mathbb{T}^2}\rho_0(x)\,\mathrm{d}x=1
\quad\text{and}\quad
\int_{\mathbb{T}^2}\rho_0(x)\u_0(x)\,\mathrm{d}x=0.
\end{align}
\smallbreak

Liquid crystals constitute a state of matter that is intermediate between solids and liquids. The three principal types of liquid crystals are nematic, cholesteric, and smectic. Nematic liquid crystals exhibit long-range orientational order: their rod-like molecules tend to align along a common preferred direction while retaining liquid-like mobility. The continuum theory of liquid crystals was pioneered by Ericksen \cite{Ericksen} and Leslie \cite{Leslie}. Based on the original equations for nematic liquid crystals, Lin \cite{Lin} derived a simplified Ericksen-Leslie system for incompressible liquid crystal flows by choosing specific material constants in the bulk energy density of the Oseen-Frank theory. 
The physical relevance and mathematical challenges of the Ericksen-Leslie model have motivated extensive research, see, e.g., \cite{Ding,Ericksen,LL, LL1, LW1,LLW,Leslie} and the references therein. We briefly review results to the system \eqref{1.1} in multiple spatial dimensions. 
\smallbreak

The analysis of global weak solutions to system \eqref{1.1} builds on the foundational theories developed by Lions \cite{PLL} and by Feireisl \cite{EFA2} for the compressible Navier-Stokes equations.
Wang and Yu \cite{WangD} proved the global existence and large-time behavior of weak solutions to a
Ginzburg--Landau approximation of system \eqref{1.1}, using a three-level approximation scheme. In \cite{JJ}, Jiang et al. 
replaced the Ginzburg-Landau penalty function from \cite{WangD} with the supercritical nonlinearity $|\nabla \d|^2\d$,
and established the global existence of weak solutions to the corresponding system
with the homogeneous Neumann boundary condition $\frac{\partial \d}{\partial n}=0$ on $\partial\Omega$, where $\Omega\subset\mathbb{R}^N$ is a bounded domain with $N=2,3$. Later, based on new analytic techniques, Jiang et al.
\cite{JJ2} established the global existence of weak solutions to system \eqref{1.1} in two cases: the 2D initial-boundary value problem with small initial energy and the 2D Cauchy problem under a geometric angle condition on the initial direction field. However, these weak solutions have limited regularity, which makes it difficult to analyze their qualitative properties. Moreover, this broad class includes nonphysical solutions, and the regularity and uniqueness of weak solutions to system \eqref{1.1} remain open even in two dimensions.

\smallbreak

We next review the results on strong and smooth solutions to the system \eqref{1.1}.
Ding et al.
\cite{DCDS}
and Huang et al.
\cite{JDE} established the local existence of strong solutions with nonnegative initial density in one and three dimensions, respectively. Wang  \cite{JMFM} proved the global well-posedness of strong solutions to the two-dimensional system \eqref{1.1} with far-field vacuum. Hu and Wu \cite{9} studied the global well-posedness of strong solutions to the three-dimensional system under a smallness condition on the initial data. Li et al. \cite{LXZ} extended the local Cauchy result of \cite{JDE} to a global one for smooth initial data with small energy and possibly vanishing density. Dai \cite{Dai} improved this result by establishing global classical solutions when the adiabatic exponent is sufficiently close to $1$ and the far-field density $\rho_\infty$ is either zero or sufficiently close to zero. Xu and Zhong \cite{small} proved the global existence and uniqueness of strong solutions to the three-dimensional Cauchy problem for the non-isothermal compressible nematic liquid crystal system under a scaling-invariant smallness condition. More recently, Zhong and Zhou \cite{DeG} established the global well-posedness of smooth solutions with large initial data and vacuum in $\mathbb{R}^2$, under a geometric angle condition on the initial director field, with constant shear viscosity $\mu$ and density-dependent bulk viscosity $\lambda(\rho)=\rho^\beta$, where $\beta>4/3$. They later extended this result to general two-dimensional bounded domains in \cite{Bound}, subject to the Navier-slip boundary condition for the velocity and the homogeneous Neumann boundary condition for the director field.

\smallbreak

Despite this progress, the smooth and strong solutions reviewed above generally require the initial density to be regular and continuous. This requirement excludes discontinuous density profiles, such as piecewise constant densities, which arise naturally in models of mixtures of two fluids.
Many recent works have been dedicated to studying solutions with discontinuous density
(so-called "shock-data") for models of viscous compressible fluids. Even though the situation is by now quite well understood for the one dimensional case (see e.g. \cite{Hoff}), the multi-dimensional case is far from
being completely elucidated.  Over the last few decades, many important advances have been made in this direction concerning compressible fluid models. For the compressible Navier-Stokes equations, we mention the papers \cite{D21,D22} by Hoff, who
greatly contributed to constructing "intermediate solutions"   allowing for discontinuity of the density,
in between the weak solutions of Lions in \cite{PLL} and the classical ones of e.g. Nash in \cite{J.N}. In Hoff's works \cite{D21,D22}, the initial density is a small $L^\infty$ perturbation of a positive reference density $\widetilde{\rho}>0$. More precisely, $\|\rho_0-\widetilde{\rho}\|_{L^\infty(\mathbb{R}^N)}<\varepsilon$ with $0<\varepsilon<\widetilde{\rho}$, so that $\rho_0\geq\widetilde{\rho}-\varepsilon>0$ almost everywhere. Thus these results exclude vacuum and require small density oscillations. Later, Hoff and Santos \cite{D24} and Hoff \cite{D23} successfully removed the condition of small oscillations in the initial density \(\rho _{0}\) in \cite{D21,D22} which allows for the presence of vacuum. More recently, Danchin and Mucha \cite{R.D} established global solutions on $\mathbb{T}^N$, $N=2,3$, with bounded nonnegative initial density. Their results allow large initial data in two dimensions, provided that the bulk viscosity is sufficiently large, and small initial data in three dimensions. They also proved uniqueness for linear pressure laws.
For the three-dimensional compressible MHD equations with small initial energy and positive, bounded initial density, Suen and Hoff \cite{A.S} extended the result of \cite{D21} and established global weak solutions. Under stronger regularity assumptions on the initial velocity and magnetic field, Suen \cite{Suen} further obtained global solutions and proved uniqueness for linear pressure laws. More recently, following the arguments in \cite{R.D}, the authors of this paper \cite{L-M-X} established the global-in-time existence of solutions to the compressible MHD equations in the case of periodic domains $\mathbb{T}^2$  with large initial data and vacuum. For the  compressible nematic liquid crystal system \eqref{1.1}, Tan and Wu \cite{GuoC} extended Hoff’s results \cite{D21} by establishing global existence of weak solutions with small initial energy and positive, essentially bounded initial density.


\smallbreak

 To our knowledge, there are few results concerning the global well-posedness of "intermediate solutions" for the 
 system \eqref{1.1} under vacuum conditions. Motivated by  \cite{R.D}, our aim in this paper is to  investigate the global existence and uniqueness of solutions to  the system \eqref{1.1} with the  bounded, nonnegative initial density. 

\medskip

Before formulating our main results, let us list some notations and conventions that will be used throughout this paper.
\textbf{Notations}: (1)\ For $1\leq q\leq \infty$ and $k\in\mathbb{Z}^{+}$, the following notation is used for the Lebesgue and Sobolev spaces:
$$
\begin{aligned}
&L^{q}=L^{q}(\mathbb{T}^{2}),\quad
W^{k,q}(\mathbb{T}^{2})
=L^{q}(\mathbb{T}^{2})\cap D^{k,q}(\mathbb{T}^{2}),\quad
H^{k}(\mathbb{T}^{2})=W^{k,2}(\mathbb{T}^{2}),\\
&D^{k,q}(\mathbb{T}^{2})=\{f\in L_{loc}^{1}(\mathbb{T}^{2}): \|\nabla^{k}f\|_{L^q(\mathbb{T}^{2})}<\infty\}.
\end{aligned}
$$
For any Banach space $X$ and $1\leq p\leq\infty$, the notation
$\|\cdot\|_{L^{p}(0,T;X)}$ denotes the norm of $X$-valued functions in the space $L^{p}(0,T;X)$.
\medskip

\noindent
(2)\ For two matrices $A=(A_{ij})$ and $B=(B_{ij})$, their double contraction is defined by $
A:B=\sum_{i,j=1}^{2}A_{ij}B_{ij}$.
The Lebesgue measure of $\mathbb{T}^{2}$ is normalized so that
$|\mathbb{T}^{2}|=1$. The following shorthand notations are used:
$$
\int f\,\mathrm{d}x
:=\int_{\mathbb{T}^{2}}f\,\mathrm{d}x,\quad
\overline{f}
:=\frac{1}{|\mathbb{T}^{2}|}\int_{\mathbb{T}^{2}}f\,\mathrm{d}x
=\int_{\mathbb{T}^{2}}f\,\mathrm{d}x,\quad
\|\cdot\|_{p}:=\|\cdot\|_{L^{p}(\mathbb{T}^{2})}.
$$
\smallbreak

\noindent
(3)\ The material derivative of $f$ is denoted by both $\dot{f}$ and $\frac{D}{Dt}f$, that is, $\dot{f}=\frac{D}{Dt}f:=\partial_t f+\u\cdot\nabla f$.
Moreover, set
$
\widetilde{f}:=f-\overline{f}$,
where $\overline{f}$ denotes the spatial average of $f$ defined above. Let $\nabla^{\perp} := (-\partial_2, \partial_1)$ denote the 2D orthogonal gradient, satisfying $\nabla^\perp \cdot f = \partial_1 f_2 - \partial_2 f_1$ for $f = (f_1, f_2)$. Throughout the paper, $M$ denotes a generic positive constant whose value may vary from line to line. Its dependence on relevant parameters will be specified whenever necessary.
\medskip

Throughout this paper, we denote the initial total energy of system \eqref{1.1} by
\begin{equation}\label{initial}
C_0:=\int\Big(\frac{1}{2}\rho_0|\u_0|^{2}+ \frac{1}{2}|\nabla\d_0|^{2} +G(\rho_0)\Big)\;\mathrm{d}x,
\end{equation}
where $G(\cdot)$ is the potential energy density given by
\begin{equation}\label{1.8}
G(\rho ):= \rho \int_1^{\rho }{\frac{P(s)-P(1)}{s^{2}}\;\mathrm{d}s}.
\end{equation}
It is clear from \eqref{1.8} that\footnote{The equivalence between $\|G(\rho)\|_{1}$ and $\|\rho-1\|_2^2$ stems from the following facts: for $0 \leq \rho \leq \rho^{*}$, there exists a positive constant $c_1$, which depends only on $\gamma$ and $\rho^{*}$, such that
\begin{align*}
\frac{1}{c_1}(\rho - 1)^2 \leq (\rho^{\gamma} - 1^{\gamma})(\rho - 1) \leq c_1(\rho - 1)^2. 
\end{align*}
Direct calculation gives
\begin{align*}
\lim_{\rho \to 1} \frac{G(\rho)}{(\rho - 1)^2} =\frac{G''(1)}{2} =\frac{P'(1)}{2}= \frac{a\gamma}{2}.
\end{align*}
Thus there exists a positive constant $c_2$, which depends only on $\gamma$, $a$ and $\rho^{*}$, such that
$$
c_2^{-1}(\rho - 1)^2 \leq G(\rho) \leq c_2(\rho - 1)^2. 
$$
}
\begin{align*}
&\rho {G(\rho)}{''}={P'(\rho)},\;G(1)=G{'}(1)=0,\\
&G(\rho)=
\begin{cases}
a\big(\rho\log\rho+1-\rho\big),\quad~  \qquad\text{if }\, \gamma=1,\\
\frac{a}{\gamma-1}\big(\rho^\gamma-\rho\gamma+\gamma-1\big),\quad\ \, \text{if }\,\gamma>1,
\end{cases}
\\
&\|G(\rho)\|_1 \;\text{is essentially equivalent to}\;\|\rho-1\|_2 ^2.
\end{align*}
Moreover, the effective viscous flux $F$ is denoted by
\begin{equation*}
F:=\nu \Dv\u-P.   
\end{equation*}
We further denote
\begin{equation*}
  \widetilde{P}:=P-\int P\,\mathrm{d}x=P-\overline{P}\,\;\text{and}\,\;\widetilde{F}:=F-\int F\,\mathrm{d}x=F-\overline{F},
\end{equation*}
where $\overline{P}$ and $\overline{F}$ stand for
the average of $P$ and $F$, respectively. Owing to the definition of  $F$ and $\nu\int \Dv\u\,\mathrm{d}x=0$, it holds
\begin{align}\label{2.1}
  \widetilde{F}=\nu\Dv \u-\widetilde{P}.
\end{align}

For the initial data $(\rho_0, \mathbf{u}_0,\d_0)$, assume that there exist two positive constants $\rho^*$ ($\rho^*=2e^{\frac{\gamma-1}{a\gamma} C_{0}}\|\rho _{0}\|_{\infty } $ determined in Proposition \ref{2-3}) and $N_0$ (independent of $\nu$ and not necessarily small) such that
\begin{equation}\label{1.10}
0\leq\rho_0\leq\rho^*\quad\text{a.e. in }\mathbb T^2,
\end{equation}
\begin{equation}\label{1.11*}
(\mathbf{u}_0,\nabla\d_0) \in H^1(\mathbb{T}^2),\;\mu\|\nabla\mathbb{P} \mathbf{u}_0\|_2^2+ \nu\|\Dv \mathbf{u}_0\|_2^2+\|\nabla^2\d_0\|_2^2  \leq N_0,
\end{equation}
where $\mathbb{P}$ $(\mathbb{P}:=\mathbb{I}+\nabla(-\Delta)^{-1}\Dv)$
denotes the $L^2$-projector onto the set of solenoidal vector fields.
Moreover, let the initial director $\d_0=({\d_0}_1, {\d_0}_2)$ satisfy the following geometric condition 
\begin{equation}\label{1.12}
{\d_0}_2\geq \epsilon_0 \,\;\text{for some}\,\;\epsilon_0>0.      
\end{equation}
\smallbreak

Now, we state our global existence results as follows.

\begin{Theorem}\label{th:main1}
Assume the initial data $(\rho_0, \mathbf{u}_0,\d_0)$ satisfy \eqref{1.6} and \eqref{1.10}--\eqref{1.12}. Then there exists a positive number $\nu_0$ depending only on $a$, $\mu$, $\gamma$, $C_0$, $N_0$, $\overline{\omega}$ and $\|\rho_0\|_\infty$
such that if $\nu\ge\nu_0$\footnote{$\nu_{0}=\max\big\{\mu, \frac{M_4^2}{\log^2\frac32}\big\}$, where $M_4$ is a positive constant depending only on $a$, $\mu$, $\gamma$, $C_0$, $N_0$, $\overline{\omega}$ and $\|\rho_0\|_\infty$. See \eqref{2.65} for further details.}, the problem \eqref{1.1}--\eqref{1.3} admits a global-in-time solution $(\rho, \u, \d)$ on $(0,\infty)\times\mathbb{T}^2$, satisfying 
\begin{equation}\label{1.15}
\left
\{ \begin{aligned}
  &  \rho \in L^{\infty }(0,\infty;L^\infty(\mathbb{T}^2))\cap C([0,\infty);L^{p}(\mathbb{T}^{2})),\;1<p<\infty ,\;(\sqrt{\rho }\u,\nabla\d)\in C([0,\infty);L^{2}(\mathbb{T}^{2})); \\
 & (\u,\nabla\d)\in L^{\infty }(0,\infty;H^{1}(\mathbb{T}^{2})),\;(\nabla ^{2}\mathbb{P}\u,\nabla F, {\nabla\d}_{t}, \nabla ^{3}\d, \sqrt{\rho }\dot{\u})\in L^{2}(0,\infty;L^{2}(\mathbb{T}^{2})); \\
 & (\sqrt{\rho t}\dot{\u}, \sqrt{t}\nabla^3\d, \sqrt{t}\nabla\d_t)\in L^{\infty} _{loc}(0,\infty;L^{2}(\mathbb{T}^{2})),\;(\sqrt{t}\nabla \dot{\u},\sqrt{t}\d_{tt},\sqrt{t}\nabla^2\d_t)\in L_{loc}^{2}(0,\infty;L^{2}(\mathbb{T}^{2}));\\
 &(\Dv\u, \operatorname{curl}\u)\in L_{loc}^{r}(0,\infty;L^{\infty}(\mathbb{T}^{2}))\;\text{with}\;r\in(1,2),
\end{aligned}
\right.
\end{equation}
where $\dot{\u}=\u_{t}+\u\cdot \nabla \u$  
 is the material derivative of $\u$.
\end{Theorem}
\begin{Remark}\label{remark1}
1. 
Compared to the results in \cite{LXZ, JGA, GuoC}, which rely on the assumption of a sufficiently small initial energy, our Theorem \ref{th:main1} remains valid for arbitrarily large initial energy, provided that the bulk viscosity coefficient is chosen to be sufficiently large. Combining techniques from \cite{R.D} for the density and  methods from \cite{DeG, Bound} for the director field $\d$, one can prove that Theorem \ref{th:main1} holds for the initial density to be discontinuous and permits vacuum.  This differs from the setting of \cite{DeG, Bound}, where the initial density satisfies $0\leq\rho_0\in W^{1,q}$ for $q>2$, the viscosity coefficients $\mu$ are constant and $\lambda = \rho^\beta$ with $\beta>\frac{4}{3}$.
\medskip



\medskip

2.\ We notice that the initial geometric condition \eqref{1.12} can be revised by the following condition: The initial direction field $\d_0 = (\d_{01},\d_{02}) : \mathbb{T}^2 \to \mathbb{S}^1$, assume that $\d_{02}$ satisfies the condition
\begin{equation}\label{G1}
\d_{02} \geq \underline{\d_{02}} \geq 0, 
\end{equation}
where $\underline{\d_{02}}$ is constant and is defined as
\begin{equation}\label{G2}
\begin{cases}
\underline{\d_{02}}>\sqrt{1 - \displaystyle\frac{\varepsilon}{A_1 A_2}} , &\,\text{if}\,\; A_1 A_2 \geq \varepsilon, \\
\underline{\d_{02}}=0, &\,\text{if}\,\; A_1 A_2 < \varepsilon.
\end{cases}
\end{equation}
Here $\varepsilon\in(0,1)$ is a sufficiently small constant, and   $A_1$ and $A_2$ are the best constants of the elliptic estimate
\begin{equation}\label{G3}
\| \nabla^2 \d \|_2^2 = A_1 \| \Delta \d \|_2^2,
\end{equation}
and the inequality
\begin{equation}\label{G4}
\| \nabla \d \|_4^4 \leq A_2\|\d-e_2\|_\infty^2 \| \nabla^2 \d \|_2^2\,\;\text{with}\,\;e_2=(0,1). 
\end{equation}
In fact, \eqref{1.12} directly implies the inequality \eqref{L-6} stated below (see Lemma \ref{2-1-4}), which further enables us to establish the boundedness of  $ \d$  in $L^2(0,T;\dot{H}^2(\mathbb{T}^2))$. As shown in our subsequent Lemma \ref{lem2.4}, under the assumption \eqref{G1}, an inequality analogous to \eqref{L-6}, namely the estimate \eqref{G5}, can still be derived. With this key bound at our disposal, the boundedness of $\d$  in $L^2(0,T;\dot{H}^2(\mathbb{T}^2))$ can be established in an identical manner. 

\end{Remark}

In the case of a linear pressure, our existence result is supplemented with the following uniqueness.

\begin{Theorem}\label{th:main2}
Assume that $P(\rho)=a\rho$ for some $a>0$ and the initial data $(\rho_0, \mathbf{u}_0,\d_0)$ satisfy \eqref{1.6} and \eqref{1.10}--\eqref{1.12}. Then, for any $T > 0$  
and for any viscosity coefficients $(\mu,\lambda)$ satisfying
$\mu>0$ and $2\mu+\lambda>0$, there exists at most one solution to system \eqref{1.1}--\eqref{1.3} supplemented with data $(\rho_0, \u_0, \d_0)$ on $[0, T ]\times\mathbb{T}^2$, with the regularity
given in Theorem \ref{th:main1} (restricted to interval $[0, T ]$).
\end{Theorem}
\begin{Remark}\label{remark4}
Due to the presence of vacuum and the fact that  the director field $\d$
brings the strong coupling term $\delta\u\cdot\nabla\d_2$, we had to introduce a coupled backward parabolic system, see \eqref{3.8} below for details.
Motivated by \cite{Uniqueness, R.D} on the compressible Navier-Stokes
equations, we add the coupling terms $-\nabla\d_2\cdot\Delta\m$ and $\nabla\d_2\cdot\m$ in the equation for $\w$,  which combine with the duality argument to control the difference $\sqrt{\rho_1}\delta\u$ and $\delta\d$ in $[0, T ]\times\mathbb{T}^2$. 
\end{Remark}

The following result illustrates the $C^{0, \alpha_t}$ of the case where
the initial density has nontrivial vacuum regions.  The detailed proofs can be similar to that of \cite{R.D,R.D2}, therefore we omit them here.
\begin{Corollary}\label{Cor1}
Let the assumptions of Theorem \ref{th:main1} be in force,
and let $(\rho,\u,\d)$ be a global solution of \eqref{1.1}
given by Theorem \ref{th:main1}.
Let $\psi$ be the flow of $\u$, given by
\begin{equation*}
\psi(t,y)
=
y+\int_0^t \u(\tau,\psi(\tau,y))\,\mathrm{d}\tau,
\qquad t\geq0,\quad y\in\mathbb{T}^2.
\end{equation*}
For the statements concerning pointwise vacuum sets,
we use the representative of the density transported
by the flow $\psi$.
Then, the following conclusions hold:

(1) Let $Y_0:=\rho_0^{-1}(\{0\})$, where
$\rho_0^{-1}(\{0\})$ denotes the preimage of $\{0\}$
under the initial density function $\rho_0$, namely, $\rho_0^{-1}(\{0\})
=
\{x\in\mathbb{T}^2:\rho_0(x)=0\}$.
Then
\begin{equation*}
\rho(t,\cdot)^{-1}(\{0\})
=
Y_t,
\qquad
Y_t:=\psi(t,Y_0).    
\end{equation*}
Furthermore, if $Y_0$ is an open set with Lipschitz
boundary, then $Y_t$ is an open set with
$C^{0,\alpha_t}$ regularity.

(2) If $\rho_0=|V_0|^{-1}\mathbf{1}_{V_0}$,
$0<|V_0|<1$,
and $V_t:=\psi(t,V_0)$, then
\begin{equation*}
\operatorname*{ess\,inf}_{x\in V_t}
\rho(t,x)>0,
\qquad\text{for all }t>0.    
\end{equation*}
Furthermore, if $V_0$ is an open Lipschitz set,
then $V_t$ has $C^{0,\alpha_t}$ regularity. Here, $\alpha_t>0$ is a continuously decreasing
function of $t$ such that $\alpha_0=1$.
\end{Corollary}
\medskip

Let us comment on the proofs of our results. In the case of
the existence of solutions, we first derive  the basic energy estimate (see \eqref{1.11}) based on the geometric condition on the initial director $\d_0$ in \eqref{1.12}. Building on this estimate, we attempt to derive 
the $H^1(\mathbb{T}^2)$-norm of $(\u,\nabla\d)$. To do it, the main difficulties lie in the following two aspects:
\smallbreak

The first aspect concerns the crucial techniques developed by  \cite{R.D2, R.D, B.D} cannot be adopted to the present setting due to the strong coupling between the velocity and director fields. Testing the momentum equation with $\u_t$ leads, through the director equation, to terms involving $\u\cdot\nabla^2\d$ that cannot be directly controlled by the available estimates. We therefore test the momentum equation with the material derivative $\dot{\u}=\u_t+\u\cdot\nabla\u$ and use the identity(see \eqref{D1} and \eqref{D2})
\begin{equation*}
\nabla\d_t+\u\cdot\nabla^2\d
=
\nabla\Delta\d+\nabla(|\nabla\d|^2\d)
-(\nabla\u)\cdot\nabla\d.    
\end{equation*}
This allows the problematic transport terms to cancel after integration by parts. The remaining highest-order terms are absorbed by the $\|\nabla^3\d\|_2^2$ dissipation, allowing us to close the coupled energy estimates. The second aspect concerns the possible presence of vacuum and the lack of spatial regularity of the density. Vacuum causes the momentum equation $\eqref{1.1}_2$ to degenerate, while the lack of density regularity prevents direct control of $\nabla P$. In estimating the pressure contribution $\int \nabla P\cdot\dot{\u}\,\mathrm{d}x$ in \eqref{2.10}, we encounter the term
$\frac{2}{\nu}\int
P(\rho)\bigl(\rho P'(\rho)-P(\rho)\bigr)
\Dv\u\,\mathrm{d}x$,
which requires special treatment.
To this end, we introduce the auxiliary function
(see \eqref{2.13})
\begin{equation*}
k(\rho)
=
P^2(\rho)
-
\rho\int_1^\rho \frac{P^2(s)}{s^2}\,\mathrm{d}s, 
\end{equation*}
which satisfies $\partial_t k(\rho)+\Dv\bigl(k(\rho)\u\bigr)
=
-2P(\rho)\bigl(\rho P'(\rho)-P(\rho)\bigr)\Dv\u$
(see \eqref{2.14}).
This allows us to incorporate the troublesome term into the energy estimate (see \eqref{2.14}). Based on these, we give the $L^\infty(0,T;H^1(\mathbb{T}^2))$-norm of $(\u,\nabla\d)$. To bound $\rho^*$ in terms of the initial data, we follow the approach in \cite{B.D} and consider the quantity $\mathbf{E}:=\log \rho-\nu^{-1}(-\Delta )^{-1} \operatorname{div}(\rho\mathbf{u})$ (see Proposition \ref{2-3}).
Combining the momentum equation $\eqref{1.1}_2$ with the mass equation and utilizing the Coifman-Lions-Meyer-Semmes inequality (see \eqref{2.54}) 
, we can obtain the key upper bound on the density by the bootstrap argument.
\smallbreak

To establish uniqueness, it suffices to show that $\delta\rho$, $\delta\u$, $\delta\d$ and $\nabla\delta\d$ are all identically zero. To this end, following the approach in \cite{D21, China, DeG}, we derive time-weighted estimates for the velocity and director fields (see Proposition \ref{2-4}) to improve the regularity of the velocity field. These time-weighted estimates provide $L^2_{loc}(0,\infty;BMO(\mathbb{T}^2))$ control of $\nabla\u$, allowing us to replace the Lipschitz condition in the uniqueness proof.
Additional difficulty is that the director equation requires control of $\|\delta\u\|_2$, which cannot be bounded by $\|\sqrt{\rho_1}\delta\u\|_2$ in the presence of vacuum. Hence, we introduce the following coupled backward parabolic system:
\begin{equation*}
\left\{ \begin{aligned}
  & \rho_1 \w_t + \rho_1 \u_1\cdot \nabla \w + \mu \Delta \w + (\mu +\lambda )\nabla \Dv \w = -\rho_1\delta \u -\nabla\d_2\cdot\Delta\m + \nabla\d_2\cdot\m, \\
  & \m_t + \Delta \m = -\delta \d + \w\cdot\nabla\d_2, \\
  & \w|_{t=T} = \m|_{t=T} = 0.
\end{aligned} \right.
\end{equation*}
Combining the duality argument in \cite{R.D,Uniqueness} with estimates for solutions to the above system and system \eqref{3.3}, we control the differences between two solutions $(\rho_1,\u_1,\d_1)$ and $(\rho_2,\u_2,\d_2)$ and establish uniqueness.
\medskip

\textbf{Organization of the paper}. In Section 2, we introduce some useful inequalities and present some lemmas. In Section 3, we give the proof of Theorem \ref{th:main1}. Section 4 is devoted to the proof of Theorem \ref{th:main2}.

\section{Preliminaries}

In this section, we recall some known facts that will be used frequently later.
The first is the Gagliardo-Nirenberg inequality which can be found in \cite{O.A,A.N}.
\begin{Lemma}\label{2-1-1}
Let $1<m,r<\infty$. For every
$\z\in W^{1,m}(\mathbb{T}^2)\cap L^r(\mathbb{T}^2)$
with $\int\z\,\mathrm{d}x=0$,
it holds that
\begin{equation}\label{L-1}
\|\z\|_q
\leq c_0\|\nabla\z\|_m^\alpha
          \|\z\|_r^{1-\alpha},
\end{equation}
where $0\leq\alpha\leq1$ and
$\alpha
=
\left(\frac1r-\frac1q\right)
\left(\frac1r-\frac1m+\frac12\right)^{-1}$.
If $m<2$, then $q$ is between $r$ and
$\frac{2m}{2-m}$, that is,
$q\in[r,\frac{2m}{2-m}]$ if $r<\frac{2m}{2-m}$,
and $q\in[\frac{2m}{2-m},r]$ if
$r\geq\frac{2m}{2-m}$.
If $m=2$, then $q\in[r,\infty)$;
if $m>2$, then $q\in[r,\infty]$.
In the case $1<m<2$ and
$q=r=\frac{2m}{2-m}$, we take $\alpha=1$.
Here $c_0=c_0(m,r,q,\mathbb{T}^2)$ is independent
of $\z$.
Consequently, for every
$\z\in W^{1,m}(\mathbb{T}^2)\cap L^r(\mathbb{T}^2)$,
one has
\begin{equation}\label{L-2}
\|\z-\overline{\z}\|_q
\leq c_0\|\nabla\z\|_m^\alpha
          \|\z-\overline{\z}\|_r^{1-\alpha},
\end{equation}
where $\overline{\z}$ is the average of $\z$
over $\mathbb{T}^2$. Moreover, for every $\z\in H^1(\mathbb{T}^2)$
and $2<p<\infty$, it holds that
\begin{equation}\label{L-3}
\|\z\|_p
\leq \widehat c_0 p^{\frac12}
     \|\z\|_2^{\frac2p}
     \|\z\|_{H^1}^{1-\frac2p},
\end{equation}
where $\widehat c_0>0$ depends only on
$\mathbb{T}^2$ and is independent of $p$ and $\z$. 
\end{Lemma}

\begin{Lemma}(\cite{R.D2})\label{2-1-2}
Let $\z \in H^1(\mathbb{T}^2)$ and let $\rho$ be in $L^{2}(\mathbb{T}^2)$. Assume that
\begin{equation*}
  \int\rho \z \;\mathrm{d}x=0\;\text{and}\;N:=\int\rho \;\mathrm{d}x>0.
\end{equation*}
There exists a positive constant $c_1$ independent of $\rho$ such that
\begin{equation}\label{2.1.3*}
  \|\z\|_2 \leq \Big (1+\frac{c_1}{N}\|\rho-N\|_{2}\Big)\|\nabla \z\|_2.
\end{equation}
\end{Lemma}
Motivated by the rigidity estimates in \cite{Lei,DeG}, we establish
the following periodic rigidity estimate to obtain coercivity of
the director dissipation in the basic energy estimate.
\begin{Lemma}\label{2-1-4} Let $\epsilon_0 > 0$ be as in \eqref{1.12} and $E > 0$. There exists a positive constant $\overline{\omega} = \overline{\omega}(\epsilon_0, E) \in (0, 1)$ such that the following holds: If $\d = (\d_1, \d_2) : \mathbb{T}^2 \to \mathbb{S}^1$, $\nabla \d \in H^1(\mathbb{T}^2)$ with $\|\nabla \d\|_{2} \leq E$ and $\d_2 \geq \epsilon_0$, then
\[
\|\nabla \d\|_4^4 \leq (1 - \overline{\omega})\|\Delta \d\|_2^2.
\]
Consequently, the associated harmonic energy is coercive, i.e.,
\begin{equation}\label{L-6}
\|\Delta \d + |\nabla \d|^2 \d\|_2^2 \geq \frac{\overline{\omega}}{2}\big(\|\Delta \d\|_{L^2}^2 + \|\nabla \d\|_4^4\big)    
\end{equation}
\end{Lemma}
\begin{proof}
First of all, since $\d_2\geq\epsilon_0>0$, we can define
$\theta=\arctan \frac{\d_1}{\d_2}$. Then $\theta$ is a bounded
periodic function in $H^2(\mathbb{T}^2)$ and
$\d=(\sin\theta,\cos\theta)$. Direct computations give
\begin{equation*}
\begin{aligned}
\partial_i\d
&=(\cos\theta,-\sin\theta)\partial_i\theta\,\;\text{and}\,\;
\Delta\d
=(\cos\theta,-\sin\theta)\Delta\theta
  -\d|\nabla\theta|^2.
\end{aligned}
\end{equation*}
Consequently,
\begin{equation}\label{rigidity-gradient-identity}
\begin{aligned}
\|\nabla\theta\|_2
=&\|\nabla\d\|_2\leq E,\,\;
\|\nabla\d\|_4^4
=\|\nabla\theta\|_4^4,\\
\|\Delta\d\|_2^2
=&\|\Delta\theta\|_2^2+\|\nabla\theta\|_4^4,\,\;
\|\Delta\d+|\nabla\d|^2\d\|_2^2=\|\Delta\theta\|_2^2.
\end{aligned}
\end{equation}
 Applying \eqref{L-1} with
$\z=\nabla\theta$, $q=4$ and $m=r=2$, and using
periodic integration by parts, we obtain
\begin{equation}\label{rigidity-theta-interpolation}
\begin{aligned}
\|\nabla\theta\|_4^4
\leq c_0^4\|\nabla\theta\|_2^2
       \|\nabla^2\theta\|_2^2
&=c_0^4
       \|\nabla\theta\|_2^2
       \|\Delta\theta\|_2^2\leq c_0^4E^2\|\Delta\theta\|_2^2.
\end{aligned}
\end{equation}
Combining \eqref{rigidity-theta-interpolation} and
\eqref{rigidity-gradient-identity}, we obtain
\begin{equation*}
\|\nabla\d\|_4^4
\leq c_0^4E^2
\bigl(\|\Delta\d\|_2^2-\|\nabla\d\|_4^4\bigr).
\end{equation*}
Hence,
\begin{equation}\label{eq-R}
\|\nabla\d\|_4^4
\leq
\frac{c_0^4E^2}{1+c_0^4E^2}\|\Delta\d\|_2^2
=
(1-\overline{\omega})\|\Delta\d\|_2^2,
\end{equation}
where
$\overline{\omega}:=(1+c_0^4E^2)^{-1}\in(0,1)$. Applying \eqref{rigidity-gradient-identity} and \eqref{eq-R} yields that
\begin{equation}\label{rigidity-basic-estimates}
\begin{aligned}
\|\Delta\d+|\nabla\d|^2\d\|_2^2
&=\|\Delta\theta\|_2^2
=\|\Delta\d\|_2^2-\|\nabla\d\|_4^4,\\
\|\nabla\d\|_4^4
&\leq(1-\overline{\omega})\|\Delta\d\|_2^2.
\end{aligned}
\end{equation}
Finally, we infer from \eqref{rigidity-basic-estimates} that
\begin{equation*}
\begin{aligned}
\|\Delta\d+|\nabla\d|^2\d\|_2^2
&=\|\Delta\d\|_2^2-\|\nabla\d\|_4^4\\
&\geq\overline{\omega}\|\Delta\d\|_2^2\\
&\geq\frac{\overline{\omega}}{2}
\bigl(\|\Delta\d\|_2^2+\|\nabla\d\|_4^4\bigr).
\end{aligned}
\end{equation*}
This gives \eqref{L-6} and completes the proof.
\end{proof}

\begin{Lemma}\label{lem2.4}
For the initial direction field $\d_0 = (\d_{01},\d_{02}) : \mathbb{T}^2 \to \mathbb{S}^1$, assume that $\d_{02}$ satisfies \eqref{G1} and \eqref{G2}. 
Then there exists a positive constant
$ \tilde c=\tilde c(A_1,A_2) $,
where \(A_1\) and \(A_2\) are defined in \eqref{G3} and \eqref{G4}, respectively, such that
\begin{equation}\label{G5}
\| \nabla \d \|_4^4+\| \nabla^2 \d \|_2^2 \leq \tilde{c} \| \Delta \d + |\nabla \d|^2 \d \|_2^2.
\end{equation}
\end{Lemma}
\begin{proof}
Applying the maximum principle to the equation of $\d_2$ (i.e., the 2-th component of $\d$) together with the geometric condition \eqref{G1} yields that for any $t > 0$ (see \cite{JJ2, Bound2}),
$$
\inf_{x\in \mathbb{T}^2} \d_2(t,x) \geq \inf_{x\in \mathbb{T}^2} \d_{02}(x) \geq \underline{\d_{02}} \geq 0.
$$
Set \(e_2=(0,1)\) and $
\delta:=\|\d-e_2\|_{\infty}$.
Since \(|\d|=1\) and \(\d_2(t,x)\geq \underline \d_{02}\), one gets
\begin{equation*}
|\d-e_2|^2=|\d|^2+{e_2}^2-2\d\cdot e_2=2(1-\d_2),    
\end{equation*}
and hence
\begin{equation*}
\delta^2\leq 2(1-\underline \d_{02}). 
\end{equation*}
It follows from \eqref{G2} that, in both cases \(A_1A_2\geq\varepsilon\) and \(A_1A_2<\varepsilon\),
\begin{equation}\label{G6}
A_1A_2\delta^2<2\varepsilon.    
\end{equation}
Since $\d-e_2$ is periodic, integration by parts over $\mathbb T^2$ gives
\begin{equation}\label{G7}
\begin{aligned}
\|\nabla \d\|_4^4
&=-
\int
(\d-e_2)
\partial_i(
|\nabla \d|^2\partial_i \d
)\,\mathrm{d}x\\
&\leq
\bar{c}\|\d-e_2\|_\infty
\int_{\mathbb T^2}
|\nabla \d|^2|\nabla^2\d|\,\mathrm{d}x\\
&\leq
\bar{c}\delta\|\nabla \d\|_4^2
\|\nabla^2\d\|_2,
\end{aligned}    
\end{equation}
where $\bar{c}>0$ is a universal constant.
Consequently, applying Young's inequality to \eqref{G7}, there exists a positive constant $A_2$ depending only on $\bar{c}$ such that 
\begin{equation}\label{G8}
\|\nabla \d\|_4^4
\leq
A_2\delta^2\|\nabla^2\d\|_2^2.    
\end{equation}
Combining \eqref{G3}, \eqref{G8} with the equality $\|\Delta \d\|_2^2
=
\|\nabla \d\|_4^4
+
\|\Delta \d+|\nabla \d|^2\d\|_2^2$ yields that
\begin{equation*}
\begin{aligned}
\|\nabla^2\d\|_2^2
&=A_1\|\Delta \d\|_2^2\\
&=A_1(
\|\nabla \d\|_4^4
+
\|\Delta \d+|\nabla \d|^2\d\|_2^2
)\\
&\leq
A_1A_2\delta^2\|\nabla^2\d\|_2^2
+
A_1\|\Delta \d+|\nabla \d|^2\d\|_2^2.
\end{aligned}    
\end{equation*}
Hence, the preceding inequality, together with \eqref{G6} and a sufficiently small choice of 
$\varepsilon>0$
, say 
$\varepsilon<\frac14$
, yields
\begin{equation*}
 \begin{aligned}
  & \|\nabla^2\d\|_2^2\leq 2A_1\| \Delta \d + |\nabla \d|^2 \d \|_2^2 .
 \end{aligned}   
\end{equation*}
This together with \eqref{G8} yields the desired \eqref{G5}. 
\end{proof}

\section{Proof of Theorem \ref{th:main1}
}\label{section3}
\setcounter{section}{3}\setcounter{equation}{0}

In this section, we always assume that $(\rho, \u,\d)$ is the smooth solution of system \eqref{1.1}-\eqref{1.3}
on $[0, T]\times \mathbb{T}^2 $, and additionally assume that
\begin{equation}\label{Sup}
0 \leq \rho(t, x) \leq \rho^* := 2e^{\frac{\gamma-1}{a\gamma} C_{0}}\|\rho _{0}\|_{\infty } \quad \text{for all } (t, x) \in [0, T]\times \mathbb{T}^2 ,
\end{equation}
where $C_0$ is determined by \eqref{initial}.
Before proving Theorem \ref{th:main1}, we establish some necessary a priori estimates, see Lemmas \ref{Lemma1}--\ref{Lemma2} and Propositions \ref{2-2}--\ref{2-4}. Let us start with the elementary energy estimate of $(\rho,\u,\d)$.
\begin{Lemma}\label{Lemma1}
Let $(\rho, \u, \d)$ be a smooth solution of system \eqref{1.1}--\eqref{1.3} on $[0,T]\times \mathbb{T}^2 $. Then,
\begin{equation}\label{1.11}
\begin{aligned}
& \sup _{0 \leq t \leq T} \int\big(G(\rho)+\frac{1}{2} \rho|\u|^2+\frac{1}{2}|\nabla\d|^2\big) \;\mathrm{d} x\\
&+\int_0^T\int\big(\mu|\nabla \mathbb{P}\u|^2+\nu(\Dv \u)^2+\frac{\overline{\omega}}{2}(|\nabla^2\d|^2+|\nabla\d|^4)\big) \;\mathrm{d} x\,\mathrm{d} t \leq 2 C_0.
\end{aligned}
\end{equation}
\end{Lemma}

\begin{proof}
First, using the
definition of $G(\rho)$ and the mass equation, we get
\begin{equation*}
  G(\rho) _t+\Dv (G(\rho)\u)+(G'(\rho)\rho-G(\rho))\Dv\u=0.
\end{equation*}
The definition of $G(\rho)$ gives
\begin{equation*}
\rho G'(\rho)=\big(G(\rho)+P(\rho)-P(1)\big),
\end{equation*}
so we can derive
\begin{equation}\label{1-1.1}
 \frac{\mathrm{d}}{\mathrm{d}t}\int G(\rho)\;\mathrm{d}x+\int P(\rho)\,\Dv\u\;\mathrm{d}x=0.
\end{equation}
Multiplying $\eqref{1.1}_2$ by $\u$, integrating the resulting equation by parts over $\mathbb{T}^2 $, and using \eqref{1-1.1}, it follows 
\begin{equation}\label{1-1-2}
\begin{aligned}
&\frac{\mathrm{d}}{\mathrm{d} t}\int\big(G(\rho)+\frac{1}{2} \rho|\u|^2\big) \;\mathrm{d} x -\int\mu\big(\Delta\u-\nabla\Dv\u\big)\cdot\u \;\mathrm{d}x+\int\nu(\Dv \u)^2\;\mathrm{d}x=\int -\u\cdot(\nabla\d\Delta\d)\;\mathrm{d}x,
\end{aligned}
\end{equation}
where we have used the fact that 
\begin{equation*}
\begin{aligned}
\nabla\d\Delta\d=\Dv\big( \nabla\d\odot\nabla\d-\frac{1}{2}|\nabla\d|^2\mathbb{I}_2\big).
\end{aligned}
\end{equation*}
Then, integrating by parts, and then using $\mathbb{P}^{2}=\mathbb{P}$ and the definition of $\Delta \mathbb{P}\u$, we gather that
\begin{equation}\label{1-1-3}
\begin{aligned}
\int \mu\big(\Delta\u-\nabla \Dv \u\big)\cdot\u \;\mathrm{d}x= \int \mu \big( \mathbb{P}\Delta\u\big)\cdot\u \;\mathrm{d}x= \int \mu \big( \mathbb{P}^2\Delta\u\big)\cdot\u \;\mathrm{d}x
=-\int \mu  |\nabla\mathbb{P}\u|^2 \;\mathrm{d}x.
\end{aligned}
\end{equation}
Next, multiplying $\eqref{1.1}_3$ by $\Delta\d+|\nabla\d|^2\d$ and integrating over $\mathbb{T}^2$, using integration by parts and the
fact that $|\d|=1$, we obtain
\begin{equation*}
\begin{aligned}
\frac{\mathrm{d}}{\mathrm{d} t}\int\frac{1}{2} |\nabla\d|^2\,\mathrm{d}x+\int |\Delta\d+|\nabla\d|^2\d|^2\,\mathrm{d}x=\int\u\cdot(\nabla\d\Delta\d)\;\mathrm{d}x,
\end{aligned}
\end{equation*}
which along with \eqref{1-1-2},  \eqref{1-1-3} and by integrating the resulting equation over $[0,T]$, ensures that
\begin{equation*}
\begin{aligned}
& \sup _{0 \leq t \leq T} \int\big(G(\rho)+\frac{1}{2} \rho|\u|^2+\frac{1}{2}|\nabla\d|^2\big) \;\mathrm{d} x\\
&+\int_0^T\int\big(\mu|\nabla \mathbb{P}\u|^2+\nu(\Dv \u)^2+|\Delta\d+|\nabla\d|^2\d|^2\big) \;\mathrm{d} x\mathrm{d} t \leq 2 C_0 .
\end{aligned}
\end{equation*}
Combining this with \eqref{L-6}, we obtain \eqref{1.11}.
The proof of Lemma
is completed.
\end{proof}

\medskip 
We now derive the following  $L^2 (0,T;L^2(\mathbb{T}^2))$-norm of the pressure and the effective viscous flux.

\begin{Lemma}\label{Lemma2}
Let $(\rho,\u,\d)$ be a smooth solution to \eqref{1.1}--\eqref{1.3} on $[0,T]\times\mathbb{T}^2 $ satisfying \eqref{Sup}. Assume that $\nu\ge \mu$,  then the following holds: 
\begin{equation}\label{2.1.6}
 \begin{aligned}
\int_0^T\frac{1}{\nu}\big(\|\widetilde{P}\|_2^2+\|\widetilde{F}\|_2^2\big)\,\mathrm{d}t\leq M,
  \end{aligned}
\end{equation}
where $M$ denotes a generic positive
constant depending on $\mu$, $a$, $\gamma$, $C_0$, $\overline{\omega}$ and  $\|\rho_0\|_{\infty}$, independent of $T$ and $\nu$.
\end{Lemma}
\begin{proof}
First, recalling that  the definition of $G(\rho)$ by \eqref{1.8},
from this, there exist two positive constants $C_1$ and $C_2$ both depending only
on $a$, $\rho^*$ and $\gamma$
such that
\begin{equation*}
 C_1(\rho-1)^2\leq G(\rho)\quad \text{and}\quad G(\rho)\leq C_2 (\rho^\gamma-1)(\rho-1),   
\end{equation*}
as well as
\begin{equation}\label{2-1}
 C_1C_2^{-1}(\rho-1)^2\leq C_2^{-1}G(\rho) \leq(\rho^\gamma-1)(\rho-1).  
\end{equation}
Multiplying ${\eqref{1.1}}_2$ by $-\nabla(-\Delta)^{-1}(\rho-1)$ and  integrating the resultant over $\mathbb{T}^2$,  one obtains after integration by parts\footnote{Throughout this paper, we will use the Einstein summation over repeated indices convention.}( for any $\eta>0$)
\begin{align*}
&\quad \, \int \big( P - P(1) \big) (\rho - 1)\, \mathrm{d}x\\
&= \Big( -\int \rho \mathbf{u} \cdot \nabla(-\Delta)^{-1}(\rho - 1) \mathrm{d}x \Big)_t + \int \rho \mathbf{u} \cdot \nabla(-\Delta)^{-1}(\rho_t)\, \mathrm{d}x + (\mu + \lambda) \int (\rho - 1) \Dv\mathbf{u} \,\mathrm{d}x \\
&\quad - \mu \int \partial_j \u^i \partial_i \partial_j (-\Delta)^{-1}(\rho - 1)\, \mathrm{d}x + \int \rho \u^i \u^j \partial_i \partial_j (-\Delta)^{-1}(\rho - 1)\, \mathrm{d}x \\
&\quad+\int (\partial_j\d\cdot\partial_i\d)\partial_i \partial_j (-\Delta)^{-1}(\rho - 1)\, \mathrm{d}x-\int\frac{1}{2}(\partial_j\d\cdot\partial_j\d)\partial_i \partial_i (-\Delta)^{-1}(\rho - 1)\,\mathrm{d}x\\
&\leq \Big( -\int \rho \mathbf{u} \cdot \nabla(-\Delta)^{-1}(\rho - 1) \mathrm{d}x \Big)_t + M (\rho^*)^2\| \mathbf{u} \|_2^2 + \nu \| \Dv\mathbf{u} \|_{2} \| \rho - 1 \|_{2} \\
&\quad + M\big( \| \nabla \mathbf{u} \|_{2} \| \rho - 1 \|_{2} + \rho^*\| \mathbf{u} \|_4^2 \| \rho - 1 \|_{2}+\|\nabla\d\|_4^2\|\rho - 1\|_2\big) \\
&\leq \big( -\int \rho \mathbf{u} \cdot \nabla(-\Delta)^{-1}(\rho - 1)\, \mathrm{d}x \Big)_t + \eta \| \rho - 1 \|_2^2\\
&\quad+ M(\eta) \big( \| \mathbf{u} \|_2^2 + \nu^2 \| \Dv\mathbf{u} \|_2^2+\|\nabla\u\|_2^2+C_0\big(\|\u\|_4^2 +\|\nabla^2\d\|_2^2)\big),
\end{align*}
where we have used $\eqref{1.1}_1$, \eqref{1.11} and $\|\u\|_4^2\| \rho - 1 \|_{2} \leq M \|\u\|_4^2C_0^{1/2}$.
Notice that $\nu=2\mu+\lambda$ and 
\begin{equation*}
\mu \|\nabla \u\|_2^{2}+(\mu +\lambda )\|\Dv \u\|_2^{2}=\mu \|\nabla \mathbb{P}\u\|_2^{2}+\nu \|\Dv \u\|_2^{2},
\end{equation*}
thus we have for $\nu\ge \mu$,
\begin{equation}\label{2.19*}
\mu \|\nabla \mathbb{P}\u\|_2^{2}+\nu \|\Dv \u\|_2^{2}\ge  \mu \|\nabla \u\|_2^{2}.
\end{equation}
Next, we aim to estimate $\|\u\|_2$ and $\|\u\|_4$. By using \eqref{2.1.3*} (due to $\int \rho \mathbf{u}\, \mathrm{d}x =0$ by assumption \eqref{1.6}) and \eqref{L-3}, we deduce
\begin{equation}\label{2.30-1}
\begin{aligned}
\|\u\|_2+\|\u\|_4&\leq M  \big(\|\u\|_{2}+\|\u\|_2^{\frac{1}{2}}\|\u\|_{H^1}^{\frac{1}{2}}\big)\\
&\leq M (1+\|\widetilde{\rho}\|_{2})\|\nabla\u\|_{2},
\end{aligned}
\end{equation}
where $\widetilde{\rho}=\rho-1$.
Combining $P(\rho)=a\rho^\gamma$ and  \eqref{2-1}, it holds 
\begin{equation*}
\begin{aligned}
C_1C_2^{-1}\int (\rho-1)^2\,\mathrm{d}x\leq& \int (\rho^\gamma-1) (\rho-1)\,\mathrm{d}x= \frac{1}{a} \int \big( P - P(1) \big) (\rho - 1)\, \mathrm{d}x.
\end{aligned}
\end{equation*}
Hence, there exists a constant $C_3>0$ depending only
on $a$, $\rho^*$ and $\gamma$
such that
\begin{equation*}
\begin{aligned}
&C_3\int (\rho-1)^2\,\mathrm{d}x\\
\leq & \Big( -\int \rho \mathbf{u} \cdot \nabla(-\Delta)^{-1}(\rho - 1) \,\mathrm{d}x \Big)_t+ \eta \| \rho - 1 \|_2^2  + M \big( \| \nabla\mathbb{P}\mathbf{u} \|_2^2 + \nu^2 \| \Dv\mathbf{u} \|_2^2+\|\nabla^2\d\|_2^2\big).
\end{aligned}
\end{equation*}
Integrating the above inequality over $(0, T)$, choosing $\eta$ sufficiently small, and applying \eqref{1.11} together with  Poincar\'{e}'s inequality,
we obtain
\begin{equation}\label{2.1.7}
\begin{aligned}
&\int_0^T \|\rho-1\|_2^2\,\mathrm{d}t\\
\leq& M\sup_{0\leq t\leq T}\big(\|\sqrt{\rho}\u\|_2\|\nabla(-\Delta)^{-1}(\rho - 1)\|_{2}\big)+M\int_0^T\big( \| \nabla\mathbb{P}\mathbf{u} \|_2^2 + \nu^2 \| \Dv\mathbf{u} \|_2^2+\|\nabla^2\d\|_2^2\big)\,\mathrm{d}t\\
=& M\sup_{0\leq t\leq T}\big(\|\sqrt{\rho}\u\|_2\|\rho-1\|_{\dot{H}^{-1}}\big)+M\int_0^T\big( \| \nabla\mathbb{P}\mathbf{u} \|_2^2 + \nu^2 \| \Dv\mathbf{u} \|_2^2+\|\nabla^2\d\|_2^2\big)\,\mathrm{d}t\\
\leq& MC_0\big(\sup_{0\leq t\leq T}\|\rho-1\|_2+\nu\big)\leq M\nu.
\end{aligned}
\end{equation}
Moreover, it follows from the definition of $\widetilde{P}=P-\overline{P}$, the mean value theorem and $\int \rho_0\,\mathrm{d}x=1$ in \eqref{1.5} that
\begin{equation}\label{4.1}
 \|\widetilde P\|_2^2
=
\frac12
\int_{\mathbb{T}^2}\int_{\mathbb{T}^2}
|P(\rho(x))-P(\rho(y))|^2
\,\mathrm{d}x\,\mathrm{d}y\leq{a^2\gamma^2(\rho^\ast)^{2(\gamma-1)}}
\|\rho-1\|_2^2.  
\end{equation}
Combining \eqref{2.1.7}, \eqref{4.1} and \eqref{2.1}, we can derive from \eqref{1.11} that
\begin{equation*}
 \int_0^T\frac{1}{\nu}\big(\|\widetilde{P}\|_2^2+\|\widetilde{F}\|_2^2\big)\,\mathrm{d}t\leq M\int_0^T \big(\frac{1}{\nu}\|\rho-1\|_2^2+\nu\|\Dv\u\|_2^2\big)\,\mathrm{d}t\leq M, 
\end{equation*}
which leads to \eqref{2.1.6} and thus completes the proof of Lemma \ref{Lemma2}.
\end{proof}
\medskip

The following proposition concerns the time-independent estimate for the $\dot{H}^1(\mathbb{T}^2)$ norm of $(\u,\nabla\d)$ in terms of the initial data.
\begin{Proposition}\label{2-2}
Let $(\rho,\u,\d)$ be a smooth solution to \eqref{1.1}--\eqref{1.3} on $[0,T] \times \mathbb{T}^2 $ satisfying \eqref{Sup}. Assume that $\nu\ge \mu$,   it holds
\begin{equation}\label{2.5}
 \begin{aligned}
\sup_{0\leq t\leq T}\Phi(t)+\int_0^{T}\Theta(t)\;\mathrm{d}t
 \le M,
\end{aligned}
\end{equation}
where
\begin{equation*}
\begin{aligned}
& \Phi(t)=1+\mu\|\nabla \mathbb{P}\u\|_2^{2}+\frac{1}{\nu }\|\widetilde{F}\|_2^{2}+\frac{1}{\nu }\|\widetilde{P}\|_2^{2}
+\|\nabla^2 \d\|_2^{2}, \\
&\Theta(t)=\|\sqrt{\rho} \dot{\u}\|_2^{2}+\|\Delta \mathbb{P}\u\|_2^{2}+\|\nabla F\|_2^{2}+\|\nabla^3\d\|_2^{2}+\frac{1}{\nu^2}\|\widetilde{P}\|_2^2,   \end{aligned}    
\end{equation*}
and $M$ is a generic positive
constant depending on $\mu$, $a$, $\gamma$, $N_0$, 
$C_0$, $\overline{\omega}$  and $\|\rho_0\|_{\infty}$, independent of $T$ and $\nu$.
\end{Proposition}
\begin{proof}\ \emph{StepI: Formulate the energy functional.}
First, we rewrite $\eqref{1.1}_{2}$ as
\begin{equation*}
 \rho \dot{\u}+\nabla P=\mu \big(\Delta \u-\nabla\Dv\u\big)+\nu \nabla \Dv \u-\Dv\big( \nabla\d\odot\nabla\d-\frac{1}{2}|\nabla\d|^2\mathbb{I}_2\big).   
\end{equation*}
Multiplying the preceding equation by $\dot{\u}$ and integrating by parts over $\mathbb{T}^2$ gives
\begin{equation}\label{2.6}
\begin{aligned}
& \int\rho |\dot{\u}|^{2}\;\mathrm{d}x+\frac{\mu}{2}\frac{\mathrm{d}}{\mathrm{d}t}\int\big(|\nabla \u|^{2}-{(\Dv \u)}^{2} \big)\;\mathrm{d}x+\int\nabla P\cdot \dot{\u}\;\mathrm{d}x \\
=&\int\mu(\Delta\u-\nabla\Dv\u)\cdot (\u\cdot \nabla \u)\;\mathrm{d}x
+\int \nu\nabla\Dv\u\cdot\dot{\u}\;\mathrm{d}x\\
&+\frac{\mathrm{d}}{\mathrm{d}t}\int\big((\nabla\d\odot\nabla\d):\nabla \u-\frac{1}{2}(|\nabla\d|^2\mathbb{I}_2):\nabla\u\big)\;\mathrm{d}x \\
 &-\int(\nabla\d\odot\nabla\d)_t:\nabla \u\;\mathrm{d}x-\int\Dv(\nabla\d\odot\nabla\d)\cdot(\u\cdot\nabla\u)\,\mathrm{d}x\\
 &+\frac{1}{2}\int (|\nabla\d|^2\mathbb{I}_2)_t:\nabla\u \,\mathrm{d}x+\int \Dv\big(\frac{1}{2}(|\nabla\d|^2\mathbb{I}_2)\big)\u\cdot\nabla\u\,\mathrm{d}x.
\end{aligned}
\end{equation}
To treat the pressure term in \eqref{2.6}, define $h(\rho)=P'(\rho)\rho-P(\rho)$. Then
\begin{equation}\label{2.7}
{P(\rho)}_{t}+\Dv(P(\rho)\u)+h(\rho)\,\Dv \u=0.
\end{equation}
Then, integrating by parts yields
\begin{equation*}
\begin{aligned}
  \int\nabla P\cdot \dot{\u}\;\mathrm{d}x
  &=-\frac{\mathrm{d}}{\mathrm{d}t}\int P\Dv \u \;\mathrm{d}x-\int{h(\rho)}{(\Dv \u)}^{2}\;\mathrm{d}x+\int P\u\cdot \nabla \Dv \u \;\mathrm{d}x-\int P\Dv(\u\cdot\nabla\u)\;\mathrm{d}x\\
  &=-\frac{\mathrm{d}}{\mathrm{d}t}\int P\Dv \u \;\mathrm{d}x-\int{h(\rho)}{(\Dv \u)}^{2}\;\mathrm{d}x-\int P\nabla\u:\nabla\u^{\top} \;\mathrm{d}x.
  \end{aligned}
\end{equation*}
Notice that from the definition of $F$, it holds
\begin{equation*}
\begin{aligned}
 {(\nu \Dv \u)}^{2}
 =P^2+F^2+2PF=-P^{2}+F^{2}+2\nu P\Dv \u ,
   \end{aligned}
\end{equation*}
from which,  and by integrating by parts once to avoid the appearance of some $\nabla P$ term, we arrive at
\begin{equation}\label{2.10}
\begin{aligned}
&\quad\, \int\nabla P\cdot \dot{\u}\;\mathrm{d}x\\
&=-\frac{\mathrm{d}}{\mathrm{d}t}\int P\Dv \u \;\mathrm{d}x+\frac{1}{\nu ^{2}}\int(P^{2}-F^2)h(\rho)\;\mathrm{d}x-\frac{2}{\nu }\int Ph(\rho)\,\Dv \u \;\mathrm{d}x-\int P\nabla\u:\nabla\u^{\top} \;\mathrm{d}x. \\
\end{aligned}
\end{equation}
Owing to the definition of $F$ and to \eqref{2.7}, one gets
\begin{equation*}
\begin{aligned}
  -\int P\,\mathrm{d}x+\nu\int \Dv\u\,\mathrm{d}x=\int F\,\mathrm{d}x\;\Rightarrow\;\overline{F}=-\overline{P}, \quad\text{and}\quad\frac{\mathrm{d}}{\mathrm{d}t}\overline{P} =-\int h(\rho)\,\Dv \u \;\mathrm{d}x,
  \end{aligned}
\end{equation*}
hence,
\begin{equation}\label{2.12}
\begin{aligned}
 \int (P^{2}-F^2)\,h(\rho)\;\mathrm{d}x 
=&\int (P+F)(P-F)\,h(\rho)\;\mathrm{d}x\\
  =&\nu ^{2}\int{(\Dv \u)}^{2}h(\rho)\;\mathrm{d}x-2\nu \int\widetilde{F}\,\Dv \u \,h(\rho) \;\mathrm{d}x-\nu \frac{\mathrm{d}}{\mathrm{d}t}{(\overline{P})}^{2}. \\
\end{aligned}
\end{equation}
Next, let the function $k(\rho)$ be given as follows:
\begin{equation}\label{2.13}
k(\rho )=P^{2}(\rho )-\rho\int_1^{\rho }{\frac{P^{2}(s)}{s^2}}\;\mathrm{d}s.
\end{equation}
In fact, $k(\rho)$ is a solution of $k(\rho)-\rho {k'(\rho)}
=2P^2(\rho)-2\rho P'(\rho) P(\rho)$. 
From the mass equation \eqref{1.1}$_1$ and \eqref{2.13}, one infers that
\begin{equation*}
\begin{aligned}
\partial _{t}k(\rho)+\Dv(k(\rho)\u)=(k(\rho)-\rho k'(\rho))\Dv\u,
\end{aligned}
\end{equation*}
and
\begin{equation*}
\begin{aligned}
k(\rho)-\rho {k'(\rho)}
=2P^2(\rho)-2\rho P'(\rho) P(\rho)
=-2P(\rho)(\rho P'(\rho)-P(\rho)).
\end{aligned}
\end{equation*}
We immediately deduce from the equalities above and the definition of $h(\rho)$ in \eqref{2.7} that 
\begin{equation}\label{2.14}
 -\frac{2}{\nu}\int P\,h(\rho)\,\Dv \u\;\mathrm{d}x
 =\frac{1}{\nu}\frac{\mathrm{d}}{\mathrm{d}t}\int{k(\rho)}\;\mathrm{d}x.
\end{equation}
Then, plugging \eqref{2.12} and \eqref{2.14} in \eqref{2.10}, we obtain
\begin{equation}\label{pressure}
\begin{aligned}
  \int\nabla P\cdot \dot{\u}\;\mathrm{d}x& =\frac{\mathrm{d}}{\mathrm{d}t}\int\Big(\frac{k(\rho)-{(\overline{P})}^{2}}{\nu }-P\Dv \u\Big)\;\mathrm{d}x-\frac{2}{\nu }\int\widetilde{F}\Dv \u \,h(\rho) \;\mathrm{d}x \\
 & \quad+\int{(\Dv \u)}^{2}h(\rho)\;\mathrm{d}x-\int P\nabla\u:\nabla\u^{\top} \;\mathrm{d}x.
\end{aligned}
\end{equation}

Let us turn to the analysis of $\int \nu \nabla\Dv\u\cdot\dot{\u}\,\mathrm{d}x$, direct calculations show that
\begin{equation*}
\begin{aligned}
\Dv\dot{\u}&=\frac{D}{Dt}\Dv\u+(\partial_1\u\cdot\nabla) \u^1+(\partial_2\u\cdot\nabla) \u^2\\
&=\frac{D}{Dt}\Dv\u+2\nabla\u^1\cdot\nabla^{\perp} \u^{2} + (\Dv \u)^2,
\end{aligned}   
\end{equation*} 
hence,
\begin{equation}\label{rot}
 \begin{aligned}
\int \nu \nabla\Dv\u\cdot\dot{\u}\,\mathrm{d}x=&-\nu\int\Dv\u\,\Dv\dot{\u}\,\mathrm{d}x\\
=&-\frac{\nu}{2}\frac{\mathrm{d}}{\mathrm{d}t}\int (\Dv\u)^2\,\mathrm{d}x-\frac{\nu}{2}\int(\Dv\u)^3\,\mathrm{d}x\\
&-2\int (\widetilde{F}+\widetilde{P})\nabla\u^1\cdot\nabla^{\perp} \u^{2}\,\mathrm{d}x, 
 \end{aligned}   
\end{equation}
where we have used the definition of $\widetilde{F}$ in \eqref{2.1} and integration by parts. Combining \eqref{2.6}, \eqref{pressure} and \eqref{rot} gives
\begin{equation}\label{2.16}
 \begin{aligned}
& \frac{\mathrm{d}}{\mathrm{d}t}\Phi _{1}(t)+\int\rho |\dot{\u}|^{2}\;\mathrm{d}x \\
=&\int\Big(\frac{2}{\nu }\widetilde{F}\Dv \u \,h(\rho)- {(\Dv \u)}^{2}h(\rho)+P\nabla\u:\nabla\u^{\top}\Big) \;\mathrm{d}x +\int\mu(\Delta\u-\nabla\Dv\u)\cdot (\u\cdot \nabla \u)\;\mathrm{d}x
\\
&-\frac{\nu}{2}\int(\Dv\u)^3\,\mathrm{d}x-2\int (\widetilde{F}+\widetilde{P})\nabla\u^1\cdot\nabla^{\perp} \u^{2}\,\mathrm{d}x -\int(\nabla\d\odot\nabla\d)_t:\nabla \u\;\mathrm{d}x\\
 &-\int\Dv(\nabla\d\odot\nabla\d)\cdot(\u\cdot\nabla\u)\,\mathrm{d}x+\frac{1}{2}\int (|\nabla\d|^2\mathbb{I}_2)_t:\nabla\u \,\mathrm{d}x+\int \Dv\big(\frac{1}{2}(|\nabla\d|^2\mathbb{I}_2)\big)\u\cdot\nabla\u\,\mathrm{d}x,
\end{aligned}
\end{equation}
where
\begin{equation*}
\begin{aligned}
  \Phi _{1}(t) =\int&\Big({\frac{\mu }{2}\big( |\nabla \u|^{2}-{(\Dv \u)}^{2} \big)+\frac{1}{2\nu }\big( {(\nu \Dv \u)}^{2}-2\nu \Dv \u\widetilde{P}+2\big( k(\rho)-{(\overline{P})}^{2} \big) \big)}
  \\
 & -\big((\nabla\d\odot\nabla\d):\nabla \u-\frac{1}{2}(|\nabla\d|^2\mathbb{I}_2):\nabla\u\big)\Big)\;\mathrm{d}x.
\end{aligned}
\end{equation*}
 Now, we need to estimate
$\int k(\rho)-(\overline{P})^2\;\mathrm{d}x$, noting that $|\mathbb{T}^2|=1$, one has
\begin{equation*}
\begin{aligned}
\int k(\rho)-(\overline{P})^2\;\mathrm{d}x
=&\int \widetilde{P}^2\;\mathrm{d}x+\int 2(P-\overline{P})\overline{P}\;\mathrm{d}x
-\int\rho \int_1^{\rho }{\frac{P^{2}(s)}{s^2}}\;\mathrm{d}s\;\mathrm{d}x\\
=&\int \widetilde{P}^2\;\mathrm{d}x+\int\rho \int_\rho^{1 }{\frac{P^{2}(s)}{s^2}}\;\mathrm{d}s\;\mathrm{d}x,\\
\end{aligned}
\end{equation*}
which together with the definitions of $\widetilde{F}$ and $\Delta \mathbb{P}\mathbf{u}$, and $\mathbb{P}^2=\mathbb{P}$, we easily see that
\begin{equation*}
\begin{aligned}
  \Phi _{1}(t)&= \frac{1}{2}\Big(\mu \|\nabla \mathbb{P}\u\|_2^{2}+\frac{1}{\nu }\|\widetilde{F}\|_2^2+\frac{1}{\nu }\|\widetilde{P}\|_2^{2}\Big)+\frac{1}{\nu }\int\Big( \rho\int_\rho ^{1}{\frac{P^{2}(s)}{s^{2}}}\;\mathrm{d}s\Big)\; \mathrm{d}x\\
  &\quad-\int\big((\nabla\d\odot\nabla\d):\nabla \u-\frac{1}{2}(|\nabla\d|^2\mathbb{I}_2):\nabla\u\big)\, \mathrm{d}x.
\end{aligned}
\end{equation*}
On the other hand, applying the operator $\nabla$ to $\eqref{1.1}_3$, one has
\begin{equation}\label{eq1}
 \nabla\d_t-\nabla\Delta\d=-\nabla(\u\cdot\nabla\d)+\nabla(|\nabla\d|^2\d).   
\end{equation}
Then, multiplying \eqref{eq1} by $-\nabla\Delta \d$ and integrating over $\mathbb{T}^2$ ensures that
\begin{align}\label{eq-H1-d}
\begin{split}
& \frac{1}{2}\frac{\mathrm{d}}{\mathrm{d}t}\int |\nabla^2 \d|^{2}\;\mathrm{d}x+\int|\nabla^3\d|^{2}\;\mathrm{d}x\\
= &\int \nabla(\u\cdot\nabla\d)\cdot\nabla\Delta\d \;\mathrm{d}x-\int \nabla (|\nabla\d|^2\d)\cdot\nabla \Delta\d\,\mathrm{d}x\\
=&\int (\nabla\u\cdot\nabla\d)\cdot\nabla\Delta\d \;\mathrm{d}x-\int (\partial_k\u^j\partial_j\partial_i\d\cdot\partial_i\partial_k\d-\frac{1}{2}\partial_j\u^j\partial_i\partial_k\d\cdot\partial_i\partial_k\d)\,\mathrm{d}x\\
&-\int \nabla (|\nabla\d|^2\d)\cdot\nabla \Delta\d\,\mathrm{d}x.
\end{split}
\end{align}
Combining \eqref{2.16} with $\kappa_0\times$\eqref{eq-H1-d}, where the positive constant $\kappa_0$ will be specified later, we arrive at   
\begin{equation}\label{2.24}
\begin{aligned}
& \frac{\mathrm{d}}{\mathrm{d}t}\big(\Phi _{1}(t)+\frac{\kappa_0}{2} \|\nabla^2\d\|_2^2\big)+\int\rho |\dot{\u}|^{2}\;\mathrm{d}x+\kappa_0\int|\nabla^3\d|^{2}\;\mathrm{d}x \\
=&\int\Big(\frac{2}{\nu }\widetilde{F}\Dv \u \,h(\rho)- {(\Dv \u)}^{2}h(\rho)+P\nabla\u:\nabla\u^{\top}\Big) \;\mathrm{d}x +\int\mu(\Delta\u-\nabla\Dv\u)\cdot (\u\cdot \nabla \u)\;\mathrm{d}x
\\
&-\frac{\nu}{2}\int(\Dv\u)^3\,\mathrm{d}x-2\int (\widetilde{F}+\widetilde{P})\nabla\u^1\cdot\nabla^{\perp} \u^{2}\,\mathrm{d}x-\int(\nabla\d\odot\nabla\d)_t:\nabla \u\;\mathrm{d}x\\
 &-\int\Dv(\nabla\d\odot\nabla\d)\cdot(\u\cdot\nabla\u)\,\mathrm{d}x +\frac{1}{2}\int (|\nabla\d|^2\mathbb{I}_2)_t:\nabla\u \,\mathrm{d}x+\int \Dv\big(\frac{1}{2}(|\nabla\d|^2\mathbb{I}_2)\big)\u\cdot\nabla\u\,\mathrm{d}x\\
 &+\kappa_0\int \big((\nabla\u\cdot\nabla\d)\cdot\nabla\Delta\d-\nabla (|\nabla\d|^2\d)\cdot\nabla \Delta\d\big) \;\mathrm{d}x\\
 &-\kappa_0\int (\partial_k\u^j\partial_j\partial_i\d\cdot\partial_i\partial_k\d-\frac{1}{2}\partial_j\u^j\partial_i\partial_k\d\cdot\partial_i\partial_k\d)\,\mathrm{d}x\\
 \triangleq &\sum\limits_{i=1}^{10}I_{i}.
\end{aligned}
\end{equation}
To get control of the right-hand side of \eqref{2.24}, let us rewrite the momentum equation in terms of the viscous effective flux $F=\nu\Dv \u-P$ and $\Dv\big( \nabla\d\odot\nabla\d-\frac{1}{2}|\nabla\d|^2\mathbb{I}_2\big)=\nabla\d\Delta\d$ as follows:
\begin{equation*}
  \mu (\Delta \u-\nabla \Dv \u)+\nabla F=\rho \dot{\u}+\nabla\d\Delta\d,
\end{equation*}
from which, \eqref{L-1} and \eqref{1.11}, we discover that
\begin{equation*}
\begin{aligned}
  \mu ^{2}\|\Delta \mathbb{P}\u\|_2^{2}+\|\nabla F\|_2^{2}\le &2\rho ^{*}\|\sqrt{\rho }\dot{\u}\|_2^{2}+2\|\nabla\d\Delta \d\|_2^{2}
\\
  \leq& 2\rho ^{*}\|\sqrt{\rho }\dot{\u}\|_2^{2}+2\|\nabla\d\|_\infty^2\|\Delta \d\|_2^{2}
  \\
  \leq& 2\rho ^{*}\|\sqrt{\rho }\dot{\u}\|_2^{2}+2c_0\|\nabla\d\|_2^2 \|\nabla^3 \d\|_2^{2}
 \\
  \leq& 2(\rho ^{*}+2C_0c_0)\big(\|\sqrt{\rho }\dot{\u}\|_2^{2}+ \|\nabla^3 \d\|_2^{2}\big),
\end{aligned}
\end{equation*}
where we have used $\|\nabla\d\|_2^2\leq 2C_0$ from the basic energy identity.
Thus there exists a constant $M_0\geq1$, depending only on $\mu$, $C_0$, $\rho^*$ and $c_0$, such that
\begin{equation}\label{2.25}
\begin{aligned}
  \frac{\min\{1,\kappa_0\}}{2M_0}\big(\|\Delta \mathbb{P}\u\|_2^{2}+\|\nabla F\|_2^{2}\big) \le &\frac{1}{2}\|\sqrt{\rho }\dot{\u}\|_2^{2}+\frac{\kappa_0}{2}\|\nabla^3\d\|_2^{2},
  \end{aligned}
\end{equation}
where $C_0$, $c_0$ and $\kappa_0$ are as before.
In what follows, let us turn to bound $I_i(i=1,2,\cdots,10)$ term by term. To it, we first decompose $\nabla \u$ as
\begin{equation}\label{2.27}
\begin{aligned}
   \nabla \u&=\nabla (\mathbb{P}\u-\nabla {(-\Delta )}^{-1}\Dv \u) \\
 & =\nabla \mathbb{P}\u-\frac{1}{\nu }\nabla ^{2}{(-\Delta )}^{-1}\widetilde{F}-\frac{1}{\nu }\nabla ^{2}{(-\Delta )}^{-1}\widetilde{P}. \\
\end{aligned}
\end{equation}
By virtue of \eqref{2.1}, \eqref{2.27}, and Poincar\'{e}'s inequality, for any sufficiently small constant $\eta > 0$ and $\nu \ge \mu$, we estimate $I_1$ as
\begin{equation*}
\begin{aligned}
|I_1|&\leq M(\rho^*)(\|\widetilde{F}\|_2\|\Dv\u\|_2+\|\Dv\u\|_2^2+\|\nabla\u\|_2^2)\\
&\leq M\big(\|\nabla F\|_2\|\Dv\u\|_2+\|\Dv\u\|_2^2+\|\nabla\u\|_2^2)\\
&\leq \eta\|\nabla F\|_2^2+M\big(\|\nabla\mathbb{P}\u\|_2^2+\frac{1}{\nu}\|\widetilde{F}\|_2^2+\frac{1}{\nu}\|\widetilde{P}\|_2^2).
\end{aligned}
\end{equation*}
In view of  integration by parts, \eqref{L-1}, \eqref{2.27} and Cauchy-Schwarz inequality, one gets
\begin{equation*}
\begin{aligned}
|I_2|  & \le M\int\big( |\nabla\u|^3+|\nabla\u|^2|\Dv\u|+ |\Dv\u|^3\,\mathrm{d}x\\
&\le M\big( \|\nabla\u\|_3^3 +\|\Dv\u\|_3^3\big)\\
&\leq M\big(\|\nabla\mathbb{P}\u\|_2^2\|\Delta\mathbb{P}\u\|_2+\frac{1}{\nu^3}(\|\widetilde{F}\|_2^2\|\nabla F\|_2+\|\widetilde{P}\|_2^2\|\widetilde{P}\|_\infty)\big)\\
&\leq \eta\big(\|\Delta \mathbb{P}\u\|_2^{2}+\|\nabla F\|_2^{2}\big)+M\big(1+\|\nabla\mathbb{P}\u\|_2^2+\frac{1}{\nu}\|\widetilde{F}\|_2^2\big)\big(\|\nabla\mathbb{P}\u\|_2^2+\frac{1}{\nu}\|\widetilde{F}\|_2^2+\frac{1}{\nu}\|\widetilde{P}\|_2^2\big),
\end{aligned}
\end{equation*}
where we have used $\|\widetilde{P}\|_\infty\leq M(\rho^*)$. To estimate $I_3$ and $I_4$, we shall use the idea presented in \cite{JMPA}. Observe that
\begin{equation*}
\nabla^\perp \cdot \nabla \u^1 = 0, \quad \Dv \nabla^\perp \u^2 = 0,    
\end{equation*}
then from Theorem II.1 in \cite{R.CP}, it is easy to see that
\begin{equation}\label{Hardy}
\|\nabla \u^1 \cdot \nabla^\perp \u^2\|_{\mathcal{H}^1} \leq M \|\nabla \u\|_2^2.    
\end{equation}
Based on the fact that $BMO$ is the dual space of $\mathcal{H}^1$, we obtain
\begin{equation*}
\begin{aligned}
|I_3|+|I_4|\leq &M\big(\nu\|\Dv\u\|_3^3+(\|\widetilde{F}\|_{BMO}\|\nabla \u^1 \cdot \nabla^\perp \u^2\|_{\mathcal{H}^1}+\|\widetilde{P}\|_\infty\|\nabla \u^1 \cdot \nabla^\perp \u^2\|_{1}\big)\\
\leq&M\frac{1}{\nu^2}(\|\widetilde{F}\|_2^2\|\nabla F\|_2+\|\widetilde{P}\|_2^2\|\widetilde{P}\|_\infty)+M\big(\|\nabla F\|_2+\|\widetilde{P}\|_\infty\big)\|\nabla\u\|_2^2\\
\leq & \eta \|\nabla F\|_2^2+M\big(1+\|\nabla\mathbb{P}\u\|_2^2+\frac{1}{\nu}\|\widetilde{F}\|_2^2+\frac{1}{\nu}\|\widetilde{P}\|_2^2\big)\big(\|\nabla\mathbb{P}\u\|_2^2+\frac{1}{\nu}\|\widetilde{F}\|_2^2+\frac{1}{\nu}\|\widetilde{P}\|_2^2\big).
\end{aligned}    
\end{equation*}
To handle $I_5$ and $I_6$, the following observation is crucial for handling the strong coupling between the director field $\d$ and the velocity field $\u$, which helps us to avoid the a priori estimates of $\u\cdot\nabla^2\d$. By means of \eqref{eq1} and repeated integration by parts, we derive
\begin{align}
&- \int (\nabla \d_t \odot \nabla \d) : \nabla \u \, \mathrm{d}x - \int (\nabla \d \odot \nabla \d_t) : \nabla \u \, \mathrm{d}x - \int \Dv(\nabla \d \odot \nabla \d) \cdot (\u \cdot \nabla) \u \, \mathrm{d}x \nonumber\\
=&  - \int \left[ (\Delta \nabla \d - \nabla \u \cdot \nabla \d + \nabla(|\nabla \d|^2\d)) \odot \nabla \d \right] : \nabla \u \, \mathrm{d}x 
+  \int \u^k \partial_k\partial_j \d \cdot \partial_i \d \partial_j \u^i \, \mathrm{d}x \label{D1}\\
&- \int \left[ \nabla \d \odot \left( \Delta \nabla \d - \nabla \u \cdot \nabla \d + \nabla(|\nabla \d|^2 \d) \right) \right] : \nabla \u \, \mathrm{d}x \nonumber\\
&+  \int \u^k \partial_j \d \cdot \partial_k \partial_i \d \partial_j \u^i \, \mathrm{d}x +  \int \partial_j \d \cdot \partial_i \d \partial_j (\u^k \partial_k \u^i) \, \mathrm{d}x \nonumber\\
=&- \int \left[ (\Delta \nabla \d - \nabla \u \cdot \nabla \d + \nabla(|\nabla \d|^2\d)) \odot \nabla \d \right] : \nabla \u \, \mathrm{d}x \nonumber\\
&- \int \left[ \nabla \d \odot \left( \Delta \nabla \d - \nabla \u \cdot \nabla \d + \nabla(|\nabla \d|^2 \d) \right) \right] : \nabla \u \, \mathrm{d}x\nonumber\\
&+ \int \big(\partial_j \d \cdot \partial_i \d \partial_j \u^k \partial_k \u^i-\partial_k \u^k \partial_j \d \cdot \partial_i \d \partial_j \u^i \big)\, \mathrm{d}x\nonumber,
\end{align}
which along with the fact $|\d|=1$ and Cauchy-Schwarz inequality, indicates that 
\begin{equation*}
\begin{aligned}
|I_5+I_6|
&\leq M\Big(\int |\nabla\u||\nabla\d||\Delta\nabla\d|\,\mathrm{d}x+\int |\nabla\d|^2|\nabla\u|^2\,\mathrm{d}x+\int |\nabla\u||\nabla\d|^4\,\mathrm{d}x\\
&\quad+\int |\nabla^2\d||\nabla\d|^2|\nabla\u|\,\mathrm{d}x+\int |\Dv\u||\nabla\d|^2|\nabla\u|\,\mathrm{d}x\Big)\\
&\leq M  \big(\|\nabla^3\d\|_2\|\nabla\d\|_\infty\|\nabla\u\|_2+\|\nabla\d\|_\infty^2\|\nabla\u\|_2^2+\|\nabla\u\|_2\|\nabla\d\|_\infty^2\|\nabla\d\|_4^2\\
&\quad+\|\nabla^2\d\|_4\|\nabla\d\|_8^2\|\nabla\u\|_2+\|\Dv\u\|_2\|\nabla\d\|_\infty^2\|\nabla\u\|_2\big).
\end{aligned}
\end{equation*}
By using \eqref{L-1} and \eqref{2.27}, we then ensure that
\begin{equation}\label{2.30}
\begin{aligned}
|I_5+I_6|
\leq & M  \big(\|\nabla^3\d\|_2^{\frac{3}{2}}\|\nabla\d\|_2^{\frac{1}{2}}\|\nabla\u\|_2+\|\nabla\d\|_2\|\nabla^3\d\|_2\|\nabla\u\|_2^2+\|\nabla\u\|_2\|\nabla\d\|_2^2\|\nabla^3\d\|_2\|\nabla^2\d\|_2\\
&\quad+\|\nabla\d\|_2^{\frac{1}{2}}\|\nabla^2\d\|_2^{2}\|\nabla^3\d\|_2^{\frac{1}{2}}\|\nabla\u\|_2+\|\Dv\u\|_2\|\nabla\d\|_2\|\nabla^3\d\|_2\|\nabla\u\|_2\big)\\
\leq& \eta\|\nabla^3\d\|_2^2+M\big(\|\nabla\u\|_2^2+\|\nabla^2\d\|_2^2\big)\big(\|\nabla\u\|_2^2+\|\Dv\u\|_2^2+\|\nabla^2\d\|_2^2\big)\\
\leq &\eta\|\nabla^3\d\|_2^2+M\big(1+\|\nabla\mathbb{P}\u\|_2^2+\frac{1}{\nu}\|\widetilde{F}\|_2^2+\frac{1}{\nu}\|\widetilde{P}\|_2^2+\|\nabla^2\d\|_2^2\big)\\
&\quad\times\big(\|\nabla\mathbb{P}\u\|_2^2+\frac{1}{\nu}\|\widetilde{F}\|_2^2+\frac{1}{\nu}\|\widetilde{P}\|_2^2+\|\nabla^2\d\|_2^2\big).
\end{aligned}    
\end{equation}
The terms $I_7$ and $I_8$ are treated similarly to $I_5$ and $I_6$:
\begin{align}
I_7+I_8 &= \int \partial_i \mathbf{u}^i \partial_j \mathbf{d} \cdot \partial_j\mathbf{d}_t \, \,\mathrm{d}x +  \int \partial_j\partial_i\d\cdot\partial_j\d\,\u^k\partial_k\u^i \,\mathrm{d}x \nonumber\\
&= \int \partial_i \mathbf{u}^i \partial_j \mathbf{d} \cdot \big( \Delta \partial_j \mathbf{d} - \partial_j \mathbf{u} \cdot \nabla \mathbf{d} + \partial_j\big( |\nabla \mathbf{d}|^2 \mathbf{d} \big) \big) \,\mathrm{d}x -\int \partial_i \mathbf{u}^i \partial_j \mathbf{d} \cdot \partial_{k}\partial_j \mathbf{d} \,\mathbf{u}^k \,\mathrm{d}x\nonumber \\
&\quad +  \int \partial_j\partial_i\d\cdot\partial_j\d\,\u^k\partial_k\u^i \,\mathrm{d}x \label{D2}\\
&= \int \partial_i \mathbf{u}^i \partial_j \mathbf{d} \cdot \big( \Delta \partial_j \mathbf{d} - \partial_j \mathbf{u} \cdot \nabla \mathbf{d} + \partial_j\big( |\nabla \mathbf{d}|^2 \mathbf{d} \big) \big) \,\mathrm{d}x\nonumber\\
&\quad -\frac{1}{2} \int\partial_j\d\cdot\partial_j\d\partial_i\u^k\partial_k\u^i \,\mathrm{d}x+\frac{1}{2}\int\partial_j\d\cdot\partial_j\d\partial_k\u^k\partial_i\u^i\,\mathrm{d}x\nonumber,
\end{align}  
from which, 
it is not difficult to derive that
\begin{equation}\label{2.31}
\begin{aligned}
|I_7+I_8|
\leq &M\big(\|\Dv \u\|_2\|\nabla \d\|_\infty\|\nabla^3\d\|_2 +\|\Dv\u\|_2\|\nabla\u\|_2\|\nabla\d\|_\infty^2+\|\Dv\u\|_2\|\nabla^2\d\|_4\|\nabla\d\|_8^2\\
&+\|\Dv\u\|_2\|\nabla\d\|_\infty^2\|\nabla\d\|_4^2+\big(\|\Dv\u\|_2+\|\nabla\u\|_2\big)\|\nabla\d\|_\infty^2\|\nabla\u\|_2\big)\\
\leq &M\big(\|\Dv \u\|_2\|\nabla \d\|_2^{\frac{1}{2}}\|\nabla^3\d\|_2^{\frac{3}{2}} +\|\Dv\u\|_2\|\nabla\u\|_2\|\nabla\d\|_2\|\nabla^3\d\|_2\\
&+\|\Dv\u\|_2\|\nabla\d\|_2^{\frac{1}{2}}\|\nabla^2\d\|_2^{2}\|\nabla^3\d\|_2^{\frac{1}{2}}\\
&+\|\Dv\u\|_2\|\nabla\d\|_2^2\|\nabla^3\d\|_2\|\nabla^2\d\|_2+\big(\|\Dv\u\|_2+\|\nabla\u\|_2\big)\|\nabla\d\|_2\|\nabla^3\d\|_2\|\nabla\u\|_2\\
\leq&\eta\|\nabla^3\d\|_2^2+M\big(\|\nabla\u\|_2^2+\|\Dv\u\|_2^2+\|\nabla^2\d\|_2^2\big)\big(\|\nabla\u\|_2^2+\|\Dv\u\|_2^2+\|\nabla^2\d\|_2^2\big)\\
\leq &\eta\|\nabla^3\d\|_2^2+M\big(1+\|\nabla\mathbb{P}\u\|_2^2+\frac{1}{\nu}\|\widetilde{F}\|_2^2+\frac{1}{\nu}\|\widetilde{P}\|_2^2+\|\nabla^2\d\|_2^2\big)\\
&\quad\times\big(\|\nabla\mathbb{P}\u\|_2^2+\frac{1}{\nu}\|\widetilde{F}\|_2^2+\frac{1}{\nu}\|\widetilde{P}\|_2^2+\|\nabla^2\d\|_2^2\big).
\end{aligned}
\end{equation}
Along the same lines as in \eqref{2.30} and \eqref{2.31},  the terms $I_9$ and $I_{10}$  satisfy
\begin{equation*}
\begin{aligned}
|I_9|+|I_{10}|
&\leq M  \big(\|\nabla\u\|_2\|\nabla\d\|_\infty\|\nabla^3\d\|_2+(\|\nabla\u\|_2+\|\Dv\u\|_2)\|\nabla^2\d\|_4^2\\
&\quad+\|\nabla\d\|_6^3\|\nabla^3\d\|_2+\|\nabla^3\d\|_2\|\nabla^2\d\|_4\|\nabla\d\|_4\big)\\
&\leq M  \big(\|\nabla\u\|_2\|\nabla\d\|_2^{\frac{1}{2}}\|\nabla^3\d\|_2^{\frac{3}{2}}+(\|\nabla\u\|_2+\|\Dv\u\|_2)\|\nabla^2\d\|_4^2\\
&\quad+\|\nabla\d\|_2\|\nabla^2\d\|_2^2\|\nabla^3\d\|_2+\|\nabla^3\d\|_2^{\frac{3}{2}}\|\nabla^2\d\|_2\|\nabla\d\|_2^{\frac{1}{2}}\big)\\
&\leq \eta\|\nabla^3\d\|_2^2+M\big(\|\nabla\u\|_2^4+\|\nabla^2\d\|_2^4+(\|\nabla\u\|_2^2+\|\Dv\u\|_2^2)\|\nabla^2\d\|_2^2\big)\\
& \leq\eta\|\nabla^3\d\|_2^2+M\big(\|\nabla\u\|_2^2+\|\Dv\u\|_2^2+\|\nabla^2\d\|_2^2\big)^2\\
& \leq \eta\|\nabla^3\d\|_2^2+M\big(1+\|\nabla\mathbb{P}\u\|_2^2+\frac{1}{\nu}\|\widetilde{F}\|_2^2+\frac{1}{\nu}\|\widetilde{P}\|_2^2+\|\nabla^2\d\|_2^2\big)\\
&\quad\times\big(\|\nabla\mathbb{P}\u\|_2^2+\frac{1}{\nu}\|\widetilde{F}\|_2^2+\frac{1}{\nu}\|\widetilde{P}\|_2^2+\|\nabla^2\d\|_2^2\big).
\end{aligned}
\end{equation*}
With $\kappa_0$ fixed as specified below (so that $\kappa_0\geq1$),
choose $\eta=1/(24M_0)$.
Inserting all estimates of $I_i$ $(i=1,\ldots,10)$
into \eqref{2.24} and using \eqref{2.25},
we obtain
\begin{equation}\label{2.24-1}
\begin{aligned}
& \frac{\mathrm{d}}{\mathrm{d}t}\big(\Phi _{1}(t)+\frac{\kappa_0}{2} \|\nabla^2\d\|_2^2\big)+\frac{1}{2}\|\sqrt{\rho}\dot{\u}\|_2^{2}+\frac{\kappa_0}{2}\|\nabla^3\d\|_2^{2} + \frac{1}{4M_0}\Big(\|\Delta \mathbb{P}\u\|_2^{2}+\|\nabla F\|_2^{2}\Big)\\
\leq &M\big(\|\nabla\mathbb{P}\u\|_2^2+\frac{1}{\nu}\|\widetilde{F}\|_2^2+\frac{1}{\nu}\|\widetilde{P}\|_2^2+\|\nabla^2\d\|_2^2\big)\big(1+\|\nabla\mathbb{P}\u\|_2^2+\frac{1}{\nu}\|\widetilde{F}\|_2^2+\frac{1}{\nu}\|\widetilde{P}\|_2^2+\|\nabla^2\d\|_2^2\big).
\end{aligned}
\end{equation}
Recalling that \eqref{1-1.1} and the effective viscous flux $\widetilde{F}=\nu\Dv \u-\widetilde{P}$, we arrive at
\begin{equation*}
\begin{aligned}
\frac{\mathrm{d}}{\mathrm{d}t}\int G(\rho)\;\mathrm{d}x+\frac{1}{\nu}\int|\widetilde{P}|^2 \;\mathrm{d}x
&=-\frac{1}{\nu}\int\widetilde{P}\widetilde{F} \;\mathrm{d}x \\
&\leq \frac{1}{2\nu}\|\widetilde{P}\|_2^2+\frac{1}{2\nu}\|\widetilde{F}\|_2^2,\\
\end{aligned}
\end{equation*}
where we have used $\int \overline{P}\Dv\u\;\mathrm{d}x=0$. The inequality above further implies
\begin{equation}\label{2.33}
\begin{aligned}
\frac{\mathrm{d}}{\mathrm{d}t}\Big(\int\frac{P^*}{\nu} G(\rho)\;\mathrm{d}x
+\frac{1}{\nu }P(1)( P^{*}-P(1) )\Big)+\frac{P^*}{2\nu^2}\int|\widetilde{P}|^2 \;\mathrm{d}x
\leq\frac{P^*}{2\nu^2}\|\widetilde{F}\|_2^2\leq  \frac{M}{\nu} \|\widetilde{F}\|_2^2(1+\frac{1}{\nu}\|\widetilde{F}\|_2^2),
\end{aligned}
\end{equation}
where $P^{*}:=a(\rho^*)^\gamma$. 
Combining \eqref{2.33} and \eqref{2.24-1}, it gives that
\begin{equation}\label{EQU-1}
\begin{aligned}
& \!\frac{\mathrm{d}}{\mathrm{d}t}\Phi_2(t)+\Theta_1(t)
\leq M\big(  
  \|\nabla\mathbb{P}\u\|_2^2+\frac{1}{\nu}\|\widetilde{F}\|_2^2+\frac{1}{\nu}\|\widetilde{P}\|_2^2+\|\nabla^2\d\|_2^2\big)\Phi(t),
\end{aligned}
\end{equation}
where 
\begin{equation*}
\begin{aligned}
\Phi _{2}(t)&= 1+\Phi_1(t)+\frac{\kappa_0}{2} \|\nabla^2\d\|_2^2+\frac{P^*}{\nu}\|G(\rho)\|_1+\frac{1}{\nu }P(1)( P^{*}-P(1) );
  \\
  \Theta_1(t) &=\frac{1}{2}\|\sqrt{\rho}\dot{\u}\|_2^{2}+\!\frac{\kappa_0}{2}\|\nabla^3\d\|_2^{2} + \frac{1}{4M_0}\Big(\|\Delta \mathbb{P}\u\|_2^{2}+\|\nabla F\|_2^{2}\Big)+\frac{P^*}{4\nu^2}\|\widetilde{P}\|_2^{2};\\
\Phi(t)&=1+\mu\|\nabla\mathbb{P}\u\|_2^2+\frac{1}{\nu}\|\widetilde{F}\|_2^2+\frac{1}{\nu}\|\widetilde{P}\|_2^2+\|\nabla^2\d\|_2^2.
  \end{aligned}
\end{equation*}

\emph{StepII}: \emph{Upgrading the functional $\Phi_2(t)$.} It is easy to see that
\begin{equation*}
\begin{aligned}
\rho \int_1^{\rho }{\frac{P^{2}(s)}{s^{2}}}\;\mathrm{d}s& \le \rho P(\rho )\int_1^{\rho }{\frac{P(s)}{s^{2}}}\;\mathrm{d}s\\
&=P(\rho )\big(G(\rho )+P(1)(\rho -1)\big) \\
&\le P^{*}\big(G(\rho )+P(1)\rho \big)-P(\rho )P(1).
\end{aligned}
\end{equation*}
Hence, owing to \eqref{1.6}, we deduce that
\begin{equation}\label{EQU-4}
\int\rho \Big(\int_1^{\rho }{\frac{P^{2}(s)}{s^{2}}}\;\mathrm{d}s\Big)\;\mathrm{d}x\le P^{*}\big( \|G(\rho)\|_{1}+P(1) \big)-\overline{P}P(1).
\end{equation}
Notice that by making use of \eqref{2.1} and \eqref{2.19*}, there exists a positive constant $M_1$ depending only on $\mu$ such that
\begin{equation}\label{EQU-3}
\begin{aligned}
 &\Big| \int\big((\nabla\d\odot\nabla\d):\nabla \u-\frac{1}{2}(|\nabla\d|^2\mathbb{I}_2):\nabla\u\big)\, \mathrm{d}x \Big|\\
\leq & M\|\nabla\d\|_4^2\|\nabla\u\|_{2}\leq \frac{5}{\mu}\|\nabla\d\|_4^{4}+\frac{\mu }{8}\|\nabla \u\|_2^2 \\
\leq  &M_{1}\|\nabla\d\|_4^{4}+\frac{1}{4}\big(\mu \|\nabla \mathbb{P}\u\|_2^{2}+\frac{1}{\nu }\|\widetilde{F}\|_2^2+\frac{1}{\nu }\|\widetilde{P}\|_2^{2}\big).
\end{aligned}
\end{equation}
Here, the desired inequality can be deduced by noting that $\frac{1}{2}f\cdot g\leq f^2+\frac{1}{16}g^2$ and $f\cdot g\leq 4f^2+\frac{g^2}{16}$.  Using \eqref{L-1}, it yields that
\begin{equation}\label{2.22}
\begin{aligned}
 M_1 \|\nabla\d\|_4^{4}\leq M \|\nabla\d\|_2^2\|\nabla^2\d \|_2^2\leq M_2 \|\nabla^2 \d\|_2^2,
\end{aligned}
\end{equation}
where $M_2$ denotes a generic positive constant depending on $C_0$ and $\mu$.
Thus, choosing $\kappa_0=1+2M_2$,
\begin{equation}\label{eq-Phi-2-2}
\begin{aligned}
\Phi _{2}(t)
& \ge 1+\frac{1 }{4}\big(\mu\|\nabla \mathbb{P}\u\|_2^{2}+\frac{1}{\nu }\|\widetilde{F}\|_2^{2}+\frac{1}{\nu }\|\widetilde{P}\|_2^{2}\big)+\frac{\kappa_0}{2}\|\nabla^2\d\|_2^{2}-M_{1}\|\nabla\d\|_4^{4}\\
& \geq1+\frac{1 }{4}\big(\mu\|\nabla \mathbb{P}\u\|_2^{2}+\frac{1}{\nu }\|\widetilde{F}\|_2^{2}+\frac{1}{\nu }\|\widetilde{P}\|_2^{2}\big)+(\frac{\kappa_0}{2}-M_2)\|\nabla^2\d\|_2^{2}
\\
& = 1
+\frac{1 }{4}\big(\mu\|\nabla \mathbb{P}\u\|_2^{2}+\frac{1}{\nu }\|\widetilde{F}\|_2^{2}+\frac{1}{\nu }\|\widetilde{P}\|_2^{2}\big)+\frac{1}{2}\|\nabla^2\d\|_2^{2}
\\
&\geq \frac{1}{4}\Phi(t)\quad \text{ for all } t\in[0,T].
\end{aligned}
\end{equation}

Now, integrating \eqref{EQU-1} with respect to the time variable $t$ over the interval  $[0,t]$ for $t\in[0,T]$, we can obtain
\begin{align}\label{eq-Phi-3}
\Phi_2(t) +\int_0^t \Theta_1(\tau)\,\mathrm{d}\tau\leq 
\Phi_2(0)+M \int_0^t \big(
  \|\nabla\mathbb{P}\u\|_2^2+\frac{1}{\nu}\|\widetilde{F}\|_2^2+\frac{1}{\nu}\|\widetilde{P}\|_2^2+\|\nabla^2\d\|_2^2\big)\Phi(\tau)\,\mathrm{d}\tau.
\end{align}
We next estimate $\Phi_2(0)$. Since $P(\rho)=a\rho^\gamma$,
$\gamma\geq1$, and
$\int\rho\,\mathrm{d}x=|\mathbb{T}^2|=1$,
Jensen's inequality gives
\begin{equation*}
\begin{aligned}
\int\rho\left(\int_\rho^1
\frac{P^2(s)}{s^2}\,\mathrm{d}s\right)\mathrm{d}x
&=\frac{a^2}{2\gamma-1}
\left(1-\int\rho^{2\gamma}\,\mathrm{d}x\right)
\leq0.
\end{aligned}
\end{equation*}
Consequently, using the definition of $\Phi_1$,
\eqref{EQU-3}, \eqref{2.22}, and
$\kappa_0=1+2M_2$, we obtain
\begin{equation*}
\begin{aligned}
\Phi_2(t)
\leq{}&1+\frac34\left(
\mu\|\nabla\mathbb{P}\u\|_2^2
+\frac1\nu\|\widetilde{F}\|_2^2
+\frac1\nu\|\widetilde{P}\|_2^2
\right)\\
&+\left(\frac12+2M_2\right)
\|\nabla^2\d\|_2^2
+\frac{P^*}{\nu}\|G(\rho)\|_1
+\frac{P(1)}{\nu}\bigl(P^*-P(1)\bigr).
\end{aligned}
\end{equation*}
Let
$\widetilde{P}_0
:=P(\rho_0)-\int P(\rho_0)\,\mathrm{d}x$.
By \eqref{2.1},
\begin{equation*}
\frac1\nu\|\widetilde{F}(0)\|_2^2
\leq
2\nu\|\Dv\u_0\|_2^2
+\frac2\nu\|\widetilde{P}_0\|_2^2.
\end{equation*}
Thus, using \eqref{1.11*} and $\nu\geq\mu$, we infer
\begin{equation}\label{EQU-5}
\begin{aligned}
\Phi_2(0)
\leq{}&M(1+N_0)
+\frac{M}{\nu}\Bigl(
\|\widetilde{P}_0\|_2^2
+P^*\|G(\rho_0)\|_1
+P(1)\bigl(P^*-P(1)\bigr)
\Bigr)
\leq M.
\end{aligned}
\end{equation}
Notice that by Lemmas \ref{Lemma1} and \ref{Lemma2}, it holds 
\begin{align*}
\int_0^T \big(
  \|\nabla\mathbb{P}\u\|_2^2+\frac{1}{\nu}\|\widetilde{F}\|_2^2+\frac{1}{\nu}\|\widetilde{P}\|_2^2+\|\nabla^2\d\|_2^2\big)\,\mathrm{d}t\leq M,
\end{align*}
where $M$ is a generic positive
constant depending only on $\eta$, $\mu$, $a$, $\gamma$, $N_0$, $C_0$, $\overline{\omega}$ and $\|\rho_0\|_\infty$. Thus, combining \eqref{eq-Phi-3}, \eqref{EQU-5} and the inequality $ 
\Phi_2(t)\geq\frac14\Phi(t)$ in \eqref{eq-Phi-2-2}
, and then applying Gronwall’s inequality, yields \eqref{2.5}.
\end{proof}

\begin{Remark}\label{2-2-1} Based on the uniform estimates obtained in Lemma \ref{Lemma1} and Proposition \ref{2-2}, we can obtain from \eqref{2.1.3*}, \eqref{2.19*} and $\widetilde{F}=\nu\Dv \u-\widetilde{P}$ that
\begin{equation*}
\begin{aligned}
&\sup_{0\leq t\leq T}\big(\|\u\|_{H^1}+\|\nabla\d\|_{H^1}+\nu\|\Dv\u\|_2^2\big)\leq M.
\end{aligned}
\end{equation*}
\end{Remark}

The following proposition provides a uniform upper bound for the density when the bulk viscosity coefficient $\nu$ is sufficiently large.
\begin{Proposition}\label{2-3}
Let $(\rho,\u,\d)$ be a smooth solution to \eqref{1.1}--\eqref{1.3} on $[0,T]\times\mathbb{T}^2 $ satisfying \eqref{Sup}, and let $(\rho_0, \u_0, \d_0)$ satisfy \eqref{1.6} and \eqref{1.10}--\eqref{1.12}. Then, there exists a positive number $\nu_{0,1}$ \big( $\nu_{0,1}=\frac{M_4^2}{\log^2\frac32}$\big) depending only on $a$, $\mu$, $\gamma$ , $C_0$, $N_0$, $\overline{\omega}$ and $\|\rho_0\|_\infty$, but independent of $T$
such that if $\nu\ge\nu_0=\max\{\mu,\nu_{0,1}\}$, then 
\begin{equation}\label{2.49}
\underset{0\le t\le T}{\mathop{\sup }}\,\|{\rho (t)}\|_{\infty }\le \frac{3}{2}e^{\frac{\gamma-1}{a\gamma} C_{0}}\|\rho _{0}\|_{\infty }.
\end{equation}
\end{Proposition}

\begin{proof}
Without loss of generality, 
we assume that $\rho$ is smooth and positive, and denote $\rho^*=2e^{\frac{\gamma-1}{a\gamma} C_{0}}\|\rho _{0}\|_{\infty }$. First, we can write according to the definition of $\widetilde{F}$:
\begin{equation*}
\partial _{t}\log \rho +\u\cdot \nabla \log \rho =-\frac{1}{\nu }\big(\widetilde{P}+\widetilde{F}\big).
\end{equation*}
Applying the divergence operator to $\eqref{1.1}_2$ yields
\begin{equation*}
\Delta \widetilde{F}=\partial _{t}(\Dv (\rho \u))+\Dv \big( \Dv (\rho \u\otimes \u)+\Dv( \nabla\d\odot\nabla\d-\frac{1}{2}|\nabla\d|^2\mathbb{I}_2) \big).
\end{equation*}
Following the argument used in \cite{B.D}, denote
\begin{equation*}
\mathbf{E}:=\log \rho -\frac{1}{\nu}{(-\Delta )}^{-1}\Dv(\rho \u).
\end{equation*}
For $\rho>0$, direct computation gives
\begin{equation}\label{2.50}
\partial _{t}\mathbf{E}+\u\cdot \nabla \mathbf{E}+\frac{1}{\nu }\widetilde{P}=-\frac{1}{\nu }[\u^{j},{(-\Delta )}^{-1}\partial _{i}\partial _{j}]\rho \u^{i}-\frac{1}{\nu }\Delta ^{-1}\Dv \Dv \big( \nabla\d\odot\nabla\d-\frac{1}{2}|\nabla\d|^2\mathbb{I}_2\big).
\end{equation}
Fix $x\in \mathbb{T}^2$ and define the corresponding particle path $X$ by
\begin{equation*}
\left\{
\begin{aligned}
&\frac{\mathrm{d} \,X(\tau;t,x)}{\mathrm{d}\tau}=\u(\tau;X(\tau;t,x)), \\
&X(\tau;t,x)\big|_{\tau=t}=x. 
\end{aligned}
\right.
\end{equation*}
Using the fact that $\frac {\mathrm {d}}{\mathrm{d}\tau}f(\tau,X(\tau;t,x))=(f_t+\mathbf u\cdot\nabla f)(\tau,X(\tau;t,x))$, it follows from \eqref{2.50} that
\begin{equation*}
\begin{aligned}
&\frac{\mathrm{d}}{\mathrm{d}\tau}\mathbf{E}(\tau;X(\tau;t,x))+\frac{1}{\nu }\widetilde{P}(\tau;X(\tau;t,x))\\
=&-\frac{1}{\nu }[\u^{j},{(-\Delta )}^{-1}\partial _{i}\partial _{j}]\rho \u^{i}(\tau;X(\tau;t,x))-\frac{1}{\nu }\Delta ^{-1}\Dv \Dv \big( \nabla\d\odot\nabla\d-\frac{1}{2}|\nabla\d|^2\mathbb{I}_2\big)(\tau;X(\tau;t,x)).
\end{aligned}
\end{equation*}
In order to handle $\frac{1}{\nu }\widetilde{P}(\tau;X(\tau;t,x))$, we present an elementary inequality as follows:
\begin{equation*}
P(\rho)=a\rho^\gamma\ge a\gamma\log\rho+a,\quad\text{for all}\;\rho>0,
\end{equation*}
thus, we get $\frac{1}{\nu }\widetilde{P}(\tau;X(\tau;t,x))\ge \frac{a\gamma}{\nu }\log\rho(\tau;X(\tau;t,x))+\frac{a}{\nu}-\frac{1}{\nu }\overline{P}(\tau)$.
\medskip

Then, since the assumed velocity field $ \u $ is smooth, setting $\mathbf{E}^+ := \max(0, \mathbf{E })$, it follows from the definition of $\mathbf{E}^+$ that
\begin{equation*}
\begin{aligned}
  & \frac{\mathrm{d}}{\mathrm{d}\tau}\mathbf{E}^+(\tau;X(\tau;t,x))+\frac{a\gamma}{\nu }\mathbf{E}^{+}(\tau;X(\tau;t,x)) \\
\le  & \frac{1}{\nu }\big| [ \u^{j},{(-\Delta )}^{-1}\partial _{i}\partial _{j} ]\rho \u^{i}(\tau;X(\tau;t,x)) \big|+\frac{a\gamma}{\nu ^{2}}\big| {(-\Delta )}^{-1}\Dv (\rho \u)(\tau;X(\tau;t,x)) \big| \\
 &+\frac{1}{\nu }\big| \Delta ^{-1}\Dv \Dv\big( \nabla\d\odot\nabla\d-\frac{1}{2}|\nabla\d|^2\mathbb{I}_2\big)(\tau;X(\tau;t,x)) \big|+\frac{1}{\nu }(\overline{P}(\tau)-a).
\end{aligned}
\end{equation*}
Integrating the above inequality over $[0,t]$ and applying the maximum principle for the transport equation yields
\begin{equation}\label{2.52}
\begin{aligned}
&\|\mathbf{E}^{+}(t)\|_{\infty } \\
\le &{\text{e}}^{-\frac{a\gamma}{\nu }t}\|\mathbf{E}^{+}(0)\|_{\infty }+\frac{\gamma-1}{a\gamma}(1-{\text{e}}^{-\frac{a\gamma}{\nu }t})C_{0}+\frac{1}{\nu }\int_0^{t}{\text{e}}^{-\frac{a\gamma}{\nu }(t-\tau )}\|[\u^{j},{(-\Delta )}^{-1}\partial _{i}\partial _{j}]\rho \u^{i}(\tau )\|_{\infty }\;\mathrm{d}\tau  \\
& +\frac{a\gamma}{\nu ^{2}}\int_0^{t}{\text{e}}^{-\frac{a\gamma}{\nu }(t-\tau )}\|{(-\Delta )}^{-1}\Dv (\rho \u)(\tau )\|_{\infty }\;\mathrm{d}\tau   \\
 &+\frac{1}{\nu }\int_0^{t}{\text{e}}^{-\frac{a\gamma}{\nu }(t-\tau )}\|\Delta ^{-1}\Dv \Dv\big( \nabla\d\odot\nabla\d-\frac{1}{2}|\nabla\d|^2\mathbb{I}_2\big)(\tau )\|_{\infty }\;\mathrm{d}\tau , \\
\end{aligned}
\end{equation}
where $\overline{P}-a=(\gamma-1)\|G\|_1\le (\gamma -1)C_{0}$.
Now we estimate the terms on the right-hand side of \eqref{2.52}, since the average value of $\rho\u$ is zero, invoking \eqref{2.30-1}, Sobolev embedding inequality $\dot{W}^{1,4}(\mathbb{T}^2)\hookrightarrow L^\infty(\mathbb{T}^2)$ and the Riesz operator, it holds 
\begin{equation}\label{2.53}
\begin{aligned}
\|{(-\Delta )}^{-1}\Dv (\rho \u)\|_{\infty }
\le& M\|{(-\Delta )}^{-1}\nabla \Dv(\rho \u)\|_{4}\\
\le& M\|\rho \u\|_{4}
\le M(\rho^*)\|\u\|_4\leq M \|\nabla\u\|_2.
\end{aligned}
\end{equation}
By the result of Coifman, Lions, Meyer and Semmes (\cite{R.CP}), the following map
\begin{equation*}
\begin{aligned}
 & \mathcal{F}:\; W^{1,r_1}(\mathbb{T}^2)\times L^{r_2}(\mathbb{T}^2)\mapsto W^{1,r_3}(\mathbb{T}^2)
\end{aligned}
\end{equation*}
defined as
\begin{equation*}
\begin{aligned}
 \mathcal{F}(a,b)=[a_j,{(-\Delta )}^{-1}\partial _{i}\partial _{j}]b_i \quad
\end{aligned}
\end{equation*}
is continuous when $\frac{1}{r_3}=\frac{1}{r_1}+\frac{1}{r_2}$, and $1<r_1,r_2,r_3<\infty$.
Then we have the following estimates
\begin{equation}\label{2.54}
\|[\u^{j},{(-\Delta )}^{-1}\partial _{i}\partial _{j}]\rho \u^{i}\|_{\infty }\leq M \|[\widetilde{\u}^{j},{(-\Delta )}^{-1}\partial _{i}\partial _{j}]\rho \u^{i}\|_{W^{1,3}}\le M\|\nabla \u\|_{12}\|\rho \u\|_{4},
\end{equation}
where $[\u^{j},{(-\Delta )}^{-1}\partial _{i}\partial _{j}]\rho \u^{i}=[\widetilde{\u}^{j},{(-\Delta )}^{-1}\partial _{i}\partial _{j}]\rho \u^{i}$ with $\widetilde{\u}=\u-\overline{\u}$.
The factor $\|\nabla u\|_{12}$
 is estimated as follows:
\begin{equation*}
\begin{aligned}
 \|\nabla\u\|_{12}  &\le M
\big(\|\nabla \mathbb{P}\u\|_{12}+\frac{1}{\nu}\|\widetilde{F}\|_{12}+\frac{1}{\nu}\|\widetilde{P}\|_{12}\big)\\
&\le M\big(\|\nabla \mathbb{P}\u\|_2^{\frac{1}{6}}\|\nabla ^{2}\mathbb{P}\u\|_2^{\frac{5}{6}}+\frac{1}{\nu} \| \widetilde{F}\|_2^{\frac{1}{6}}\|\nabla F\|_2^{\frac{5}{6}}+\frac{1}{\nu}\|\widetilde{P}\|_{\infty } \big), \\
 &\le M\big(\|\nabla ^{2}\mathbb{P}\u\|_{2}+\frac{1}{\nu} \|\nabla F\|_{2}+\frac{1}{\nu}\|\widetilde{P}\|_{\infty } \big), 
\end{aligned}
\end{equation*}
where we have used Poincar\'{e}'s inequality.
Thus, in view of \eqref{2.54} and $\frac{1}{\nu}\leq \frac{1}{\mu}$, one may bound
\begin{equation}\label{2.56}
\begin{aligned}
\|[\u^{j},{(-\Delta )}^{-1}\partial _{i}\partial _{j}]\rho \u^{i}\|_{\infty }
\le &M \Big(\|\nabla ^{2}\mathbb{P}\u\|_{2}+\frac{1}{\nu}\big( \|\nabla F\|_{2}+\|\widetilde{P}\|_{\infty } \big)\Big) \|\nabla\u\|_2\\
\leq& M \big(\frac{1}{\mu}\Theta^{\frac{1}{2}}+\frac{1}{\nu}\|\widetilde{P}\|_{\infty}\big)
\|\nabla\u\|_2,
\end{aligned}
\end{equation}
where $\Theta(t)$ is defined in \eqref{2.5}.
Owing to $\int\,\Dv\big( \nabla\d\odot\nabla\d-\frac{1}{2}|\nabla\d|^2\mathbb{I}_2\big)\;\mathrm{d}x=0$, using Riesz operator and Sobolev embedding  inequality $\dot{W}^{1,4}(\mathbb{T}^2)\hookrightarrow L^\infty(\mathbb{T}^2)$, the last term on the right-hand side of \eqref{2.52} can be estimated as follows
\begin{equation}\label{2.58}
\begin{aligned}
   \|\Delta ^{-1}\Dv \Dv\big( \nabla\d\odot\nabla\d-\frac{1}{2}|\nabla\d|^2\mathbb{I}_2\big)\|_{\infty}
 \le & M\|\Delta ^{-1}\nabla \Dv \Dv\big( \nabla\d\odot\nabla\d-\frac{1}{2}|\nabla\d|^2\mathbb{I}_2\big)\|_{4} \\
 \le & M\|\nabla\d\nabla^2 \d\|_{4}\le M\|\nabla\d\|_{\infty}\|\nabla^2 \d\|_{4} \\
\le & M \|\nabla \d\|_2^{\frac{1}{2}}\|\nabla^2\d\|_2^{\frac{1}{2}}\|\nabla^3\d\|_2\le M \Theta^{\frac{1}{2}}.
\end{aligned}
\end{equation}

Now, let $[0, T_0]$ be the maximal subinterval of $[0, T]$ on which the bootstrap bound \eqref{Sup} holds. Plugging the inequalities \eqref{2.53}, \eqref{2.56}, and \eqref{2.58} into \eqref{2.52}, we conclude that
\begin{equation}\label{2.62}
\begin{aligned}
   \|\mathbf{E}^{+}(t)\|_{\infty }
  & \le \|\mathbf{E}^{+}(0)\|_{\infty }+\frac{\gamma-1}{a\gamma} C_{0}+M\frac{a\gamma}{\nu ^{2}}\int_0^{t}e^{-\frac{a\gamma}{\nu }(t-\tau )}\|\nabla \u(\tau )\|_{2}\;\mathrm{d}\tau  \\
 & \quad+M\frac{1}{\nu }\int_0^{t}e^{-\frac{a\gamma}{\nu }(t-\tau )}\big(\mu^{-1}\Theta^{\frac{1}{2}}+\nu^{-1}\|\widetilde{P}\|_{\infty}\big)
\|\nabla \u\|_2 \;\mathrm{d}\tau +M\frac{1}{\nu }\int_0^{t}e^{-\frac{a\gamma}{\nu }(t-\tau )}\Theta^{\frac{1}{2}}\;\mathrm{d}\tau . \\
 &\leq \|\mathbf{E}^{+}(0)\|_{\infty }+\frac{\gamma-1}{a\gamma} C_{0}+M\frac{1+a\gamma}{\nu^2}\Big(\frac{\nu}{a\gamma}\Big)^{\frac{1}{2}}\Big(\int_0^{t}\|\nabla \u(\tau )\|_2^{2}\;\mathrm{d}\tau\Big)^{\frac{1}{2}}\\
&\quad
+M\frac{1}{\nu}\Big(\frac{\nu}{a\gamma}\Big)^{\frac{1}{2}}\Big(\int_0^{t}\Theta(\tau)\;\mathrm{d}\tau\Big)^{\frac{1}{2}}
 +
 M\frac{1}{\nu}\Big(\frac{\nu}{a\gamma}\Big)^{\frac{1}{4}}
\Big(\int_0^{t}\Theta(\tau)\;\mathrm{d}\tau\Big)^{\frac{1}{2}}\Big(\int_0^{t}\|\nabla \u(\tau)\|_2^{4}\;\mathrm{d}\tau\Big)^{\frac{1}{4}}.
\end{aligned}
\end{equation}
Since the terms on the right-hand side of \eqref{2.62}  may be bounded in terms of the data according to the basic energy estimate \eqref{1.11} and to \eqref{2.5}, we eventually get if $\nu $ is large enough:
\begin{equation}\label{2.63}
\|\mathbf{E}^{+}(t)\|_{\infty }\le \|\mathbf{E}^{+}(0)\|_{\infty }+M\nu ^{-\frac{1}{2}}+\frac{\gamma-1}{a\gamma} C_{0}.
\end{equation}
The definition of $\mathbf{E}^+$ together with \eqref{2.5} and \eqref{2.53} gives
\begin{equation}\label{2.64}
\log \rho \le \mathbf{E}^{+}+\frac{1}{\nu}\|\Delta ^{-1}\Dv (\rho \u)\|_{\infty }\le \mathbf{E}^{+}+\nu^{-1}M.
\end{equation}
Combining \eqref{2.63} and \eqref{2.64}, we can find a positive constant $M_4$ depending on $a$, $\mu$, $\gamma$ , $C_0$, $N_0$, $\overline{\omega}$ and $\|\rho_0\|_\infty$ such that
\begin{equation}\label{2.65}
\log\|\rho(t)\|_\infty
\leq
\log\|\rho_0\|_\infty
+M_4\nu^{-1/2}
+\frac{\gamma-1}{a\gamma}C_0,
\qquad t\in[0,T_0].
\end{equation}
If $\nu$ is large enough and satisfies $M_4\nu^{-1/2}\leq\log\frac32$,
then \eqref{2.65} gives
\begin{equation*}
\sup_{0\leq t\leq T_0}\|\rho(t)\|_\infty
\leq
\frac32
\exp\left(\frac{\gamma-1}{a\gamma}C_0\right)
\|\rho_0\|_\infty
=
\frac34\rho^*
<
\rho^*.    
\end{equation*}
Consequently, the bootstrap bound is improved, and
the continuation argument gives $T_0=T$.
This completes the proof of Proposition \ref{2-3}.
\end{proof}
\medskip

The following proposition is devoted to the proof of the weighted estimates, which will be fundamental in the proof of
uniqueness.
\begin{Proposition}\label{2-4}
Let $T>0$, and let $(\rho,\u,\d)$ be a smooth solution to
\eqref{1.1}--\eqref{1.3} on $[0,T]\times\mathbb{T}^{2}$ satisfying
\eqref{Sup}. Assume that $\nu\geq\nu_{0}$, where $\nu_{0}$ is
determined in Proposition \ref{2-3}. Then
\begin{equation}\label{2.66}
\begin{aligned}
  & \underset{0\le t\le T}{\mathop{\sup }}\,t(\|\sqrt{\rho }\dot{\u}\|_2^{2}+\|\nabla\d_t\|_2^2+\|\nabla^3\d\|_2^2)+\int_0^{T}t( \|\nabla \dot{\u}\|_2^{2}+ \|\d_{tt}\|_2^2+\|\nabla^2\d_t\|_2^2)\;\mathrm{d}t  \le M(T),
\end{aligned}
\end{equation}
where $M(T)$ depends on $a$, $\mu$, $\gamma$ , $C_0$, $N_0$, $T$, $\overline{\omega}$ and $\|\rho_0\|_\infty$ and is independent of $\nu$.
\end{Proposition}

\begin{proof}
 \emph{StepI:} First, rewrite the momentum equation as follows:
\begin{equation}\label{2.67}
\rho \dot{\u}-\mu \Delta \u-\frac{(\nu -\mu )}{\nu}\nabla \widetilde{F}+\frac{\mu}{\nu}\nabla \widetilde{P}+\Dv\big( \nabla\d\odot\nabla\d-\frac{1}{2}|\nabla\d|^2\mathbb{I}_2\big)=0.
\end{equation}
Applying $\partial/\partial t + \text{div}(\mathbf{u}\cdot)$ to the $i$-th component of \eqref{2.67}, and multiplying the resultant by $\dot{\mathbf{u}}^i$, one gets by some simple calculations that
\begin{equation}\label{2.68}
    \begin{aligned}
\frac{1}{2} \Big( \int \rho|\dot{\mathbf{u}}^i|^2 \,\mathrm{d}x \Big)_t 
=& \mu \int \dot{\mathbf{u}}^i \big( \Delta \mathbf{u}_t^i + \text{div}(\mathbf{u}\Delta \mathbf{u}^i) \big) \,\mathrm{d}x \\
&+ \frac{\nu-\mu}{\nu}\int \dot{\mathbf{u}}^i \big( \partial_t \partial_i \widetilde{F} + \text{div}( \mathbf{u}\partial_i \widetilde{F} )\big) \,\mathrm{d}x\\
&-\frac{\mu}{\nu} \int \dot{\mathbf{u}}^i \big( \partial_i \widetilde{P}_t + \text{div}(\mathbf{u}\partial_i \widetilde{P}) \big) \,\mathrm{d}x\\
&+ \frac{1}{2} \int \dot{\mathbf{u}}^i \big( \partial_t \partial_i |\nabla \mathbf{d}|^2 + \text{div}(\mathbf{u}\partial_i |\nabla \mathbf{d}|^2) \big) \,\mathrm{d}x\\
&- \int \dot{\mathbf{u}}^i \big( \partial_t \text{div}(\partial_j \mathbf{d} \cdot \partial_i \mathbf{d}) + \text{div}(\mathbf{u}\text{div}\,(\partial_i \mathbf{d} \cdot \partial_j \mathbf{d})) \big) \,\mathrm{d}x \triangleq \sum_{k=1}^5 S_k.         \end{aligned}
\end{equation}
We shall estimate each term in \eqref{2.68}. In what follows, $M$ denotes a positive constant depending on $a$, $\mu$, $\gamma$, $C_0$, $N_0$, $\overline{\omega}$ and $\|\rho_0\|_\infty$, but independent of  $T$ and $\nu$. Integration by parts gives the following estimate for $S_1$:
\begin{equation*}
 \begin{aligned}
S_1 &= \mu \int \dot{\mathbf{u}}^i \big[ \partial_t \Delta \mathbf{u}^i + \text{div}(\mathbf{u}\Delta \mathbf{u}^i) \big] \,\mathrm{d}x \\
&= -\mu \int \big( |\nabla \dot{\mathbf{u}}|^2 + \partial_j \dot{\mathbf{u}}^i \partial_k \mathbf{u}^k \partial_j \mathbf{u}^i - \partial_j \dot{\mathbf{u}}^i \partial_j \mathbf{u}^k \partial_k \mathbf{u}^i - \partial_j \mathbf{u}^i \partial_j \mathbf{u}^k \partial_k \dot{\mathbf{u}}^i \big) \,\mathrm{d}x \\
&\leq -\mu \|\nabla \dot{\mathbf{u}}\|_2^2+2\mu\|\nabla\dot{\mathbf{u}}\|_{2}\|\nabla\u\|_4\big(\|\Dv\u\|_4+\|\nabla\u\|_4\big)\\
&\leq -\mu \|\nabla \dot{\mathbf{u}}\|_2^2 + \eta \|\nabla \dot{\mathbf{u}}\|_2^2+M\big(\|\nabla \mathbf{u}\|_4^4+\|\Dv \mathbf{u}\|_4^4\big),
\end{aligned}   
\end{equation*}
where $\eta>0$ is a small constant to be fixed later.
Based on \eqref{2.7} and integration by parts, we deduce that
\begin{equation*}
\begin{aligned}
 S_3&= - \frac{\mu}{\nu}\int \dot{\mathbf{u}}^i \big( \partial_i P_t + \text{div}(\mathbf{u}\partial_i P) \big) \,\mathrm{d}x \\
 &=\frac{\mu}{\nu}\int \partial_i\dot{\u}^iP_t\,\mathrm{d}x+\frac{\mu}{\nu}\int \partial_k\dot{\u}^i\u^k\partial_i P\,\mathrm{d}x\\
 &=-\frac{\mu}{\nu}\int \partial_k\dot{\u}^i\partial_i\u^k P\,\mathrm{d}x-\frac{\mu}{\nu}\int \Dv\dot{\u}\,h(\rho)\Dv\u\,\mathrm{d}x\\
 &\leq \eta \|\nabla\dot{\u}\|_2^2+M\frac{1}{\nu^2}\big( \|\nabla\u\|_2^2+\|\Dv\u\|_2^2\big)\\
  &\leq \eta \|\nabla\dot{\u}\|_2^2+M.
\end{aligned}    
\end{equation*}
In order to estimate $S_2$, we write
\begin{align}\label{eq:3.34}
\sum_{i=1}^2 \partial_i \u\cdot \nabla \dot{\u}^i = \Dv \u\, \Dv \dot{\u} + \nabla \u^1\cdot \nabla^\perp \dot{\u}^2 - \nabla \u^2\cdot \nabla^\perp \dot{\u}^1.
\end{align}
Integrating by parts and using \eqref{2.25}, \eqref{Hardy} and \eqref{eq:3.34} yields that (for any $\eta>0$)
\begin{align}
S_2 &= \frac{\nu - \mu}{\nu} \int \dot{\u}^i \big[\partial_i \partial_t \widetilde{F} + \Dv(\u\partial_i \widetilde{F})\big] \;\mathrm{d}x \nonumber\\
&= -\frac{\nu - \mu}{\nu} \int \Dv\dot{\u}\, \dot{\widetilde{F}} \,\;\mathrm{d}x + \frac{\nu - \mu}{\nu} \int -\widetilde{F} \Dv\dot{\u}\Dv\u+ \widetilde{F}\partial_i\u\cdot\nabla\dot{\u}^i\;\mathrm{d}x \label{eq:2.39}\\
&= -\frac{\nu - \mu}{\nu} \int \Dv\dot{\u}\, \dot{\widetilde{F}} \,\;\mathrm{d}x + \frac{\nu - \mu}{\nu} \int \widetilde{F} \left( \nabla \u^1\cdot \nabla^\perp \dot{\u}^2 - \nabla \u^2\cdot \nabla^\perp \dot{\u}^1 \right) \;\mathrm{d}x\nonumber\\
&\leq -\frac{\nu - \mu}{\nu} \int \Dv\dot{\u}\, \dot{\widetilde{F}} \,\;\mathrm{d}x + M\|\nabla F\|_{2} \|\nabla \u\|_{2} \|\nabla \dot{\u}\|_{2} \nonumber\\
&\leq -\frac{\nu - \mu}{\nu} \int \Dv\dot{\u}\, \dot{\widetilde{F}} \,\;\mathrm{d}x + \frac{\eta}{3} \|\nabla \dot{\u}\|_2^2 + M\big(\|\sqrt{\rho}\dot{\u}\|_2^2+\|\nabla^3\d\|_2^2\big) \|\nabla \u\|_2^2
\nonumber,
\end{align}
where the last inequality follows from \eqref{L-1} and 
\begin{equation}\label{DEU-1}
 \|\nabla\d\nabla^2\d\|_2\leq \|\nabla\d\|_\infty\|\nabla^2\d\|_2\leq c_0\|\nabla\d\|_2\|\nabla^3\d\|_2\leq M \|\nabla^3\d\|_2.   
\end{equation}
To estimate the first term in the last line of the above inequality, we use the definition of $\widetilde{F}$ and \eqref{2.7}  to obtain
\begin{align*}
\Dv\dot{\u} &= (\Dv \u)_{t} + \partial_j \u^i \partial_i \u^j + \u\cdot \nabla \Dv \u \nonumber\\
&= \frac{1}{\nu}\dot{\widetilde{F}} + \frac{1}{\nu}(P_{t} + \u\cdot \nabla P) - \frac{1}{\nu}\overline{P}_{t} + \partial_j \u^i \partial_i \u^j \nonumber\\
&= \frac{1}{\nu}\dot{\widetilde{F}} - \frac{1}{\nu}\gamma P \Dv \u + \frac{\gamma - 1}{\nu} \int P \Dv \u \,\;\mathrm{d}x + \partial_j \u^i \partial_i \u^j,
\end{align*}
which together with Young’s inequality yields
\begin{align}\label{eq:3.36}
&\quad-\frac{\nu - \mu}{\nu} \int \Dv\dot{\u}\, \dot{\widetilde{F}} \,\;\mathrm{d}x\nonumber\\
&= -\frac{\nu - \mu}{\nu^2} \|\dot{\widetilde{F}}\|_2^2 - \frac{\nu - \mu}{\nu} \int \dot{\widetilde{F}}\partial_j \u^i \partial_i \u^j \,\;\mathrm{d}x \\
&\quad + \frac{\nu - \mu}{\nu^2}\gamma \int \dot{\widetilde{F}}P \Dv \u \,\;\mathrm{d}x - \frac{\nu - \mu}{\nu^2}(\gamma - 1) \int P \Dv \u \,\;\mathrm{d}x \int \dot{\widetilde{F}} \,\;\mathrm{d}x \nonumber\\
&\leq -\frac{\nu - \mu}{2\nu^2} \|\dot{\widetilde{F}}\|_2^2 + M\frac{1}{\nu} \|\Dv \u\|_2^2 - \frac{\nu - \mu}{\nu} \int \dot{\widetilde{F}}\partial_j \u^i \partial_i \u^j \,\;\mathrm{d}x\nonumber.
\end{align}
The last term in \eqref{eq:3.36} can be estimated as follows:
\begin{align*}
\int \dot{\widetilde{F}}\partial_j \u^i \partial_i \u^j \,\;\mathrm{d}x
&= \int (\widetilde{F}_{t} + \u\cdot \nabla \widetilde{F})\partial_j \u\cdot \nabla \u^j \,\;\mathrm{d}x \nonumber\\
&=\frac{\mathrm{d}}{\mathrm{d}t}\Big( \int \widetilde{F}\partial_j \u\cdot \nabla \u^j \,\;\mathrm{d}x \Big) - 2\int \widetilde{F}\partial_j \u\cdot \nabla \dot{\u}^j \,\;\mathrm{d}x \nonumber\\
&\quad + 2\int \widetilde{F}\partial_j \u\cdot \nabla \u\cdot \nabla \u^j \,\;\mathrm{d}x - \int \widetilde{F}\Dv \u \partial_j \u\cdot \nabla \u^j \,\;\mathrm{d}x,
\end{align*}
which further gives
\begin{align}\label{eq:3.38}
&-\frac{\nu - \mu}{\nu} \int \dot{\widetilde{F}}\partial_j \u^i \partial_i \u^j \,\;\mathrm{d}x\nonumber\\
=& -\frac{\nu - \mu}{\nu} \frac{\mathrm{d}}{\mathrm{d}t}\Big( \int \widetilde{F}\partial_j \u\cdot \nabla \u^j \,\;\mathrm{d}x \Big) + \frac{2(\nu - \mu)}{\nu} \int \widetilde{F}\partial_j \u\cdot \nabla \dot{\u}^j \,\;\mathrm{d}x \nonumber\\
& - \frac{2(\nu - \mu)}{\nu} \int \widetilde{F}\partial_j \u\cdot \nabla \u\cdot \nabla \u^j \,\;\mathrm{d}x + \frac{\nu - \mu}{\nu} \int \widetilde{F}\Dv \u \partial_j \u\cdot \nabla \u^j \,\;\mathrm{d}x\nonumber\\
:=&-\frac{\nu - \mu}{\nu} \frac{\mathrm{d}}{\mathrm{d}t}\Big( \int \widetilde{F}\partial_j \u\cdot \nabla \u^j \,\;\mathrm{d}x \Big)+S_{2,1}+S_{2,2}+S_{2,3}.
\end{align}
For $S_{2,1}$ in \eqref{eq:3.38}, it follows from \eqref{2.7}, \eqref{eq:3.34} and the definition of $\widetilde{F}$  that
\begin{equation*}
\begin{aligned}
&\frac{2(\nu - \mu)}{\nu} \int \widetilde{F}\partial_j \u\cdot \nabla \dot{\u}^j \,\;\mathrm{d}x\\
\leq& \Big|2\int \widetilde{F} \big( \Dv \u \Dv \dot{\u} + \nabla \u^1\cdot \nabla^\perp \dot{\u}^2 - \nabla \u^2\cdot \nabla^\perp \dot{\u}^1 \big) \;\mathrm{d}x\Big| \\
\leq& M\|\nabla \dot{\u}\|_{2} \|\widetilde{F}\|_{4} \|\Dv \u\|_{4} + M\|\nabla F\|_{2} \|\nabla \dot{\u}\|_{2} \|\nabla \u\|_{2} \\
\leq &\frac{\eta}{3} \|\nabla \dot{\u}\|_2^2 + M\big(\frac{1}{\nu^2}\|\widetilde{F}\|_4^4 + \frac{1}{\nu^2}\|\widetilde{P}\|_4^4\big) + M\big(\|\sqrt{\rho}\dot{\u}\|_2^2+\|\nabla^3\d\|_2^2\big) \|\nabla \u\|_2^2.
\end{aligned}
\end{equation*}
Next, consider $S_{2,2}$ in \eqref{eq:3.38}. The following identity will be used:
\begin{align*}
\sum_{j=1}^2 \partial_j \u \cdot \nabla \u \cdot \nabla \u^j = (\Dv \u)^3 + 3\Dv \u \, \nabla \u^1 \cdot \nabla^\perp \u^2.
\end{align*}
Consequently,
\begin{align}\label{eq:3.40}
-\frac{2(\nu - \mu)}{\nu} \int \widetilde{F} \partial_j \u \cdot \nabla \u \cdot \nabla \u^j \;\mathrm{d}x
&\leq 2 \Big|\int \widetilde{F} \big( (\Dv \u)^3 + 3\Dv \u \, \nabla \u^1 \cdot \nabla^\perp \u^2 \big) \;\mathrm{d}x \Big|\nonumber\\
&\leq M\|\widetilde{F}\|_{4} \|\Dv \u\|_{4} \|\nabla \u\|_4^2 \nonumber\\
&\leq M\big(\frac{1}{\nu^2} \|\widetilde{F}\|_4^4 + \frac{1}{\nu^2} \|\widetilde{P}\|_4^4 +\|\nabla \u\|_4^4\big).
\end{align}
Analogously to the derivation of \eqref{eq:3.40}, we easily see that
\begin{align*}
S_{2,3}
\leq M\big(\frac{1}{\nu^2} \|\widetilde{F}\|_4^4 + \frac{1}{\nu^2} \|\widetilde{P}\|_4^4 + \|\nabla \u\|_4^4\big).
\end{align*}
Substituting the estimates for $S_{2,1}$, $S_{2,2}$ and $S_{2,3}$ into \eqref{eq:3.38} gives
\begin{align}\label{eq:3.43}
&-\frac{\nu - \mu}{\nu} \int \dot{\widetilde{F}} \partial_j \u^i \partial_i \u^j \;\mathrm{d}x\nonumber\\
&\leq -\frac{\nu - \mu}{\nu} \frac{\mathrm{d}}{\mathrm{d}t}\Big( \int \widetilde{F} \partial_j \u \cdot \nabla \u^j \;\mathrm{d}x \Big) + \frac{\eta}{3} \|\nabla \dot{\u}\|_2^2+M\big(\frac{1}{\nu^2} \|\widetilde{F}\|_4^4 + \frac{1}{\nu^2} \|\widetilde{P}\|_4^4 +\|\nabla \u\|_4^4\big)\nonumber\\
&\quad+M\big(\|\sqrt{\rho}\dot{\u}\|_2^2+\|\nabla^3\d\|_2^2\big) \|\nabla \u\|_2^2.
\end{align}
Combining \eqref{eq:2.39}, \eqref{eq:3.36} and \eqref{eq:3.43} yields
\begin{align*}
S_2
&\leq -\frac{\nu - \mu}{\nu} \frac{\mathrm{d}}{\mathrm{d}t}\Big( \int \widetilde{F} \partial_j \u \cdot \nabla \u^j \;\mathrm{d}x \Big) - \frac{\nu - \mu}{2\nu^2} \|\dot{\widetilde{F}}\|_2^2+ \eta \|\nabla \dot{\u}\|_2^2\nonumber\\
&\quad+M\big(\frac{1}{\nu^2} \|\widetilde{F}\|_4^4 + \frac{1}{\nu^2} \|\widetilde{P}\|_4^4 +\frac{1}{\nu}\|\Dv\u\|_2^2+\|\nabla \u\|_4^4\big)
+M\big(\|\sqrt{\rho}\dot{\u}\|_2^2+\|\nabla^3\d\|_2^2\big) \|\nabla \u\|_2^2.
\end{align*}
To estimate $S_4$ and $S_5$, one derives from \eqref{L-1} and Remark \ref{2-2-1} that
\begin{equation*}
\begin{aligned}
S_4+S_5 =& -\int \partial_i \dot{\mathbf{u}}^i \partial_j \mathbf{d} \cdot \partial_j\mathbf{d}_t \, \,\mathrm{d}x- \frac{1}{2} \int \partial_k\dot{\mathbf{u}}^i \mathbf{u}^k \cdot \partial_i |\nabla \mathbf{d}|^2\,\mathrm{d}x \\
&+\int \partial_j \dot{\mathbf{u}}^i \partial_j \mathbf{d}_t \cdot \partial_i\mathbf{d} \, \,\mathrm{d}x+\int \partial_j \dot{\mathbf{u}}^i \partial_j \mathbf{d} \cdot \partial_i\mathbf{d}_t \, \,\mathrm{d}x + \int \partial_k\dot{\mathbf{u}}^i \mathbf{u}^k \cdot \partial_j (\partial_j \mathbf{d} \cdot \partial_i\d) \,\mathrm{d}x \\
\leq & M\big(\|\nabla\d\|_4\|\nabla\d_t\|_4\|\nabla\dot{\u}\|_2+\|\nabla\d\|_8\|\u\|_8\|\nabla^2\d\|_4\|\nabla\dot{\u}\|_2\big)\\
\leq & M\big(\|\nabla\d_t\|_2^{\frac{1}{2}}\|\nabla^2\d_t\|_2^{\frac{1}{2}} \|\nabla\dot{\u}\|_2+\|\nabla^2\d\|_2^{\frac{1}{2}}\|\nabla^3\d\|_2^{\frac{1}{2}}\|\nabla\dot{\u}\|_2\big)\\
\leq & \eta\big(\|\nabla\dot{\u}\|_2^2+\|\nabla^2\d_t\|_2^2\big)+M\big(1+\|\nabla\d_t\|_2^2+\|\nabla^3\d\|_2^2\big).
\end{aligned}    
\end{equation*}
Substituting the estimates of $S_1$--$S_5$ into \eqref{2.68} gives that
\begin{equation}\label{2.73}
\begin{aligned}
&\frac{\mathrm{d}}{\mathrm{d}t}\Big( \int \frac{1}{2}\rho |\dot{\u}|^2 \;\mathrm{d}x + \frac{\nu - \mu}{\nu}\int \widetilde{F} \partial_j \u \cdot \nabla \u^j \;\mathrm{d}x \Big)+ \mu\|\nabla \dot{\u}\|_2^2 + \frac{\nu - \mu}{2\nu^2} \|\dot{\widetilde{F}}\|_2^2 \\
 \leq& 4\eta\big(\|\nabla\dot{\u}\|_2^2+\|\nabla^2\d_t\|_2^2\big)\\
&+M\big(1+ \frac{1}{\nu^2} \|\widetilde{F}\|_4^4 + \frac{1}{\nu^2} \|\widetilde{P}\|_4^4 +\| \nabla \u\|_4^4+\| \Dv \u \|_4^4+\|\sqrt{\rho}\dot{\u}\|_2^2+\|\nabla\d_t\|_2^2+\|\nabla^3\d\|_2^2\big).
\end{aligned}
\end{equation}
\smallbreak

\noindent\emph{StepII:}
By virtue of \eqref{L-1}, \eqref{1.11}, \eqref{2.5}, \eqref{eq1} and the inequality $\|\nabla\d\|_2\leq M  \|\nabla^2\d\|_2$, it is easy to see that
\begin{equation}\label{3.11}
\begin{aligned}
\int_0^T \|\nabla\d_t\|_2^2\,\mathrm{d}t
&\le M\int_0^T
\Bigl(
\|\nabla^3\d\|_2^2
+\|\nabla(\u\cdot\nabla\d)\|_2^2
+\|\nabla(|\nabla\d|^2\d)\|_2^2
\Bigr)\,\mathrm{d}t\\
&\le M+M\int_0^T
\Bigl(
\|\nabla\d\|_2^2
+\|\nabla^2\d\|_2^2
+\|\nabla^3\d\|_2^2
\Bigr)\,\mathrm{d}t
\le M.
\end{aligned}
\end{equation}
Then, differentiating $\eqref{1.1}_3$ with respect to $t$, we obtain
\begin{equation}\label{DEu}
\d_{tt} - \Delta \d_t = \big( |\nabla \d|^2 \d - \u \cdot \nabla \d \big)_t.    
\end{equation}
By direct computations and the fact that $|\d| = 1$, one obtains from \eqref{DEu} and $\dot{\u}=\u_t+\u\cdot\nabla\u$ that
\begin{equation}\label{Equ2}
\begin{aligned}
&\frac{\mathrm{d}}{\mathrm{d}t} \int |\nabla \d_t|^2\,\mathrm{d}x + \int\bigl( |\d_{tt}|^2 + |\nabla^2 \d_t|^2 \bigr)\,\mathrm{d}x \\
\leq& M\int |\nabla \d_t|^2\,\mathrm{d}x 
+ M \int  |\dot{\u}|^2 |\nabla \d|^2 \,\mathrm{d}x
+ M \int  |\u|^2 |\nabla \u|^2 |\nabla \d|^2 \,\mathrm{d}x \\
&+  M\int |\u|^2 |\nabla \d_t|^2 \,\mathrm{d}x 
+ M \int |\nabla \d|^2 |\nabla \d_t|^2 \,\mathrm{d}x 
+ M \int |\d_t|^2 |\nabla\d|^4 \,\mathrm{d}x\\
\triangleq& \sum_{i=1}^6 D_i.
\end{aligned}
\end{equation}
For $D_1$ and $D_2$, by using \eqref{L-1} and Poincar\'{e}'s inequality, one gets
\begin{equation*}
\begin{aligned}
|D_1|+|D_2|
\leq& M\big(\|\nabla\d_t\|_2^2+\|\dot{\u}\|_4^2\|\nabla\d\|_4^2\big)\\
\leq & M\big(\|\nabla\d_t\|_2^2+\|\nabla\dot{\u}\|_2^2\big),
\end{aligned}    
\end{equation*}
where we have used the fact that
\begin{equation}\label{Dot}
\|\dot{\u}\|_4^2\leq M \|\dot{\u}\|_2\|\dot{\u}\|_{H^1}\leq M (1+\|\widetilde{\rho}\|_2)\|\nabla\dot{\u}\|_2^2\leq M(\rho^*) \|\nabla\dot{\u}\|_2^2,   
\end{equation}
due to \eqref{L-3} and \eqref{2.1.3*} with $\int \rho \dot{\u}\,\mathrm{d}x=0$. To bound $D_3$-$D_5$, it follows from \eqref{L-1}, \eqref{2.30-1} and Remark \ref{2-2-1} that
\begin{equation*}
 \begin{aligned}
&|D_3|+|D_4|+|D_5|\\
\leq & M\big(\|\u\|_8^2\|\nabla\u\|_4^2\|\nabla\d\|_8^2+\|\u\|_4^2\|\nabla\d_t\|_4^2+\|\nabla\d\|_4^2\|\nabla\d_t\|_4^2 \big)\\
\leq & M\big(\|\nabla\u\|_4^2+\|\nabla\d_t\|_2\|\nabla^2\d_t\|_2 \big)\\
\leq & \frac{1}{2} \|\nabla^2\d_t\|_2^2+M\big(1+\|\nabla\u\|_4^4+\|\nabla\d_t\|_2^2 \big).
 \end{aligned}   
\end{equation*}
Using $\eqref{1.1}_3$, \eqref{L-1} and $|\d|=1$, we deduce that
\begin{equation*}
\begin{aligned}
|D_6|\leq &M\|\d_t\|_2^2\|\nabla\d\|_\infty^4\\
\leq &M\big(\|\nabla^2\d\|_2^2+\|\nabla\d\|_4^4+\|\u\|_4^2\|\nabla\d\|_4^2)\|\nabla^3\d\|_2^2\\
\leq &M\|\nabla^3\d\|_2^2.
\end{aligned}    
\end{equation*}
Substituting the estimates of $D_1$--$D_6$ into $\eqref{Equ2}$ results in (taking $\eta$ suitably small)
\begin{equation}\label{eq:3.46}
 \begin{aligned}
 &\frac{\mathrm{d}}{\mathrm{d}t}\|\nabla\d_t\|_2^2+\big(\|\d_{tt}\|_2^2+\frac{1}{2}\|\nabla^2\d_t\|_2^2\big)\\
 \leq& M_5\|\nabla\dot{\u}\|_2^2+ M\big(1+\|\nabla\d_t\|_2^2+\|\nabla^3\d\|_2^2+\|\nabla\u\|_4^4\big), 
 \end{aligned}   
\end{equation}
where $M_5$ is a positive constant which depends only on $a$, $\mu$, $\gamma$ , $C_0$, $N_0$, $\overline{\omega}$ and $\|\rho_0\|_\infty$. Additionally, combining \eqref{eq1} with the standard $L^2$--estimate of elliptic system, we obtain from \eqref{L-1} and Remark \ref{2-2-1} that 
\begin{equation}\label{eq:3.47}
\begin{aligned}
\|\nabla^3\d\|_2^2&\leq M\big(\|\nabla\d_t\|_2^2+\|\nabla\u\nabla\d\|_2^2+\|\u\nabla^2\d\|_2^2+\|\nabla|\nabla\d|^2\|_2^2+\|\nabla\d\|_6^6\big)\\
&\leq M\big(\|\nabla\d_t\|_2^2+\|\nabla\u\|_2^2\|\nabla\d\|_\infty^2+\|\u\|_4^2\|\nabla^2\d\|_4^2+\|\nabla\d\|_4^2\|\nabla^2\d\|_4^2+\|\nabla\d\|_6^6\big)\\
&\leq \frac{1}{2} \|\nabla^3\d\|_2^2+M\big(1+\|\nabla\d_t\|_2^2\big).
\end{aligned}    
\end{equation}
Multiplying \eqref{2.73} by $2M_5/\mu$ and combining the resulting inequality with \eqref{eq:3.46}, we obtain, upon choosing $0<\eta\leq\min\{\mu/16,\mu/(32M_5)\}$,
\begin{equation}\label{eq:3.48}
\begin{aligned}
\frac{\mathrm{d}}{\mathrm{d}t}& \Big( \int\big(\frac{M_5}{\mu}\rho |\dot{\u}|^2+|\nabla\d_t|^2\big) \;\mathrm{d}x +\frac{2M_5}{\mu}\frac{\nu - \mu}{\nu} \int \widetilde{F} \partial_j \u \cdot \nabla \u^j \;\mathrm{d}x \Big)\\
&+ \frac{M_5}{\mu}\frac{\nu - \mu}{\nu^2} \|\dot{\widetilde{F}}\|_2^2+ \frac{M_5}{2}\|\nabla \dot{\u}\|_2^2+\|\d_{tt}\|_2^2+\frac{1}{4}\|\nabla^2\d_t\|_2^2\\
\leq &M\big(1+ \frac{1}{\nu^2} \|\widetilde{F}\|_4^4 + \frac{1}{\nu^2} \|\widetilde{P}\|_4^4 +\|\nabla \u\|_4^4+\|\Dv\u\|_4^4+ \| \nabla^3 \d\|_2^2 +\|\sqrt{\rho}\dot{\u}\|_2^2+\|\nabla\d_t\|_2^2\big).
\end{aligned}
\end{equation}
Moreover, we observe from \eqref{2.25}, \eqref{eq:3.34}, \eqref{DEU-1}, \eqref{eq:3.47} and the definition of $\widetilde{F}$ that
\begin{align}\label{eq:3.49}
&\frac{(\nu - \mu)}{\nu} \Big|\int \widetilde{F} \partial_i \u \cdot \nabla \u^i \;\mathrm{d}x\Big|\nonumber\\
\leq &\frac{(\nu - \mu)}{\nu}\Big|\int \widetilde{F} (\Dv \u)^2 \;\mathrm{d}x\Big| +  2\frac{(\nu - \mu)}{\nu}\Big|\int \widetilde{F} \nabla \u^1 \cdot \nabla^\perp \u^2 \;\mathrm{d}x\Big| \nonumber\\
\leq &M\big(\frac{1}{\nu^2} \|\widetilde{F}\|_3^3 + \frac{1}{\nu^2} \|\widetilde{P}\|_3^3 + \|\nabla F\|_{2} \|\nabla \u\|_2^2\big) \nonumber\\
\leq& M\big( \frac{1}{\nu^2} \|\widetilde{F}\|_2^2 \big(\|\nabla F\|_{2}+\|\nabla^3\d\|_2\big) + \frac{1}{\nu^2} \|\widetilde{P}\|_3^3 + (\|\sqrt{\rho}\dot{\u}\|_{2}+\|\nabla^3\d\|_2) \|\nabla \u\|_2^2\big) \nonumber\\
\leq &\frac{1}{4}\|\sqrt{\rho}\dot{\u}\|_2^2+\frac{\mu}{4M_5}\|\nabla\d_t\|_2^2+M.
\end{align}
\smallbreak

\noindent\emph{StepIII:} 
To close the estimate in \eqref{eq:3.48}, we use $\nu\ge\nu_0$, \eqref{L-1}, \eqref{2.1.6}, \eqref{2.5}, the Helmholtz decomposition \eqref{2.27} and the boundedness
of the Riesz operators,  to obtain
\begin{equation}\label{2.688}
\begin{aligned}
\int_0^T
\bigl(\|\nabla\u\|_4^4+\|\Dv\u\|_4^4\bigr)\,\mathrm{d}t
\le& M\int_0^T
\Big(
\|\nabla\mathbb{P}\u\|_4^4
+\frac{\|\widetilde P\|_4^4+\|\widetilde F\|_4^4}{\nu^2}
\Big)\,\mathrm{d}t\\
\le& M\big(\|\nabla \mathbb{P}\u\|_{L^{\infty }(0,T;L^{2})}^{2}\|\nabla^2 \mathbb{P}\u\|_{L^{2 }(0,T;L^{2})}^2+\|\widetilde{P}\|_{L^{\infty }(0,T;L^{\infty})}^{2} \frac{1}{\nu^2}\|\widetilde{P}\|_{L^{2}(0,T;L^{2})}^{2}\\
&+\frac{1}{\nu}\|\widetilde{F}\|_{L^{\infty }(0,T;L^{2})}^{2}\frac{1}{\nu}\|\nabla F\|_{L^{2}(0,T;L^{2})}^{2}\big)\\
\leq & M.
\end{aligned}
\end{equation}
Multiplying \eqref{eq:3.48} by $t$ gives that
\begin{equation*}
\begin{aligned}
\frac{\mathrm{d}}{\mathrm{d}t}&\Big( t\big(\frac{M_5}{\mu}\|\sqrt{\rho }\dot{\u}\|_2^{2}+\|\nabla\d_t\|_2^2\big) +\frac{2M_5}{\mu}\frac{\nu - \mu}{\nu} t\int \widetilde{F} \partial_j \u \cdot \nabla \u^j \;\mathrm{d}x \Big)\\
&+ t\big(\frac{M_5}{\mu}\frac{\nu - \mu}{\nu^2} \|\dot{\widetilde{F}}\|_2^2+ \frac{M_5}{2}\|\nabla \dot{\u}\|_2^2+\|\d_{tt}\|_2^2+\frac{1}{4}\|\nabla^2\d_t\|_2^2\big)\\
\leq &\big(\frac{M_5}{\mu}\|\sqrt{\rho }\dot{\u}\|_2^{2}+\|\nabla\d_t\|_2^2\big)+\frac{2M_5}{\mu}\frac{\nu - \mu}{\nu} \Big|\int \widetilde{F} \partial_j \u \cdot \nabla \u^j \;\mathrm{d}x\Big|\\ 
&+Mt\big(1+ \frac{1}{\nu^2} \|\widetilde{F}\|_4^4 + \frac{1}{\nu^2} \|\widetilde{P}\|_4^4 +\|\nabla \u\|_4^4+\|\Dv\u\|_4^4+ \| \nabla^3 \d\|_2^2 +\|\sqrt{\rho}\dot{\u}\|_2^2+\|\nabla\d_t\|_2^2\big).    
\end{aligned}    
\end{equation*}
Integrating the above inequality over $(0, T)$ and applying \eqref{2.5}, \eqref{3.11}, \eqref{eq:3.49} and \eqref{2.688} yields that
\begin{equation}\label{Equ3}
 \begin{aligned}
&\sup_{0\leq t\leq T} t\big(\|\sqrt{\rho} \dot{\u}\|_2^2+\|\nabla\d_t\|_2^2\big)\\
&\quad+\int_0^Tt\big(\|\nabla \dot{\u}\|_2^2 + \frac{\nu - \mu}{\nu^2} \|\dot{\widetilde{F}}\|_2^2+\|\d_{tt}\|_2^2+\|\nabla^2\d_t\|_2^2\big)\,\mathrm{d}t\leq  M(T),
 \end{aligned}   
\end{equation}
where $M(T)$ depends on $a$, $\mu$, $\gamma$ , $C_0$, $N_0$, $T$, $\overline{\omega}$ and $\|\rho_0\|_\infty$ and is independent of $\nu$.
Combining \eqref{Equ3} and \eqref{eq:3.47}  yields  $\sup_{0\leq t\leq T}t\|\nabla^3\d\|_2^2\leq M(T)$, which together with \eqref{Equ3} implies \eqref{2.66}, and hence we finish the proof of Proposition \ref{2-4}.
\end{proof}

\begin{Corollary}\label{2-5}
Assume that $(\rho, \u, \d)$ is a smooth solution of \eqref{1.1}--\eqref{1.3} on $[0, T ] \times \mathbb{T}^2$ fulfilling \eqref{Sup}, with $\nu\ge\nu_0$, where $\nu_0$ is determined
in Proposition \ref{2-3}. For any $\varepsilon\in(0,1)$, we obtain \footnote{Estimates \eqref{2.86} together with the identity $\operatorname{curl}\operatorname{curl}\u=-\Delta \u+\nabla\operatorname{div}\u$ yield that $\operatorname{curl}\u \in L^r_{loc}(\mathbb{R}^+;L^\infty)$ with $1<r<2$.}
\begin{equation}\label{2.86}
\|\Dv \u\|_{L^{(2-\varepsilon )}(0,T;L^{\infty })}^{2-\varepsilon}\le M(T)\nu ^{2(\varepsilon-1)},\;\|\nabla \mathbb{P}\u\|_{L^{(2-\varepsilon )}(0,T;L^{\infty })}^{2-\varepsilon}\le M(T),
\end{equation}
where $M(T)$ depends on $a$, $\mu$, $\gamma$ , $C_0$, $N_0$, $T$, $\varepsilon$, $\overline{\omega}$ and $\|\rho_0\|_\infty$ and is independent of $\nu$.
\end{Corollary}

\begin{proof}
Fix $\varepsilon\in(0,1)$ and set $p=2+\varepsilon\in(2,3)$.
The momentum equation gives
\begin{equation}\label{D-3}
\nabla(\sqrt t F)=\mathbb Q
 \big(\sqrt t\,\rho\dot{\u}+\sqrt t\,\nabla\d\Delta\d\big),
\qquad
\mu\Delta(\sqrt t\,\mathbb P\u)=\mathbb P
 \big(\sqrt t\,\rho\dot{\u}+\sqrt t\,\nabla\d\Delta\d\big),
\end{equation}
where $\mathbb Q:=\mathbb I-\mathbb P$.
Since $\int\rho\dot{\u}\,\mathrm{d}x=0$, using \eqref{Sup}, \eqref{2.1.3*}
and \eqref{L-3}, we obtain
\begin{equation}\label{D-4}
\|\sqrt t\,\rho\dot{\u}\|_p
\leq M_p
\|\sqrt{\rho t}\dot{\u}\|_2^{(4-p)/p}
\|\sqrt t\,\nabla\dot{\u}\|_2^{2(p-2)/p}.
\end{equation}
We infer from \eqref{DEU-1} and \eqref{2.66} that
\begin{equation}\label{D-2}
\|\sqrt t\,\nabla\d\Delta\d\|_2
\leq M\|\sqrt t\,\nabla^3\d\|_2\leq M(T).
\end{equation}
Moreover, \eqref{L-1}, \eqref{1.11}, \eqref{2.5}
and \eqref{2.66} give
\begin{equation}\label{D-1}
\begin{aligned}
\int_0^T\|\sqrt t\,\nabla(\nabla\d\Delta\d)\|_2^2\,\mathrm{d}t
&\leq M\sup_{0<t\leq T}
 \big(t\|\nabla^3\d\|_2^2\big)
 \int_0^T\big(\|\nabla^2\d\|_2^2
+\|\nabla\d\|_\infty^2\big)\,\mathrm{d}t\\
&\leq M(T),
\end{aligned}
\end{equation}
where we used
$\|\nabla\d\|_\infty^2
\leq M\|\nabla\d\|_2\|\nabla^3\d\|_2$.
Combining \eqref{D-2} and \eqref{D-1} with \eqref{L-3},
and then using \eqref{D-3}, \eqref{D-4} and \eqref{2.66}, yields
\begin{equation}\label{2.88}
\int_0^T\Big(
 \|\nabla(\sqrt t F)\|_p^{p/(p-2)}
+\|\nabla^2(\sqrt t\,\mathbb P\u)\|_p^{p/(p-2)}
\Big)\,\mathrm{d}t\leq M_p(T),
\end{equation}
and
\begin{equation}\label{2.89}
\int_0^T\|\sqrt t\,\nabla\d\Delta\d\|_p^{2p/(p-2)}
\,\mathrm{d}t\leq M_p(T).
\end{equation}
In view of \eqref{L-1}, we have
\[
\begin{aligned}
\|\widetilde F\|_\infty
&\leq M_p
 \|\widetilde F\|_2^{(p-2)/(2p-2)}
 \|\nabla F\|_p^{p/(2p-2)},\\
\|\nabla\mathbb P\u\|_\infty
&\leq M_p
 \|\nabla\mathbb P\u\|_2^{(p-2)/(2p-2)}
 \|\nabla^2\mathbb P\u\|_p^{p/(2p-2)}.
\end{aligned}
\]
Since $p=2+\varepsilon$, the time integral required by
H\"older's inequality is finite:
\[
\int_0^T t^{-(4-\varepsilon^2)/(4+2\varepsilon^2)}
\,\mathrm{d}t<\infty.
\]
Hence \eqref{2.5}, \eqref{2.88}, H\"older's inequality
and $\nu\geq\nu_0>0$ imply
\[
\begin{aligned}
\int_0^T\|\widetilde F\|_\infty^{2-\varepsilon}\,\mathrm{d}t
&\leq M(T)\nu^{(2-\varepsilon)(p-2)/(4p-4)}
\leq M(T)\nu^\varepsilon,\\
\int_0^T\|\nabla\mathbb P\u\|_\infty^{2-\varepsilon}
\,\mathrm{d}t&\leq M(T).
\end{aligned}
\]
Finally, since $\nu\Dv\u=\widetilde F+\widetilde P$
and $\widetilde P$ is bounded,
\[
\int_0^T\|\Dv\u\|_\infty^{2-\varepsilon}\,\mathrm{d}t
\leq M(T)\nu^{\varepsilon-2}(\nu^\varepsilon+1)
\leq M(T)\nu^{2\varepsilon-2}.
\]
This proves \eqref{2.86}.
\end{proof}

The following corollary is proved by an argument similar to that used in the proof of Corollary 3.8 in \cite{L-M-X}.
\begin{Corollary}\label{2-7}
Let $(\rho, \u, \d)$ be a smooth solution of \eqref{1.1} on $[0, T ] \times \mathbb{T}^2$ fulfilling \eqref{Sup}. Let $\nu\ge\nu_0$, where $\nu_0$ is determined
in Proposition \ref{2-3}. Then, 
we have
\begin{equation*}
\int_0^T (1+|\log t|)\|\nabla\u \|_{BMO}\;\mathrm{d}t\leq M(T),
\end{equation*}
where $M(T)$ depends on $a$, $\mu$, $\gamma$ , $C_0$, $N_0$, $T$, $\overline{\omega}$ and $\|\rho_0\|_\infty$ and is independent of $\nu$.
\end{Corollary}
\begin{proof}
Remembering \eqref{2.27} and using the embedding $L^\infty(\mathbb{T}^2)\hookrightarrow BMO(\mathbb{T}^2)$ and the boundedness of Riesz operator in $BMO(\mathbb{T}^2)$, we gather that
\begin{equation}\label{2.90}
\begin{aligned}
\|\nabla\u \|_{BMO}&=\|\nabla\mathbb{P}\u-\nabla^2(-\Delta)^{-1}\Dv\u\|_{BMO}\\    &\leq M \big( \|\nabla\mathbb{P}\u\|_{BMO}+\|\Dv\u\|_{BMO}\big)\\
&\leq M \big( \|\nabla\mathbb{P}\u\|_{\infty}+\|\Dv\u\|_{\infty}\big).\\
\end{aligned}
\end{equation}
It is not difficult to see that for $1<\alpha<2$,
\begin{equation}\label{2.91}
\begin{aligned}
&\int_0^T (1+|\log t|)\|\nabla\u \|_{BMO}\;\mathrm{d}t\\
&=\int_0^T \|\nabla\u \|_{BMO}\;\mathrm{d}t+\int_0^T |\log t\cdot t^s|\,t^{-s}\|\nabla\u \|_{BMO}\;\mathrm{d}t\\
&\leq T^{\frac{\alpha-1}{\alpha}}\Big(\int_0^T \|\nabla\u \|_{BMO}^{\alpha}\;\mathrm{d}t\Big)^{\frac{1}{\alpha}}
+\sup_{0\leq t\leq T}|\log t\cdot t^s|\Big(\int_0^T t^{-{s}\frac{\alpha}{\alpha-1}}\;\mathrm{d}t\Big)^{\frac{\alpha-1}{\alpha}}\Big(\int_0^T\|\nabla\u \|_{BMO}^{\alpha}\;\mathrm{d}t\Big)^{\frac{1}{\alpha}}
\end{aligned}
\end{equation}
Then, for fix $s\in(0,\frac{\alpha-1}{\alpha})$ such that $0<{s}\frac{\alpha}{\alpha-1}<1$ (our condition on $s$ makes it possible), one infers
\begin{equation*}
  \lim_{t\rightarrow 0^+} t^{s}\log t=0\,\;\text{and}\,\;
\int_0^T t^{-{s}\frac{\alpha}{\alpha-1}}\;\mathrm{d}t\leq M(T).
\end{equation*}
Combining \eqref{2.86}, \eqref{2.90}, and \eqref{2.91} completes the proof of Corollary \ref{2-7}.
\end{proof}

With all the estimates in Section 3 at hand, we now prove Theorem \ref{th:main1}.
\medskip

\begin{proof}[Proof of Theorem \ref{th:main1}]
\noindent\textit{Step 1. Approximation of the initial data and common bounds.}
Let $j_\delta$ be a nonnegative periodic mollifier of unit mass.
For $0<\delta<1$, set
\begin{equation*}
\begin{aligned}
\rho_0^\delta&=(1-\delta)(j_\delta*\rho_0)+\delta,
&\d_0^\delta&=\frac{j_\delta*\d_0}{|j_\delta*\d_0|},\\
\u_0^\delta&=j_\delta*\u_0-\mathbf c_\delta,
&\mathbf c_\delta&=\int\rho_0^\delta(j_\delta*\u_0)\,\mathrm{d}x.
\end{aligned}
\end{equation*}
By construction,
\begin{equation*}
\begin{gathered}
\int\rho_0^\delta\,\mathrm{d}x=1,\qquad
\int\rho_0^\delta\u_0^\delta\,\mathrm{d}x=0,\qquad
\delta\leq\rho_0^\delta\leq\|\rho_0\|_\infty,\\
|\d_0^\delta|=1,\qquad(\d_0^\delta)_2\geq\epsilon_0.
\end{gathered}
\end{equation*}
The normalization of the director is legitimate since
$\epsilon_0\leq|j_\delta*\d_0|\leq1$.
Moreover, $\mathbf c_\delta\to0$ and
\begin{equation*}
\rho_0^\delta\to\rho_0\ \text{in }L^p\ (1\leq p<\infty),\qquad
\u_0^\delta\to\u_0\ \text{in }H^1,\qquad
\d_0^\delta\to\d_0\ \text{in }H^2.
\end{equation*}
Let $C_0^\delta$ and $N_0^\delta$ denote, respectively, the
energy in \eqref{initial} and the left-hand side of
\eqref{1.11*}, evaluated at the approximate data.
The preceding convergences imply $C_0^\delta\to C_0$ and
$\lim_{\delta\to0}N_0^\delta\leq N_0$. Hence, for all sufficiently
small $\delta$, the fixed quantities
\begin{equation*}
\widehat C_0:=C_0+1,\qquad
\widehat N_0:=N_0+1,\qquad
\widehat R:=\|\rho_0\|_\infty
\end{equation*}
are common upper bounds for $C_0^\delta$, $N_0^\delta$ and
$\|\rho_0^\delta\|_\infty$, respectively.
We use these bounds throughout the approximation argument.
In particular, we fix the common density bootstrap bound
\begin{equation*}
\widehat\rho^*
:=2\exp\!\left(\frac{\gamma-1}{a\gamma}\widehat C_0\right)
\widehat R
\end{equation*}
according to \eqref{2.65}, and choose the rigidity constant
in Lemma \ref{2-1-4} with $E=(1+2\widehat C_0)^{1/2}$.
We then choose $M_4$ in Proposition \ref{2-3} with these
common bounds and set
\begin{equation*}
\nu_0=\max\left\{\mu,\frac{M_4^2}{\log^2(3/2)}\right\}.
\end{equation*}
This is the threshold used in Theorem \ref{th:main1}.
All these constants depend only on the original data and the
fixed parameters, and are independent of $\delta$ and $\nu$.

\medskip
\noindent\textit{Step 2. Global approximate solutions and passage to the limit.}
Fix $\nu\geq\nu_0$. The local theory for smooth periodic data
with positive density provides a maximal classical solution
$(\rho^\delta,\u^\delta,\d^\delta)$, with lifespan $T_\delta$.
The maximum principle preserves $|\d^\delta|=1$ and
$(\d^\delta)_2\geq\epsilon_0$.
The estimates of Propositions \ref{2-2}--\ref{2-4} and
Corollaries \ref{2-5}--\ref{2-7}, applied with the common bounds
from Step 1, are uniform in $\delta$.
In particular, the continuity equation ensures that
$\rho^\delta$ stays bounded away from zero on finite time
intervals, for each fixed $\delta$.

Let us briefly justify that $T_\delta=\infty$.
If $T_\delta<\infty$, then, for $0<\tau<T_\delta$ and $p>2$,
\eqref{2.66}, \eqref{Dot} and Sobolev embedding give
\begin{equation*}
\begin{aligned}
&\|\nabla^2\u^\delta\|_p
\leq C\|\nabla\rho^\delta\|_p+g_\delta(t),
\qquad g_\delta\in L^2(\tau,T_\delta),
\\
&\d^\delta\in L^\infty(\tau,T_\delta;H^3).
\end{aligned}
\end{equation*}
Here $g_\delta$ is bounded by a fixed multiple of
$\|\rho^\delta\dot{\u}^\delta\|_p
+\|\nabla\d^\delta\Delta\d^\delta\|_p$.
The transport estimate for $\nabla\rho^\delta$, combined
with the logarithmic Sobolev inequality and the BMO bound
of Corollary \ref{2-7}, yields boundedness of
$\|\nabla\rho^\delta\|_p$ and integrability of
$\|\nabla\u^\delta\|_\infty$ up to $T_\delta$, as in
\cite[Section~4]{R.D}.
Parabolic regularity for the velocity and director equations
then allows continuation beyond $T_\delta$, a contradiction.
Thus the approximate solutions are global.

Fix $T>0$ and put $Q_T=(0,T)\times\mathbb T^2$.
By \eqref{Dot}, \eqref{2.688} and \eqref{Equ3}, for each
$0<\tau<T$,
\begin{equation*}
\|\u_t^\delta\|_{L^2(\tau,T;L^2)}
+\|F_t^\delta\|_{L^2(\tau,T;H^{-1})}
\leq C_{\tau,T,\nu},
\end{equation*}
uniformly in $\delta$. For the flux, we use
$F_t^\delta=\nu\Dv\u_t^\delta-P_t^\delta$ and
$P_t^\delta=-\Dv(P^\delta\u^\delta)
-h(\rho^\delta)\Dv\u^\delta$.
The director equation and \eqref{3.11} also bound
$\d_t^\delta$ in $L^2(0,T;H^1)$.
The Aubin--Lions lemma, followed by $\tau\downarrow0$ using
the uniform $L^\infty(0,T;L^2)$ bounds, therefore gives,
up to extraction,
\begin{equation*}
\begin{aligned}
&\rho^\delta\rightharpoonup^\ast\rho
&&\text{in }L^\infty(Q_T),\\
&(\u^\delta,F^\delta)\to(\u,F)
&&\text{in }L^2(Q_T),\\
&\d^\delta\to\d
&&\text{in }C([0,T];H^1)\cap L^2(0,T;H^2).
\end{aligned}
\end{equation*}
To identify the pressure, we use the strong convergence of
$F^\delta$ and the renormalized continuity equation.
Let $\mathcal R$ be the weak-$\ast$ limit of
$\rho^\delta\log\rho^\delta$.
As in \cite[Section~4]{R.D}, the monotonicity of $P$ shows
that the integral of the nonnegative entropy defect
$\mathcal R-\rho\log\rho$ cannot increase.
It vanishes initially; hence strict convexity yields
\begin{equation*}
\rho^\delta\to\rho,\qquad P(\rho^\delta)\to P(\rho)
\quad\text{strongly in }L^p(Q_T),\qquad 1\leq p<\infty.
\end{equation*}
The director convergence likewise gives strong convergence
of $\nabla\d^\delta\odot\nabla\d^\delta$ and
$|\nabla\d^\delta|^2\d^\delta$ in $L^2(Q_T)$.
We may thus pass to the limit in all the nonlinear terms.
The estimates pass to the limit by weak lower semicontinuity.
The renormalized continuity equation and the energy balance,
justified by time regularization, give the initial traces
and the time continuity in \eqref{1.15}.
Finally, a diagonal extraction on finite time intervals
provides the required global solution.
\end{proof}

\section{Proof of Theorem \ref{th:main2} }
Based on the result of global existence established in Theorem \ref{th:main1}, the goal of this section is to prove the uniqueness of the solution $(\rho, \u, \d)$ under additional assumptions on the pressure, i.e., $P=a\rho$ with $a>0$. Although we miss the regularity property that $\nabla\u\in L^1_{loc}(0,T;L^{\infty }(\mathbb{T}^{2}))$, and $\u$ may not have a Lipschitz flow, it is sufficient to prove the uniqueness of the solution given that $\nabla \u$ belongs to  $L^r_{loc}(0,T;BMO(\mathbb{T}^{2}))$ with $r\in(1,2)$. In fact, the uniqueness of $(\rho, \u, \d)$ in Theorem \ref{th:main2} is a direct result of the following proposition.


\begin{Proposition}\label{2-6}
Assume that $P(\rho)=a\rho$ for some $a > 0$ in \eqref{1.1}, and let $T_0\in (0,1/2)$. Let $(\rho_1, \u_1, \d_1)$ and $(\rho_2, \u_2, \d_2)$ be two finite energy
solutions  to \eqref{1.1} on $[0, T_0]\times \mathbb{T}^2$, 
which satisfy
\eqref{1.6} and share the same initial data. If, in addition\footnotemark \footnotetext{Note that $(\nabla \d_1, \nabla \d_2) \in L^4(0, T_0; L^\infty(\mathbb{T}^2))$ follows from \eqref{L-1} and \eqref{2.5}, $(\sqrt{t}\dot{\u}_1, \sqrt{t}\dot{\u}_2) \in L^2(0, T_0; L^4(\mathbb{T}^2))$ follows from \eqref{Dot} and \eqref{2.66}, and $(\nabla \u_1, \nabla \u_2) \in L^2(0, T_0; L^4(\mathbb{T}^2))$ follows from \eqref{L-1}, \eqref{1.11} and \eqref{2.5}.
}
\begin{equation}\label{3.1}
\left\{
\begin{aligned}
  &(\rho_1,\rho_2) \in L^{\infty }(0,T_0;L^\infty(\mathbb{T}^2))\cap C([0,T_0];L^{p}(\mathbb{T}^{2})),\;1<p<\infty ,
  \; 0\leq (\rho_1,\rho_2)\leq \rho^*, \\
  & (\u_1, \u_2, \nabla \d_1, \nabla \d_2) \in L^\infty(0, T_0; H^1(\mathbb{T}^2)),
    \quad (\sqrt{t}\nabla \dot{\u}_1, \sqrt{t}\nabla \dot{\u}_2,\nabla{\d_1}_t,\nabla{\d_2}_t) \in L^2([0, T_0]\times \mathbb{T}^2), \\
  & (\sqrt{t\rho_1}\dot{\u}_1, \sqrt{t\rho_2}\dot{\u}_2) \in L^\infty(0, T_0; L^2(\mathbb{T}^2)),
    \quad (\nabla \u_1, \nabla \u_2) \in L^2(0, T_0; L^4(\mathbb{T}^2)), \\
  & (\nabla \d_1, \nabla \d_2) \in L^4(0, T_0; L^\infty(\mathbb{T}^2)),\quad (\sqrt{t}\dot{\u}_1, \sqrt{t}
  \dot{\u}_2) \in L^2(0, T_0; L^4(\mathbb{T}^2)),
\end{aligned}
\right.
\end{equation}
then $(\rho_2, \u_2, \d_2)\equiv(\rho_1, \u_1, \d_1)$ on $[0, T_0]\times \mathbb{T}^2$.
\end{Proposition}

\begin{proof}
For simplicity, assume that $a = 1$, and consider the two solutions $(\rho_1, \u_1, \d_1)$ and  $(\rho_2, \u_2, \d_2)$ satisfying \eqref{3.1} with the same initial data $(\rho_0, \u_0, \d_0)$. The system for the difference
\begin{equation*}
\delta \rho =\rho_1 -\rho _{2},\delta \u=\u_1-\u_{2},\delta \d=\d_1-\d_{2},
\end{equation*}
reads
\begin{equation}\label{3.3}
\left\{ \begin{aligned}
  & \delta \rho _{t}+\Dv (\delta \rho \u_{2}+\rho_1 \delta \u)=0, \\
 & \rho_1 \delta \u_{t}+\rho_1 \u_1\cdot \nabla \delta \u-\mu \Delta \delta \u-(\mu +\lambda )\nabla \Dv \delta \u+\nabla \delta \rho =-\delta \rho \dot{\u}_{2}-\rho_1 \delta \u\cdot \nabla \u_{2} \\
 & \quad-\nabla\d_2 \cdot \Delta\delta\d -\nabla\delta\d\cdot \Delta\d_1, \\
 & {\delta\d}_{t}- \Delta \delta\d=\nabla\delta\d:(\nabla\d_2+\nabla\d_1)\d_1+|\nabla\d_2|^2\delta\d-\delta\u\cdot\nabla\d_2-\u_1\cdot\nabla\delta\d, \\
 & |\d_1|=|\d_2|=1. 
\end{aligned} \right.
\end{equation}

\noindent$\mathbf{Step I.}$\, \textit{A priori estimates
for $\delta\rho$}.  In order to show that $\delta \rho \equiv 0,\;\delta \u\equiv 0$ and $\delta \d\equiv 0$,  we shall estimate the $L^\infty(0,T_0;\dot{H}^{-1}(\mathbb{T}^2))$ norm of $\delta\rho$, and the $L^2([0,T_0]\times\mathbb{T}^2)$ norms of $\sqrt{\rho_1}\delta \u$, $\delta \d$ and $\nabla\delta \d$.
To this end, we set $\phi =-{(-\Delta )}^{-1}\delta \rho$. Then,
\begin{equation}\label{3.4}
\|\nabla \phi \|_{2}=\|\delta \rho \|_{\dot{H}^{-1}}.
\end{equation}
Testing the first equation of \eqref{3.3} by $\phi$ yields
\begin{equation}\label{3.4-1}
\frac{1}{2}\frac{\mathrm{d}}{\mathrm{d}t}\|\nabla \phi \|_2^{2}\le \Big| \int{(\u_{2}\cdot \nabla \phi \delta \rho +\rho_1 \delta \u\cdot \nabla \phi )\;\mathrm{d}x} \Big|.
\end{equation}
The last term on the right-hand side of the inequality above can be estimated as follows:
\begin{equation}\label{3.4-2}
\Big| \int\rho_1 \delta \u\cdot \nabla \phi \;\mathrm{d}x \Big|\le \sqrt{\rho^*}\|\sqrt{\rho_1} \delta \u\|_{2}\|\nabla \phi \|_{2}.
\end{equation}
As for the first term on the right-hand side of \eqref{3.4-1}. Integration by parts yields
\begin{equation*}
\int\u_{2}\cdot \nabla \phi \delta \rho \;\mathrm{d}x=\int\u_2^{j}\partial _{j}\phi \Delta \phi \;\mathrm{d}x=-\int\partial _{k}\u_2^{j}\partial _{j}\phi \partial _{k}\phi \;\mathrm{d}x+\frac{1}{2}\int{\Dv\u_{2}|\nabla \phi |^{2}\;\mathrm{d}x}.
\end{equation*}
Thus, by virtue of the duality between the Hardy $\mathcal{H}^1$ and $BMO$ spaces, we obtain
\begin{equation*}
\Big| \int\u_{2}\cdot \nabla \phi \delta \rho \;\mathrm{d}x \Big|\le c\|\nabla \u_{2}\|_{BMO}\|\nabla \phi \otimes \nabla \phi \|_{\mathcal{H}^{1}},
\end{equation*}
with a universal positive constant $c$.
Employing the following inequality (see, e.g., \cite{PBM})
\begin{equation*}
\|f\|_{\mathcal{H}^{1}}\le c\|f\|_{1}\big( |\log \|f\|_{1}|+\log (e+\|f\|_{\infty}) \big),
\end{equation*}
we get that
\begin{equation}\label{3.5}
\Big| \int\u_{2}\cdot \nabla \phi \delta \rho \;\mathrm{d}x \Big|\le c\|\nabla \u_{2}\|_{BMO}\|\nabla \phi \|_2^{2}\big( |\log \|\nabla \phi \|_2^{2}|+\log (e+\|\nabla \phi \|_{\infty}^{2}) \big).
\end{equation}
Since $\rho_1$ and $\rho_2$ are bounded, for all $t\in[0,T_0]$, we can obtain $\|\nabla \phi \|_{\infty }\le c\rho^*$.
Thus, inserting \eqref{3.4-2} and \eqref{3.5} into \eqref{3.4-1} implies that
\begin{equation*}
\frac{1}{2}\frac{\mathrm{d}}{\mathrm{d}t}\|\nabla \phi \|_2^{2}\le C\Big( \|\sqrt{\rho_1 }\delta \u\|_{2}+\|\nabla \u_{2}\|_{BMO}\|\nabla \phi \|_{2}\big( 1+\big|\log \|\nabla \phi \|_{2}\big| \big) \Big)\|\nabla \phi \|_{2}.
\end{equation*}
Here, $C$ denotes a positive constant depending on $\rho^*$. Due to the condition $\nabla \phi(0) = 0$, we obtain for all $t \in [0, T_0]$ 
\begin{equation*}
\|\nabla \phi \|_{2}\le C\Big( \int_0^{t_0}{\|\sqrt{\rho_1 }\delta \u\|_{2}}\;\mathrm{d}\tau +\int_0^{t_0}{\|\nabla \u_{2}\|_{BMO}\|\nabla \phi \|_{2}\big( 1+|\log \|\nabla \phi \|_{2} \big)\;\mathrm{d}\tau } \Big).
\end{equation*}
From now on, fix $t_0\in(0,T_0]$ and let $t\in(0,t_0]$. Define
$ Z(t):=\sup_{0<\tau\leq t}\tau^{-1/2} \|\delta\rho(\tau)\|_{\dot H^{-1}}$.
Then, \eqref{3.4} together with the Cauchy-Schwarz inequality yields
\begin{equation}\label{3.7}
Z(t_0)\le C\Big( \| \sqrt{\rho_1}\delta \u\|_{L^{2}((0,t_0)\times \mathbb{T}^{2})}+\int_0^{t_0}{\|\nabla \u_{2}(\tau)\|_{BMO}Z(\tau)\big( 1+|\log \tau |+|\log Z(\tau)| \big)\;\mathrm{d}\tau } \Big).
\end{equation}
\smallbreak

In what follows, we use a duality argument to control the difference $\sqrt{\rho_1}\delta \u$, $\delta \d$ and $\nabla\delta \d$ in $L^{2}((0,t_0)\times \mathbb{T}^{2})$. Let  $(\mathbf{w}(t,x),\mathbf{m}(t,x))$  be the solution to the following backward parabolic system:
\begin{equation}\label{3.8}
\left\{ \begin{aligned}
  & \rho_1 \w_{t}+\rho_1 \u_1\cdot \nabla \w+\mu \Delta \w+(\mu +\lambda )\nabla \Dv \w=-\rho_1\delta \u -\nabla\d_2\cdot\Delta\m+\nabla\d_2\cdot\m,\\
 & \m_{t}+ \Delta \m=-\delta \d+\w\cdot\nabla\d_2, \\
  & \w|_{t=t_0}=\m|_{t=t_0}=0,
\end{aligned} \right.
\end{equation}
where $\rho_1$, $\rho_2$,  $\u_1$, $\u_2$, $\d_1$, $\d_2$ satisfy \eqref{3.1}, and $\delta\u=\u_1-\u_2$ and $\delta\d=\d_1-\d_2$ are given by \eqref{3.3}. We emphasize that the above system fails for classical linear parabolic theory when the coefficients
are rough and may vanish. However, if $\rho_1$, $\u_1$ and $\d_2$ are smooth and $\rho_1$ is bounded away from zero, the
well-posedness is standard. For this system, the existence theory can be achieved by a regularizing
procedure together with the uniform estimates \eqref{3.9}.

\smallbreak

\noindent$\mathbf{Step II.}$\, \textit{Uniform estimates
for $(\w,\m)$}. Note that the mean of $\w$ need not be constant in time, we set
\[
\overline{\w}(t):=\int\w(t,x)\,\mathrm{d}x,
\qquad
N_1:=\int\rho_1(t,x)\,\mathrm{d}x=1.
\]
Then,
\begin{equation*}
N_1\overline{\w}=\int\rho_1\w\,\mathrm{d}x+\int(N_1-\rho_1)(\w-\overline{\w})\,\mathrm{d}x,
\end{equation*}
from this, Cauchy-Schwarz and Poincar\'{e}'s inequalities, we can obtain
\begin{equation*}
\begin{aligned}
N_1|\overline{\w}| &\leq\Big|\int\rho_1 \w \;\mathrm{d}x\Big|+\|N_1-\rho_1\|_{2}\|\w-\overline{\w}\|_2\\
&\leq C\big(\|\sqrt{\rho_1} \w\|_2 +\|\nabla\w\|_2\big).
\end{aligned}
\end{equation*}
Here, the constant $C$ is defined  as before. In particular,
\begin{equation}\label{3.16}
|\overline{\w}(t)|
\leq C\big(\|\sqrt{\rho_1}\w\|_2
   +\|\nabla\w\|_2\big).
\end{equation}
Applying the embeddings \(H^1(\mathbb T^2)\hookrightarrow L^4(\mathbb T^2)\) and \(H^2(\mathbb T^2)\hookrightarrow L^\infty(\mathbb T^2)\) to \(\mathbf w-\overline{\mathbf w}\), and using \eqref{3.16} and Poincaré’s inequality, we obtain
\begin{equation}\label{step2-sobolev}
\begin{aligned}
\|\w\|_4
&\leq C\big(
  \|\sqrt{\rho_1}\w\|_2
 +\|\nabla\w\|_2\big),\\
\|\w\|_\infty
&\leq C\big(
  \|\sqrt{\rho_1}\w\|_2
 +\|\nabla\w\|_2
 +\|\nabla^2\w\|_2\big).
\end{aligned}
\end{equation}
Since the Helmholtz projection is orthogonal in \(L^2(\mathbb{T}^2)\) and commutes with spatial derivatives, we get
\begin{equation}\label{EQU-9}
\begin{aligned}
\|\nabla\w\|_2^2
&=\|\nabla\mathbb{P}\w\|_2^2
  +\|\Dv\w\|_2^2,\\
\|\nabla^2\w\|_2^2
&=\|\nabla^2\mathbb{P}\w\|_2^2
  +\|\nabla\Dv\w\|_2^2.
\end{aligned}
\end{equation}

First, testing $\eqref{3.8}_1$ by $\w$ and $\eqref{3.8}_2$
by $\m-\Delta\m$, respectively, summing up and integrating the resulting equation over $\mathbb{T}^2$, we infer that
\begin{equation}\label{step2-low-identity}
\begin{aligned}
&-\frac12\frac{\mathrm{d}}{\mathrm{d}t}\big(
 \|\sqrt{\rho_1}\w\|_2^2
+\|\m\|_2^2
+\|\nabla\m\|_2^2\big)\\
&\quad+\mu\|\nabla\mathbb{P}\w\|_2^2
+\nu\|\Dv\w\|_2^2
+\|\nabla\m\|_2^2
+\|\nabla^2\m\|_2^2\\
&=\int\rho_1\delta\u\cdot\w
 +\delta\d\cdot\m
 +\nabla\delta\d:\nabla\m\,\mathrm{d}x-2I(t),
\end{aligned}
\end{equation}
where
\begin{equation}\label{step2-cross-term}
\begin{aligned}
I(t):=&\int
 (\w\cdot\nabla\d_2)\cdot(\m-\Delta\m)\,\mathrm{d}x=\int
 (\w\cdot\nabla\d_2)\cdot\m
 +\nabla(\w\cdot\nabla\d_2):\nabla\m\,\mathrm{d}x.
\end{aligned}
\end{equation}
For $I(t)$, we use \eqref{3.1} and \eqref{step2-sobolev}  to obtain
\begin{align}
|I(t)|\leq{}&
 \|\nabla\d_2\|_4
 \|\w\|_4\|\m\|_2+\big(\|\nabla\d_2\|_4
 \|\nabla\w\|_2+\|\nabla^2\d_2\|_2
 \|\w\|_4\big)
 \|\nabla\m\|_4\label{It}\\
 \leq{}&\varepsilon\big(
 \|\nabla\w\|_2^2
+\|\nabla^2\m\|_2^2\big)+M_0\big(
 \|\sqrt{\rho_1}\w\|_2^2
+\|\m\|_2^2
+\|\nabla\m\|_2^2\big)\notag,
\end{align}
where $M_0$ is a positive constant depending only on $\varepsilon$  and on the norms of
solutions $(\rho_1,\u_1,\d_1)$ and $(\rho_2,\u_2,\d_2)$. Using Cauchy-Schwarz inequality, it holds
\begin{equation}\label{Cauchy}
\begin{aligned}
&\Big|\int\rho_1\delta\u\cdot\w
 +\delta\d\cdot\m
 +\nabla\delta\d:\nabla\m\,\mathrm{d}x\Big|\\
 \leq & c\big(
 \|\sqrt{\rho_1}\delta\u\|_2^2
+\|\delta\d\|_2^2
+\|\nabla\delta\d\|_2^2+
 \|\sqrt{\rho_1}\w\|_2^2
+\|\m\|_2^2
+\|\nabla\m\|_2^2\big).
\end{aligned}
\end{equation}
Taking $\varepsilon$ sufficiently small
and putting \eqref{It} and \eqref{Cauchy} into \eqref{step2-low-identity}, we deduce
that, for a fixed constant $M_1>0$,
\begin{equation}\label{2}
\begin{aligned}
&-\frac{\mathrm{d}}{\mathrm{d}t}\big(
 \|\sqrt{\rho_1}\w\|_2^2
+\|\m\|_2^2
+\|\nabla\m\|_2^2\big)\\
&\quad+M_1\big(
 \|\nabla\w\|_2^2
+\|\nabla\m\|_2^2
+\|\nabla^2\m\|_2^2\big)\\
&\leq \widetilde{M}\big(
 \|\sqrt{\rho_1}\delta\u\|_2^2
+\|\delta\d\|_2^2
+\|\nabla\delta\d\|_2^2+
 \|\sqrt{\rho_1}\w\|_2^2
+\|\m\|_2^2
+\|\nabla\m\|_2^2\big),
\end{aligned}
\end{equation}
where $\widetilde{M}$ is a positive constant depending only on $\varepsilon$, $\mu$, $\lambda$, $\gamma$  and on the norms of
solutions $(\rho_1,\u_1,\d_1)$ and $(\rho_2,\u_2,\d_2)$.

Next, applying the elliptic estimates to  $\eqref{3.8}_1$ gives
\begin{align*}
\|\nabla^2\w\|_2^2
\leq C\big(&
 \|\sqrt{\rho_1}\w_t\|_2^2
+\|\sqrt{\rho_1}\delta\u\|_2^2+\|\u_1\|_4^2
  \|\nabla\w\|_4^2+\|\nabla\d_2\|_\infty^2
  (\|\m\|_2^2
  +\|\nabla^2\m\|_2^2)\big).
\end{align*}
Since $\nabla\w$ has zero spatial mean,
\[
\|\nabla\w\|_4^2
\leq c\|\nabla\w\|_2
       \|\nabla^2\w\|_2.
\]
Hence,
\begin{equation}\label{3.8*}
\begin{aligned}
\|\nabla^2\w\|_2^2
\leq{}& M_2\big(
 \|\sqrt{\rho_1}\w_t\|_2^2
+\|\sqrt{\rho_1}\delta\u\|_2^2+\|\u_1\|_4^4
\|\nabla\w\|_2^2+\|\nabla\d_2\|_\infty^2
 (\|\m\|_2^2
     +\|\nabla^2\m\|_2^2)\big),
\end{aligned}
\end{equation}
where $M_2$ depends only on $\mu$, $\lambda$ and $\rho^*$.
For the higher energy estimate, testing
$\eqref{3.8}_1$ by $\w_t$ and integrating over $\mathbb{T}^2$
gives
\begin{equation}\label{step2-velocity-test}
\begin{aligned}
&-\frac12\frac{\mathrm{d}}{\mathrm{d}t}
 \big(\mu\|\nabla\mathbb{P}\w\|_2^2
      +\nu\|\Dv\w\|_2^2\big)
 +\|\sqrt{\rho_1}\w_t\|_2^2\\
&=-\int
 \big(\rho_1\delta\u+\rho_1\u_1\cdot\nabla\w\big)
 \cdot\w_t\,\mathrm{d}x+\int
 (\w_t\cdot\nabla\d_2)\cdot(\m-\Delta\m)\,\mathrm{d}x.
\end{aligned}
\end{equation}
To rewrite the last integral, differentiate \eqref{step2-cross-term} with respect to $t$:
\begin{equation*}
\begin{aligned}
\frac{\mathrm{d}}{\mathrm{d}t}I(t)
={}&\int
 (\w_t\cdot\nabla\d_2)\cdot(\m-\Delta\m)\,\mathrm{d}x+\int
 (\w\cdot\nabla\d_{2t})\cdot(\m-\Delta\m)\,\mathrm{d}x\\
&+\int
 (\w\cdot\nabla\d_2)\cdot\m_t\,\mathrm{d}x-\int
 (\w\cdot\nabla\d_2)\cdot\Delta\m_t\,\mathrm{d}x,
\end{aligned}
\end{equation*}
which together with \eqref{step2-velocity-test} yields
\begin{equation}\label{step2-velocity-rewritten}
\begin{aligned}
&-\frac12\frac{\mathrm{d}}{\mathrm{d}t}
 \big(\mu\|\nabla\mathbb{P}\w\|_2^2
      +\nu\|\Dv\w\|_2^2+2I(t)\big)
 +\|\sqrt{\rho_1}\w_t\|_2^2\\
&=-\int
 \big(\rho_1\delta\u+\rho_1\u_1\cdot\nabla\w\big)
 \cdot\w_t\,\mathrm{d}x-\int
 (\w\cdot\nabla\d_{2t})\cdot(\m-\Delta\m)\,\mathrm{d}x\\
&\quad-\int
 (\w\cdot\nabla\d_2)\cdot\m_t\,\mathrm{d}x+\int
 (\w\cdot\nabla\d_2)\cdot\Delta\m_t\,\mathrm{d}x.
\end{aligned}
\end{equation}
On the other hand, multiplying $\eqref{3.8}_2$ by $\m_t$ and $-\Delta\m_t$, respectively, and integrating
the resulting equations by parts over $\mathbb{T}^2$, we infer that
\begin{equation}\label{step2-director-test}
\begin{aligned}
&-\frac12\frac{\mathrm{d}}{\mathrm{d}t}
 \big(\|\nabla\m\|_2^2+\|\nabla^2\m\|_2^2\big)
 +\|\m_t\|_2^2+\|\nabla\m_t\|_2^2\\
&=-\int
 \big(\delta\d\cdot\m_t
     +\nabla\delta\d:\nabla\m_t\big)\,\mathrm{d}x+\int
 (\w\cdot\nabla\d_2)\cdot\m_t\,\mathrm{d}x-\int
 (\w\cdot\nabla\d_2)\cdot\Delta\m_t\,\mathrm{d}x.
\end{aligned}
\end{equation}
Combining \eqref{step2-velocity-rewritten} and \eqref{step2-director-test}, we arrive at
\begin{equation}\label{step2-high-identity}
\begin{aligned}
&-\frac12\frac{\mathrm{d}}{\mathrm{d}t}\big(
 \mu\|\nabla\mathbb{P}\w\|_2^2
+\nu\|\Dv\w\|_2^2
+\|\nabla\m\|_2^2
+\|\nabla^2\m\|_2^2+2I(t)\big)\\
&\quad+\|\sqrt{\rho_1}\w_t\|_2^2
+\|\m_t\|_2^2
+\|\nabla\m_t\|_2^2\\
&=-\int
 \big(\rho_1\delta\u+\rho_1\u_1\cdot\nabla\w\big)
 \cdot\w_t\,\mathrm{d}x-\int
 (\w\cdot\nabla\d_{2t})\cdot(\m-\Delta\m)\,\mathrm{d}x\\
&\quad-\int
 \big(\delta\d\cdot\m_t
     +\nabla\delta\d:\nabla\m_t\big)\,\mathrm{d}x.
\end{aligned}
\end{equation}
We estimate every term on the right-hand side of \eqref{step2-high-identity}.
For any $\varepsilon>0$, it is easy to see that
\begin{align*}
&\Big|\int
 \rho_1\delta\u\cdot\w_t
+\delta\d\cdot\m_t+\nabla\delta\d:\nabla\m_t\,\mathrm{d}x\Big|\\
\leq&\varepsilon\big(
 \|\sqrt{\rho_1}\w_t\|_2^2
+\|\m_t\|_2^2
+\|\nabla\m_t\|_2^2\big)+M_\varepsilon\big(
 \|\sqrt{\rho_1}\delta\u\|_2^2
+\|\delta\d\|_2^2
+\|\nabla\delta\d\|_2^2\big).
\end{align*}
For the second term on the right-hand side of \eqref{step2-high-identity}, one gets that
\begin{align*}
&\Big|\int
 \rho_1(\u_1\cdot\nabla\w)\cdot\w_t\,\mathrm{d}x\Big|\\
\leq& c\sqrt{\rho^*}
 \|\sqrt{\rho_1}\w_t\|_2
 \|\u_1\|_4
 \|\nabla\w\|_2^{1/2}
 \|\nabla^2\w\|_2^{1/2}\\
\leq&\varepsilon\big(
 \|\sqrt{\rho_1}\w_t\|_2^2
+\|\nabla^2\w\|_2^2\big)
+\widetilde{M}\|\u_1\|_4^4
 \|\nabla\w\|_2^2.
\end{align*}
For the last term in \eqref{step2-high-identity}, it follows from \eqref{step2-sobolev} that
\begin{equation*}
\begin{aligned}
&\Big|\int
 (\w\cdot\nabla\d_{2t})\cdot(\m-\Delta\m)\,\mathrm{d}x\Big|\\
\leq &\|\nabla\d_{2t}\|_2
 \|\w\|_\infty
 \big(\|\m\|_2
     +\|\nabla^2\m\|_2\big)\\
\leq&\varepsilon\|\nabla^2\w\|_2^2
+\widetilde{M}\big(1+\|\nabla\d_{2t}\|_2^2\big)\\
&\times\big(
 \|\sqrt{\rho_1}\w\|_2^2
+\|\nabla\w\|_2^2
+\|\m\|_2^2
+\|\nabla^2\m\|_2^2\big).
\end{aligned}
\end{equation*}
Inserting these estimates into \eqref{step2-high-identity},
adding $\frac{1}{4M_2}$ times \eqref{3.8*},
and taking $0<\varepsilon<\frac{1}{16(1+M_2)}$,
we obtain
\begin{equation}\label{step2-high-estimate}
\begin{aligned}
&-\frac{1}{2}\frac{\mathrm{d}}{\mathrm{d}t}\big(
 \mu\|\nabla\mathbb{P}\w\|_2^2
+\nu\|\Dv\w\|_2^2
+\|\nabla\m\|_2^2
+\|\nabla^2\m\|_2^2+2I(t)\big)\\
&\quad+\frac{1}{2}\big(
 \|\sqrt{\rho_1}\w_t\|_2^2
+\|\m_t\|_2^2
+\|\nabla\m_t\|_2^2\big)
+\frac{1}{8M_2}\|\nabla^2\w\|_2^2\\
&\leq \widetilde{M}\big(
 \|\sqrt{\rho_1}\delta\u\|_2^2
+\|\delta\d\|_2^2
+\|\nabla\delta\d\|_2^2\big)\\
&\quad+\widetilde{M}B(t)\big(
 \|\sqrt{\rho_1}\w\|_2^2
+\|\nabla\w\|_2^2
+\|\m\|_2^2
+\|\nabla\m\|_2^2
+\|\nabla^2\m\|_2^2\big),
\end{aligned}
\end{equation}
where
\begin{equation*}
B(t):=1+\|\u_1(t)\|_4^4
       +\|\nabla\d_2(t)\|_\infty^2
       +\|\nabla\d_{2t}(t)\|_2^2.
\end{equation*}
Indeed, \eqref{3.1} implies $B\in L^1(0,T_0)$.

We now define the modified energy, with a fixed sufficiently large
constant $\vartheta>0$, by
\begin{equation*}
\begin{aligned}
E(t):={}&\mu\|\nabla\mathbb{P}\w\|_2^2
+\nu\|\Dv\w\|_2^2
+\|\nabla\m\|_2^2
+\|\nabla^2\m\|_2^2+2I(t)\\
&+\vartheta\big(
 \|\sqrt{\rho_1}\w\|_2^2
+\|\m\|_2^2
+\|\nabla\m\|_2^2\big).
\end{aligned}
\end{equation*}
By \eqref{It}, we choose $\varepsilon$ sufficiently small so that $0<\varepsilon<\min\{\frac\mu4,\frac\nu4,\frac{1}{16(1+M_2)}\}$, and then choose $\vartheta$ so that $\vartheta>2M_0+1$. Using \eqref{EQU-9}, we obtain that
\begin{equation}\label{step2-energy-equivalence}
\begin{aligned}
&\xi_1\big(
 \|\sqrt{\rho_1}\w\|_2^2
+\|\nabla\w\|_2^2
+\|\m\|_2^2
+\|\nabla\m\|_2^2
+\|\nabla^2\m\|_2^2\big)\\
\leq &E(t)\\
\leq &\xi_2\big(
 \|\sqrt{\rho_1}\w\|_2^2
+\|\nabla\w\|_2^2
+\|\m\|_2^2
+\|\nabla\m\|_2^2
+\|\nabla^2\m\|_2^2\big),
\end{aligned}
\end{equation}
where the constants $\xi_1$ and $\xi_2$ depend on $\varepsilon$, $\mu$, $\lambda$, $\gamma$  and on the norms of
solutions $(\rho_1,\u_1,\d_1)$ and $(\rho_2,\u_2,\d_2)$.
Multiply \eqref{step2-high-estimate} by $2$ and \eqref{2} by $\vartheta$ and adding the resulting inequalities, we obtain
\begin{equation*}
\begin{aligned}
&-\frac{\mathrm{d}}{\mathrm{d}t}E(t)
+M_3\big(
 \|\sqrt{\rho_1}\w_t\|_2^2
+\|\m_t\|_2^2
+\|\nabla\m_t\|_2^2
+\|\nabla^2\w\|_2^2\\
&\hspace{38mm}
+\|\nabla\w\|_2^2
+\|\nabla\m\|_2^2
+\|\nabla^2\m\|_2^2\big)\\
&\leq \widetilde{M}\big(
 (\|\sqrt{\rho_1}\delta\u\|_2^2
+\|\delta\d\|_2^2
+\|\nabla\delta\d\|_2^2)+B(t)E(t)\big).
\end{aligned}
\end{equation*}
Since $\w_{|t=t_0}=\m_{|t=t_0}=0$,  Gr\"onwall's inequality gives, for
$0\leq t\leq t_0$,
\begin{equation*}
\begin{aligned}
E(t)\leq \widetilde{M}\int_t^{t_0}
 &\exp\!\Big(\widetilde{M}\int_t^s B(\tau)\,\mathrm{d}\tau\Big)\,\big(
 \|\sqrt{\rho_1(s)}\delta\u(s)\|_2^2
+\|\delta\d(s)\|_2^2
+\|\nabla\delta\d(s)\|_2^2\big)\,\mathrm{d}s,
\end{aligned}
\end{equation*}
which together with \eqref{step2-energy-equivalence} and $B\in L^1(0,T_0)$ for $0\leq t_0\leq T_0$ gives
\begin{equation}\label{step2-energy-bound}
\begin{aligned}
&\operatorname*{sup}_{0\leq t\leq t_0}\big(
 |\overline{\w}(t)|^2
+\|\w\|_2^2
+\|\sqrt{\rho_1}\w\|_2^2
+\|\nabla\mathbb{P}\w\|_2^2
+\|\Dv\w\|_2^2
+\|\m\|_2^2
+\|\nabla\m\|_2^2
+\|\nabla^2\m\|_2^2\big)\\
&\quad+\|\sqrt{\rho_1}\w_t\|_{L^2(0,t_0;L^2)}^2
+\|\m_t\|_{L^2(0,t_0;L^2)}^2
+\|\nabla\m_t\|_{L^2(0,t_0;L^2)}^2+\|\nabla^2\w\|_{L^2(0,t_0;L^2)}^2\\
&\quad+\|\nabla\w\|_{L^2(0,t_0;L^2)}^2
+\|\nabla\m\|_{L^2(0,t_0;L^2)}^2
+\|\nabla^2\m\|_{L^2(0,t_0;L^2)}^2\\
&\leq \widetilde{M}\big(
 \|\sqrt{\rho_1}\delta\u\|_{L^2(0,t_0;L^2)}^2
+\|\delta\d\|_{L^2(0,t_0;L^2)}^2
+\|\nabla\delta\d\|_{L^2(0,t_0;L^2)}^2\big).
\end{aligned}
\end{equation}
Moreover, applying $\nabla$ to $\eqref{3.8}_2$, we get that
\begin{equation*}
\nabla\Delta\m=-\nabla\m_t-\nabla\delta\d
                   +\nabla(\w\cdot\nabla\d_2),
\end{equation*}
which leads to
\begin{equation}\label{step2-third-derivatives}
\begin{aligned}
\|\nabla^3\m\|_2^2
\leq3\|\nabla\m_t\|_2^2
+3\|\nabla\delta\d\|_2^2+6\|\nabla\d_2\|_\infty^2
     \|\nabla\w\|_2^2+6\|\nabla^2\d_2\|_2^2
     \|\w\|_\infty^2.
\end{aligned}
\end{equation}
The last two terms in \eqref{step2-third-derivatives} are integrable. More explicitly,
\begin{align}
\|\nabla(\w\cdot\nabla\d_2)\|_{L^2(0,t_0;L^2)}^2
\leq & C\big(
 \|\nabla\d_2\|_{L^2(0,T_0;L^\infty)}^2
+T_0\|\nabla^2\d_2\|_{L^\infty(0,T_0;L^2)}^2\big)\notag\\
&\times\big(
\|\sqrt{\rho_1}\w\|_{L^\infty(0,t_0;L^2)}^2
+\|\nabla\w\|_{L^\infty(0,t_0;L^2)}^2\big)\label{nabla3}\\
&\quad+C\|\nabla^2\d_2\|_{L^\infty(0,T_0;L^2)}^2
\|\nabla^2\w\|_{L^2(0,t_0;L^2)}^2\notag\\
\leq & \widetilde{M}\big(
 \|\sqrt{\rho_1}\delta\u\|_{L^2(0,t_0;L^2)}^2
+\|\delta\d\|_{L^2(0,t_0;L^2)}^2
+\|\nabla\delta\d\|_{L^2(0,t_0;L^2)}^2\big)\notag.
\end{align}
Thus, combining \eqref{step2-energy-bound}, \eqref{step2-third-derivatives} and \eqref{nabla3}  yields that
\begin{equation}\label{3.9}
\begin{aligned}
&\operatorname*{sup}_{0\leq t\leq t_0}\big(
 |\overline{\w}(t)|^2
+\|\w\|_2^2
+\|\sqrt{\rho_1}\w\|_2^2
+\|\nabla\mathbb{P}\w\|_2^2
+\|\Dv\w\|_2^2
+\|\m\|_2^2
+\|\nabla\m\|_2^2
+\|\nabla^2\m\|_2^2\big)\\
&\quad+\|\sqrt{\rho_1}\w_t\|_{L^2(0,t_0;L^2)}^2
+\|\m_t\|_{L^2(0,t_0;L^2)}^2
+\|\nabla\m_t\|_{L^2(0,t_0;L^2)}^2+\|\nabla^2\w\|_{L^2(0,t_0;L^2)}^2\\
&\quad+\|\nabla\w\|_{L^2(0,t_0;L^2)}^2
+\|\nabla\m\|_{L^2(0,t_0;L^2)}^2
+\|\nabla^2\m\|_{L^2(0,t_0;L^2)}^2+\|\nabla^3\m\|_{L^2(0,t_0;L^2)}^2\\
&\leq \widetilde{M}\big(
 \|\sqrt{\rho_1}\delta\u\|_{L^2(0,t_0;L^2)}^2
+\|\delta\d\|_{L^2(0,t_0;L^2)}^2
+\|\nabla\delta\d\|_{L^2(0,t_0;L^2)}^2\big).
\end{aligned}
\end{equation}
\smallbreak

\noindent$\mathbf{Step III.}$\, \textit{The proof of uniqueness}. For convenience, we set $$X(t_0):=\big(\|\sqrt{\rho_1}\delta \u(t)\|_{L^{2}(0,t_0;L^{2})}^2+ \|\delta \d(t)\|_{L^{2}(0,t_0;L^{2})}^2+\|\nabla\delta \d(t)\|_{L^{2}(0,t_0;L^{2})}^2\big)$$ 
Multiplying $\eqref{3.3}_2$ by $\w$ and integrating the resulting equation over $[0,t_0]\times\mathbb{T}^2$, by use of integration by parts, we can obtain
\begin{equation*}
\begin{aligned}
&\int_0^{t_0}\int(\rho_1\delta\u\cdot\w)_t-\big((\rho_1 \w)_{t}+\Dv(\rho_1 \u_1\otimes\w)+\mu \Delta \w+(\mu +\lambda )\nabla \Dv \w\big)\cdot\delta\u\;\mathrm{d}x\;\mathrm{d}t\\
=&\int_0^{t_0}\int \big(-\nabla\delta\rho-\delta \rho \dot{\u}_{2}-\rho_1 \delta \u\cdot \nabla \u_{2} -\nabla\d_2\cdot\Delta\delta\d-\nabla\delta\d\cdot\Delta\d_1\big)\cdot\w\;\mathrm{d}x\;\mathrm{d}t.
\end{aligned}
\end{equation*}
Thus, we infer from $\eqref{3.8}_1$ that
\begin{equation}\label{3.13}
\begin{aligned}
& \int_0^{t_0}\|\sqrt{\rho_1}\delta \u\|_2^2\;\mathrm{d}t+\int_0^{t_0}\int\big(\delta\u\cdot\nabla\d_2\cdot\Delta\m-\delta\u\cdot\nabla\d_2\cdot\m\big)\,\mathrm{d}x\;\mathrm{d}t \\
=&\int_0^{t_0}\int\big(-\nabla\delta\rho-\delta \rho \dot{\u}_{2}-\rho_1 \delta \u\cdot \nabla \u_{2} -\nabla\d_2\cdot\Delta\delta\d-\nabla\delta\d\cdot\Delta\d_1\big)\cdot\w\;\mathrm{d}x\;\mathrm{d}t\\
=& \int_0^{t_0}{\int\delta \rho\, \Dv \w \;\mathrm{d}x\;\mathrm{d}t}-\int_0^{t_0}{\int\delta \rho \dot{\u}_{2}\cdot \w \;\mathrm{d}x\;\mathrm{d}t}-\int_0^{t_0}{\int\rho_1 (\delta \u\cdot \nabla \u_{2})\cdot \w \;\mathrm{d}x\;\mathrm{d}t} \\
 &+\int_0^{t_0}{\int{\big( \w\nabla^2\d_2\cdot\nabla\delta\d+\nabla\w\nabla\d_2\cdot\nabla\delta\d-\w\nabla\delta\d\cdot\Delta\d_1 \big)\;\mathrm{d}x\;\mathrm{d}t}}\\
 \triangleq& \sum\limits_{i=1}^{4}J_i,
\end{aligned}
\end{equation}
where we have used $\w \big|_{t=t_0}=\delta\u\big|_{t=0}=0$.
For $J_1$, using \eqref{3.9} gives that
\begin{equation}\label{3.14}
\begin{aligned}
|J_1|  & \le t_0\|t^{-1/2}\delta \rho \|_{L^{\infty }(0,t_0;\dot{H}^{-1})}\|\nabla \Dv\w\|_{L^{2}(0,t_0;L^{2})} \\
 & \le \widetilde{M} t_0\|t^{-1/2}\delta \rho \|_{L^{\infty }(0,t_0;\dot{H}^{-1})}X(t_0)^{1/2}. \\
\end{aligned}
\end{equation}
One can bound $J_2$ as follows:
\begin{equation}\label{3.15}
\begin{aligned}
  |J_2|
  &=\Big| \int_0^{t_0}{\int t^{-1/2}}\delta \rho t^{1/2}\dot{\u}_{2}\cdot \w \;\mathrm{d}x\;\mathrm{d}t \Big| \\
 & \le \|t^{-1/2}\delta \rho \|_{L^{\infty }(0,t_0;\dot{H}^{-1})}\|t^{1/2}\nabla (\dot{\u}_{2}\cdot \w)\|_{L^{1}(0,t_0;L^{2})} \\
 & \le \|t^{-1/2}\delta \rho \|_{L^{\infty }(0,t_0;\dot{H}^{-1})}\Big( \|t^{1/2}\nabla \dot{\u}_{2}\|_{L^{2}(0,t_0;L^{2})}\|\w\|_{L^{2}(0,t_0;L^{\infty })}
 \\&\quad +\|t^{1/2}\dot{\u}_{2}\|_{L^{2}(0,t_0;L^{4})}\|\nabla \w\|_{L^{2}(0,t_0;L^{4})} \Big). \\
\end{aligned}
\end{equation}
According to \eqref{L-2} and \eqref{3.9}, one gets
\begin{equation}\label{b3.19}
\begin{aligned}
\|\w\|_{L^{2}(0,t_0;L^{\infty })}\leq& \|\w-\overline{\w}\|_{L^{2}(0,t_0;L^{\infty })}+{t_0}^{1/2}|\overline{\w}|\\
\leq &c \big({t_0}^{1/4}\|\nabla\w\|_{L^\infty(0,t_0;L^{2})}^{1/2}
\|\nabla^2\w\|_{L^{2}(0,t_0;L^{2})}^{1/2}+{t_0}^{1/2}|\overline{\w}|\big)\\
\le&c {t_0}^{1/4}\big(\|\nabla\w\|_{L^\infty(0,t_0;L^{2})}+
\|\nabla^2\w\|_{L^{2}(0,t_0;L^{2})}+|\overline{\w}|\big)\\
\leq &\widetilde{M} {t_0}^{1/4}X(t_0)^{1/2}, 
\end{aligned}
\end{equation}
Similarly,
\begin{equation}\label{3.20}
\begin{aligned}
\|\nabla \w\|_{L^{2}(0,t_0;L^{4})}&\le c {t_0}^{1/4}\|\nabla\w\|_{L^{\infty}(0,t_0;L^{2})}^{1/2}
\|\nabla^2\w\|_{L^{2}(0,t_0;L^{2})}^{1/2}\\ 
&\leq\widetilde{M}{t_0}^{1/4}X(t_0)^{1/2},
\end{aligned}
\end{equation}
By means of \eqref{3.1}, \eqref{b3.19} and \eqref{3.20}, the inequality \eqref{3.15} is transformed into
\begin{equation}\label{3.21}
\begin{aligned}
  |J_2|
  \le\widetilde{M}\|t^{-1/2}\delta \rho \|_{L^{\infty }(0,t_0;\dot{H}^{-1})}{t_0}^{1/4}X(t_0)^{1/2}.
\end{aligned}
\end{equation} 
For $J_3$, we can derive from \eqref{3.1}, \eqref{3.9} and Sobolev embedding inequality $H^{1}(\mathbb{T}^2)\hookrightarrow L^4(\mathbb{T}^2)$ that
\begin{equation}\label{3.22}
\begin{aligned}
 |J_3|
 & \le \sqrt{\rho ^{*}}\|\sqrt{\rho_1}\delta \u\|_{L^{2}(0,t_0;L^{2})}\|\nabla \u_{2}\|_{L^{2}(0,t_0;L^{4})}\|\w\|_{L^{\infty }(0,t_0;L^{4})}\\
 & \le C_1\|\sqrt{\rho_1}\delta \u\|_{L^{2}(0,t_0;L^{2})}\|\w\|_{L^{\infty }(0,t_0;H^{1})}\\
 &\leq C_1X(t_0),
\end{aligned}
\end{equation}
where $C_1$ depends on $\widetilde{M}$ and $\|\nabla \u_{2}\|_{L^{2}(0,t_0;L^{4})}$, and satisfies $\lim_{t_0\to 0^+}C_1(t_0)=0$.
To bound $J_4$, we use \eqref{3.1}, \eqref{b3.19}  and  \eqref{3.20} to obtain
\begin{equation}\label{3.23}
\begin{aligned}
|J_4|
&\leq \int_0^{t_0} \|\nabla\delta\d\|_2\big(\|\nabla^2\d_1\|_2+\|\nabla^2\d_2\|_2\big)\|\w\|_\infty+\|\nabla\delta\d\|_2\|\nabla\d_2\|_4\|\nabla\w\|_4\;\mathrm{d}t\\
&\leq \|\nabla\delta\d\|_{L^{2}(0,t_0;L^{2})}\big(\|\nabla^2\d_1\|_{L^{\infty}(0,t_0;L^{2})}+\|\nabla^2\d_2\|_{L^{\infty}(0,t_0;L^{2})}+\|\nabla\d_2\|_{L^{\infty}(0,t_0;L^{4})}\big)\\
&\quad\times\big( \|\w\|_{L^2(0,t_0;L^\infty)}+\|\nabla \w\|_{L^{2}(0,t_0;L^{4})}\big) \\
&\leq \widetilde{M}\|\nabla\delta\d\|_{L^{2}(0,t_0;L^{2})}\big( \|\w\|_{L^2(0,t_0;L^\infty)}+\|\nabla \w\|_{L^{2}(0,t_0;L^{4})}\big) \\
&\leq\widetilde{M}{t_0}^{1/4}X(t_0).
\end{aligned}
\end{equation}
Substitute \eqref{3.14}, \eqref{3.21}-- \eqref{3.23} into \eqref{3.13} to obtain (for $t_0<{t_0}^{1/4}$)
\begin{equation}\label{3.24}
\begin{aligned}
  & \int_0^{t_0}\|\sqrt{\rho_1}\delta \u\|_2^2\;\mathrm{d}t+\int_0^{t_0}\int\big(\delta\u\cdot\nabla\d_2\cdot\Delta\m-\delta\u\cdot\nabla\d_2\cdot\m\big)\,\mathrm{d}x\;\mathrm{d}t \\
\le &  \widetilde{M}{t_0}^{1/4}\big(\|t^{-1/2}\delta \rho \|_{L^{\infty }(0,t_0;\dot{H}^{-1})}X(t_0)^{1/2}+X(t_0)\big)+ C_1X(t_0).
\end{aligned}
\end{equation}

On the other hand, multiplying $\eqref{3.3}_3$ by $\m$ and integrating the resultant over $[0,t_0]\times \mathbb{T}^2$, we obtain that
\begin{equation*}
\begin{aligned}
&\int_0^{t_0}\int(\delta\d\cdot\m)_{t}-\big(\m_t+\Delta\m)\cdot
\delta\d\;\mathrm{d}x\;\mathrm{d}t\\
&=\int_0^{t_0}\int\big(\nabla\delta\d:(\nabla\d_2+\nabla\d_1)\d_1+|\nabla\d_2|^2\delta\d-\delta\u\cdot\nabla\d_2-\u_1\cdot\nabla\delta\d\big)\cdot\m \;\mathrm{d}x\;\mathrm{d}t.
\end{aligned}
\end{equation*}
Then, by applying the operator $\nabla$ to $\eqref{3.3}_3$, multiplying the resulting equation by $\nabla\m$ and using integration by parts over $[0,t_0] \times\mathbb{T}^2 $, we get 
\begin{equation*}
\begin{aligned}
&\int_0^{t_0}\int(\nabla\delta\d\cdot\nabla\m)_{t}-\big({\nabla\m}_t+\nabla\Delta\m)
\cdot\nabla\delta\d\;\mathrm{d}x\;\mathrm{d}t\\
&=\int_0^{t_0}\int\nabla\big(\nabla\delta\d:(\nabla\d_2+\nabla\d_1)\d_1+|\nabla\d_2|^2\delta\d-\delta\u\cdot\nabla\d_2-\u_1\cdot\nabla\delta\d\big)\cdot\nabla\m \;\mathrm{d}x\;\mathrm{d}t\\
&=\int_0^{t_0}\int\big(-\nabla\delta\d:(\nabla\d_2+\nabla\d_1)\d_1\cdot\Delta\m-|\nabla\d_2|^2\delta\d\cdot\Delta\m+\delta\u\cdot\nabla\d_2\cdot\Delta\m\\
&\quad+\u_1\cdot\nabla\delta\d\cdot\Delta\m\big)\;\mathrm{d}x\;\mathrm{d}t.
\end{aligned}
\end{equation*}
Combining the two equalities above and utilizing the fact that  $\m\big|_{t=t_0}=\delta\d\big|_{t=0}=0$, it follows from $\eqref{3.8}_2$ that 
\begin{equation*}
\begin{aligned}
& \int_0^{t_0}\int-\big(\m_t+ \Delta\m\big)\cdot\delta\d-\big({\nabla\m}_t+ \nabla\Delta\m\big)\cdot\nabla\delta\d\;\mathrm{d}x\;\mathrm{d}t\\
=& \int_0^{t_0}\big(\|\delta \d\|_2^{2}+\|\nabla\delta \d\|_2^{2}\big)\;\mathrm{d}t-\int_0^{t_0}\int\big(\w\cdot\nabla\d_2\cdot\delta\d+\nabla(\w\cdot\nabla\d_2)\nabla\delta\d\big)\,\mathrm{d}x\;\mathrm{d}t\\ 
=&\int_0^{t_0}\int\big(\nabla\delta\d:(\nabla\d_2+\nabla\d_1)\d_1+|\nabla\d_2|^2\delta\d-\delta\u\cdot\nabla\d_2-\u_1\cdot\nabla\delta\d\big)\cdot\m \;\mathrm{d}x\;\mathrm{d}t\\
&+\int_0^{t_0}\int-\nabla\delta\d:(\nabla\d_2+\nabla\d_1)\d_1\cdot\Delta\m-|\nabla\d_2|^2\delta\d\cdot\Delta\m+(\delta\u\cdot\nabla\d_2+\u_1\cdot\nabla\delta\d)\cdot\Delta\m \,\mathrm{d}x\;\mathrm{d}t,
\end{aligned}
\end{equation*}
which further implies that 
\begin{equation}\label{3.25}
\begin{aligned}
& \int_0^{t_0}\big(\|\delta \d\|_2^{2}+\|\nabla\delta \d\|_2^{2}\big)\;\mathrm{d}t-\int_0^{t_0}\int\big(\delta\u\cdot\nabla\d_2\cdot\Delta\m-\delta\u\cdot\nabla\d_2\cdot\m\big)\,\mathrm{d}x\;\mathrm{d}t\\
\leq&\int_0^{t_0}\int\big(|\nabla\delta\d||\nabla\d_2+\nabla\d_1|+|\nabla\d_2|^2|\delta\d|+|\u_1||\nabla\delta\d|\big)|\m|+|\w||\nabla\d_2||\delta\d| \;\mathrm{d}x\;\mathrm{d}t\\ 
&+\int_0^{t_0}\int\big(|\nabla\delta\d||\nabla\d_2+\nabla\d_1||\Delta\m|+|\nabla\d_2|^2|\delta\d||\Delta\m|+|\u_1||\nabla\delta\d||\Delta\m|\\
&+|\nabla(\w\cdot\nabla\d_2)||\nabla\delta\d|\big)\;\mathrm{d}x\;\mathrm{d}t\\
\leq&\|\nabla\delta\d\|_{L^{2}(0,t_0;L^{2})}\|(\nabla \d_1+\nabla\d_2)\|_{L^{2}(0,t_0;L^{\infty})}\big(\|\m\|_{L^{\infty}(0,t_0;L^{2})}+\|\nabla^2\m\|_{L^{\infty}(0,t_0;L^{2})}\big) \\
&+\|\nabla\d_2\|_{L^{4}(0,t_0;L^{\infty})}^{2}\|\delta\d\|_{L^{2}(0,t_0;L^{2})}\big(\|\m\|_{L^{\infty}(0,t_0;L^{2})}+\|\nabla^2\m\|_{L^{\infty}(0,t_0;L^{2})}\big)\\
&+\|\u_1\|_{L^{\infty}(0,t_0;L^{4})}\|\nabla\delta\d\|_{L^{2}(0,t_0;L^{2})}\big({t_0}^{1/2}\|\m\|_{L^{\infty}(0,t_0;L^{4})}+c\|\nabla^2\m\|_{L^{2}(0,t_0;L^{4})}\big)\\
&+\|\w\|_{L^{2}(0,t_0;L^{\infty})}\|\nabla\d_2\|_{L^{\infty}(0,t_0;L^{2})}\|\delta\d\|_{L^{2}(0,t_0;L^{2})}\\
&+\big(\|\nabla\w\|_{L^{2}(0,t_0;L^{4})}\|\nabla\d_2\|_{L^{\infty}(0,t_0;L^{4})}+\|\w\|_{L^{2}(0,t_0;L^{\infty})}\|\nabla^2\d_2\|_{L^{\infty}(0,t_0;L^{2})}\big)\|\nabla\delta\d\|_{L^{2}(0,t_0;L^{2})}\\
=&\sum_{i=1}^{5}LL_i.
\end{aligned}
\end{equation}
Now we estimate the terms on the right-hand side of \eqref{3.25}.  The estimate for $LL_1$ and $LL_2$ are straightforward, by \eqref{3.1} and \eqref{3.9}, there exists a positive
constant $C_2$ depends on  $\widetilde{M}$,
$\|\nabla \d_1\|_{L^{2}(0,t_0;L^{\infty})}$, $\|\nabla\d_2\|_{L^{2}(0,t_0;L^{\infty})}$, $\|\nabla\d_2\|_{L^{4}(0,t_0;L^{\infty})}$ and $\varepsilon$, and satisfies $\lim_{t_0\to 0^+}C_2(t_0)=0$ such that 
\begin{equation*}
\begin{aligned}
&LL_1+LL_2\\
\leq  &  \varepsilon\big(\|\delta\d\|_{L^{2}(0,t_0;L^{2})}^{2}+\|\nabla\delta\d\|_{L^{2}(0,t_0;L^{2})}^{2}\big)+C_2\big(\|\m\|_{L^{\infty}(0,t_0;L^{2})}^2+\|\nabla^2\m\|_{L^{\infty}(0,t_0;L^{2})}^2\big)\\
\leq &(\varepsilon+C_2) X(t_0).
\end{aligned}    
\end{equation*}
By \eqref{3.9} and the Sobolev embedding $H^1(\mathbb T^2)\hookrightarrow L^4(\mathbb T^2)$, we obtain
\begin{equation}\label{2.0}
 t_0\|\mathbf m\|_{L^\infty(0,t_0;L^4)}^2
\leq c t_0\|\mathbf m\|_{L^\infty(0,t_0;H^1)}^2
\leq \widetilde M t_0 X(t_0).   
\end{equation}
For $LL_3$, $LL_4$ and $LL_5$, in view of \eqref{L-1}, \eqref{3.1}, \eqref{3.9}, \eqref{b3.19}, \eqref{2.0} and Poincar\'{e}'s inequality, it is clear that (for any $t_0<1$)
\begin{equation*}
\begin{aligned}
 &LL_3+LL_4+LL_5\\
 \leq & \varepsilon\big(\|\delta\d\|_{L^{2}(0,t_0;L^{2})}^{2}+\|\nabla\delta\d\|_{L^{2}(0,t_0;L^{2})}^{2}\big)+\widetilde{M}\big(t_0\|\m\|_{L^{\infty}(0,t_0;L^{4})}^2+{t_0}^{1/2}\|\nabla^2\m\|_{L^{\infty}(0,t_0;L^{2})}\|\nabla^3\m\|_{L^{2}(0,t_0;L^{2})}\\
 &+\|\w\|_{L^{2}(0,t_0;L^{\infty})}^2+{t_0}^{1/2}\|\nabla\w\|_{L^{\infty}(0,t_0;L^{2})}\|\nabla^2\w\|_{L^{2}(0,t_0;L^{2})}\big)\\
 \leq & \varepsilon\big(\|\delta\d\|_{L^{2}(0,t_0;L^{2})}^{2}+\|\nabla\delta\d\|_{L^{2}(0,t_0;L^{2})}^{2}\big)+\widetilde{M}\big(t_0X(t_0)+{t_0}^{1/2}\|\nabla^2\m\|_{L^{\infty}(0,t_0;L^{2})}\|\nabla^3\m\|_{L^{2}(0,t_0;L^{2})}\\
 &+{t_0}^{1/2}X(t_0)+{t_0}^{1/2}\|\nabla\w\|_{L^{\infty}(0,t_0;L^{2})}\|\nabla^2\w\|_{L^{2}(0,t_0;L^{2})}\big)\\
 \leq & \big(\varepsilon+\widetilde{M}{t_0}^{1/2}\big)X(t_0).
\end{aligned}    
\end{equation*}
Putting the estimates of $LL_1$--$LL_5$ into \eqref{3.25}, it is easy to see that
\begin{equation}\label{4}
\begin{aligned}
& \int_0^{t_0}\big(\|\delta \d\|_2^{2}+\|\nabla\delta \d\|_2^{2}\big)\;\mathrm{d}t-\int_0^{t_0}\int\big(\delta\u\cdot\nabla\d_2\cdot\Delta\m-\delta\u\cdot\nabla\d_2\cdot\m\big)\,\mathrm{d}x\;\mathrm{d}t\\ 
 \leq& \big(2\varepsilon+\widetilde{M}{t_0}^{1/2}+C_2\big)X(t_0).
\end{aligned}    
\end{equation}
Clearly, for $t_0\in[0,T_0]$, invoking  \eqref{3.24}, \eqref{4} and the definition of $X(t_0)$, we infer that if $t_0$ and $\varepsilon$ are small enough that $\widetilde{M}{t_0}^{1/4}+C_{1}+2\varepsilon+\widetilde{M}{t_0}^{1/2}+C_2<\frac{1}{2}$, then
\begin{equation}\label{3.27}
\begin{aligned}
\frac{1}{2}X(t_0)
&\le \widetilde{M}{t_0}^{1/4}\|t^{-1/2}\delta \rho \|_{L^{\infty }(0,t_0;\dot{H}^{-1})}X(t_0)^{1/2},
\end{aligned}
\end{equation}
which together with Young's inequality implies that
\begin{equation}\label{3.28}
\begin{aligned}
\|\sqrt{\rho_1}\delta \u\|_{L^{2}(0,t_0;L^{2})}
&\le 2\widetilde{M}{t_0}^{1/4}\|t^{-1/2}\delta \rho \|_{L^{\infty }(0,t_0;\dot{H}^{-1})}.
\end{aligned}
\end{equation}
For $t\in[0,t_0]$, plugging \eqref{3.28} in \eqref{3.7} and applying the property that $t_0$ is sufficiently small that $2C\widetilde{M}{t_0}^{1/4}<\frac{1}{2}$, we derive that
\begin{equation*}
\begin{aligned}
\frac{1}{2}Z(t_0)  & \le C\int_0^{t_0}{\|\nabla \u_{2}(\tau)\|_{BMO}\big( 1+|\log \tau |+|\log Z(\tau )| \big)Z(\tau )\;\mathrm{d}\tau } \\
 & \le C\int_0^{t_0}{\big( 1+|\log \tau | \big)\|\nabla \u_{2}(\tau )\|_{BMO}\big( 1+|\log Z(\tau )| \big)Z(\tau )\;\mathrm{d}\tau }, \\
\end{aligned}
\end{equation*}
then, Osgood's lemma \cite{HBJ} (Lemma 3.4) and Corollary \ref{2-7} imply that $Z \equiv0$ (i.e., $\delta \rho\equiv 0$) on $[0, t_0 ]$, and thus, owing to \eqref{3.27} and $\eqref{3.28}$, that $\sqrt{\rho_1}\delta \u\equiv 0$ and $\delta\d\equiv 0$ on $[0, t_0 ]$. 
\medskip

To see $\delta\u\equiv 0$,  observe that after using $\delta \rho\equiv 0$, $\sqrt{\rho_1}\delta \u\equiv 0$ and $\delta\d\equiv 0$ in \eqref{3.3},  one has
\begin{equation*}
\rho_1 \delta \u_{t}+\rho_1 \u_1\cdot \nabla \delta \u-\mu \Delta \delta \u-(\mu +\lambda )\nabla \Dv \delta \u=0,     
\end{equation*}
which further leads to 
\begin{equation*}
 \frac{1}{2}\|\sqrt{\rho_1}\delta \u\|_{L^{\infty}(0,t_0;L^{2})}^2+\int_0^{t_0}\big(\mu\|\nabla\mathbb{P} \delta\u\|_2^{2}+\nu\|\Dv\delta\u\|_2^{2} \big)\;\mathrm{d}t=0.   
\end{equation*}
In view of \eqref{2.1.3*} together with $\int\rho_1\delta\u\,\mathrm{d}x=0$, one can obtain $\|\delta\u\|_2\leq \rho^*\|\nabla\delta\u\|_2=0$, i.e.,    $\delta\u=0$.  This
completes the proof of Proposition \ref{2-6}.
\end{proof}


\section*{Declarations}

\vskip .1in
\noindent{\bf Competing interests  }\

On behalf of all authors, the corresponding author states that there are no potential conflicts of interest with respect to the research of this article.

\vskip .1in
\noindent{\bf Authors' contributions   }\

Qiao Liu, Qingkui Ma and  Fuyi Xu  contributed equally to this work.

\vskip .1in
\noindent{\bf Funding}\

Q. Liu and Q. Ma are partially supported by the National Natural Science Foundation of China (12571252,12071122), the Natural Science Foundation of Hunan Province (2023JJ10059)  and the Fundamental Research Funds for the Central
Universities of Central South University (2026ZZTS0230).  F. Xu is partially supported by the
National Natural Science Foundation of China (12326430), and the  Natural Science Foundation of Shandong Province (ZR2026MS0022).

\vskip .1in
\noindent{\bf Data availability statement}

All data generated or analysed during this study are included in this article.

\end{document}